\documentclass[12pt]{amsart}

\usepackage{mathtools, amsthm, amsfonts, enumerate, comment, fancyhdr, mathrsfs, tikz-cd, xparse, amscd,bm}

\usepackage[margin=3cm]{geometry}

\usepackage{tikz, adjustbox, mathtools}
\usepackage[hidelinks]{hyperref}
\usetikzlibrary{matrix,calc}

\usepackage{cite}

\renewcommand{\MR}[1]{}

\newcommand{\owedge}{\mathrel{\begin{picture}(8,8)\put(4,2){\circle{8}}\put(0.3,-0.6){\scalebox{0.92}{$\wedge$}}\end{picture}}}
\newcommand{\loplus}{\mathrel{\begin{picture}(8,8)
\qbezier(4.03,8)(2.34,8)(1.17,6.83)
\qbezier(1.17,6.83)(0,5.66)(0,4)
\qbezier(0,4)(0,2.34)(1.17,1.17)
\qbezier(1.17,1.17)(2.34,0)(4.03,0)
\put(-0.55,1){$+$}\end{picture}}}

\newcounter{teller}

\newenvironment{tabel}{\begin{list}
{\rm  (\alph{teller})\hfill}{\usecounter{teller} \leftmargin=1.1cm
\labelwidth=1.1cm \labelsep=0cm \parsep=0cm}
                      }{\end{list}}

\newcounter{tellerr}

\newcounter{tellerrr}

\newcounter{thms} \numberwithin{thms}{section}
\numberwithin{equation}{section}

\newtheorem{theorem}[thms]{Theorem}
\newtheorem{proposition}[thms]{Proposition}
\newtheorem{lemma}[thms]{Lemma}
\newtheorem{corollary}[thms]{Corollary}
\newtheorem{remark}[thms]{Remark}

\theoremstyle{remark}

\theoremstyle{definition}

\newtheorem{definition}[thms]{Definition}

\allowdisplaybreaks 

\author{Simon C. W. Goodwin}
\address[Goodwin]{Department of Mathematics, The University of Auckland, Auckland 1142, New Zealand}
\email{simon.goodwin@auckland.ac.nz}
\author{A. Rod Gover} 
\address[Gover]{Department of Mathematics, The University of Auckland, Auckland 1142, New Zealand}
\email{r.gover@auckland.ac.nz}
\author{Jonathan M. Kress}
\address[Kress]{School of Mathematics and Statistics, University of New South Wales, Sydney 2052, Australia}
\email{j.kress@unsw.edu.au}
\begin{document}
\title{Superintegrable Systems and the Tractor Connection}
\begin{abstract}
Quadratic constants of motion for natural pseudo-Riemannian Hamiltonians with potential correspond to conformal Killing tensors admitting multiple Killing scales. Taking this correspondence further, we develop a dimension-independent ($n\ge 3$) framework for studying second-order superintegrable systems using conformal tractor calculus, and invariant BGG prolongation of the rank two conformal Killing equation. We apply this framework to classify second-order maximally superintegrable systems within certain natural families of Hamiltonians generated by pairs of conformally related Einstein metrics in arbitrary dimensions. We recover the second-order analogue of Bertrand’s theorem as a special case.
\end{abstract}

\thanks{
This work was supported in part by the Marsden Fund of the Royal Society Te Apārangi (grants 19-UOA-008 and 24-UOA-005) and by the Australian Research Council (Discovery Project DP190102360). This work was also supported in part by the Simons Foundation (grant SFI-MPS-T-Institutes-00010825) and by the State Treasury funds under the project ``Organization of the Simons Semesters at the Banach Center -- New Energies in 2026--2028'' (MNiSW/2025/DAP/491). Part of this work was carried out during visits by the first and second named authors to UNSW Sydney and by the third named author to the University of Auckland; they thank their respective hosts for their hospitality.}
\maketitle

\section{Introduction}

Symmetry is fundamental to formulating and understanding our physical models. Many of the models which have been most influential to our understanding of mechanics often possess an abundance of conserved quantities. The Kepler two body problem and its quantisation, the hydrogen atom, have an additional ``hidden'' symmetry, encoded in the Laplace-Runge-Lenz vector. For classical Hamiltonian systems, hidden symmetries correspond to higher order Killing tensors, which are also known to be required for separability of the Hamilton-Jacobi and Schr\"odinger equations \cite{kalnins_killing_1980}.

A \textit{Hamiltonian} is a choice of a smooth function $\mathcal{H}\colon N\to \mathbb{R} $ defined on a symplectic manifold $(N,\omega)$. Hamilton's equations are then determined by the \textit{Hamiltonian vector field} $X_\mathcal{H} = \omega^{-1}(d \mathcal{H})$, whereby a solution of Hamilton's equations is an integral curve $\gamma\colon \mathbb{R}\to N $ such that $\dot{\gamma}=X_\mathcal{H}\circ \gamma$.

In the context of classical mechanics on a manifold $M$, Hamiltonians are usually defined on $N = T^*M$, equipped with the canonical symplectic form $\omega$ and Poisson bracket $\left\{,\right\}$. Hamilton's equations then become the familiar
\begin{equation}\label{Eqn:HsE}
	\frac{dq}{dt} = \frac{\partial\mathcal{H}}{\partial p}\quad \text{and} \quad \frac{dp}{dt} = -\frac{\partial\mathcal{H}}{\partial q},
\end{equation}  
where $q = (q_{1},\dots ,q_n)$ are coordinates on $M$ and $(q_{1},\dots ,q_n,p_{1},\dots ,p_n)$ are the induced coordinates on $T^*M$. We will restrict to the classical set of Hamiltonians of the form 
\[
\mathcal{H}(q,p) = g_q^{-1}(p,p)+V(q),
\] 
for $q\in M$, $p\in T_q^*M$, where $g$ is a pseudo-Riemannian metric on $M$. A function $f\colon T^*M\to \mathbb{R}$ is a \textit{constant of motion}, or \textit{first integral} if 
\[
\omega^{-1}(d \mathcal{H},df) = \langle X_\mathcal{H}, df \rangle = \left\{\mathcal{H},f\right\} = 0,
\]
or equivalently, $f$ is constant along the trajectory of any solution of Hamilton's equations. A Hamiltonian is \textit{superintegrable} if for some $m\ge 1$, there are $n+m$ constants of motion $f_{1},\dots ,f_{n+m}$ which are \textit{functionally independent}, meaning that $df_{1},\dots ,df_{n+m}$ are linearly independent on an open dense subset of $T^*M$. The Hamiltonian $\mathcal{H}$ is \textit{maximally superintegrable} if it has $2n-1$ functionally independent constants of motion - this is the maximum number since $\langle df, X{}_{\mathcal{H}} \rangle=0$ for all constants of motion $f$. A Hamiltonian $\mathcal{H}$ is \textit{Liouville integrable} if it has $n$ functionally independent constants of motion in involution, that is, $\left\{f_i ,f{}_{j} \right\}=0$. However we follow the assumptions in \cite{miller_classical_2013}, which do not require any subset of the constants of motion in a superintegrable system to be in involution. Each constant of motion restricts solutions of \eqref{Eqn:HsE} to a hypersurface in $T^*M$. So $2n-1$ functionally independent constants of motion restrict trajectories to a one dimensional embedded submanifold, meaning each compact connected regular common level is a periodic orbit.

The prototypical examples of maximally superintegrable Hamiltonian systems are the isotropic harmonic oscillator, 
\begin{equation}\label{Eqn:I:HO}
	\mathcal{H}^{HO}(q,p) = |p|^2+\lambda r^2+C,
\end{equation} 
and the Kepler-Coulomb System
\begin{equation}\label{Eqn:I:Kep}
	\mathcal{H}^{KC}(q,p) = |p|^2+\frac{\lambda}{r}+C,
\end{equation}
both defined on the cotangent bundle of ${\mathbb{R}}{^{n}}$, with the flat Euclidean metric $g^{-1}_q(p,p) = |p|^2$. Both Hamiltonians are radially symmetric, that is of the form
\begin{equation}\label{Eqn:RSP}
\mathcal{H}(q,p) = |p|^2+V(r).
\end{equation} 
For any Hamiltonian of the form \eqref{Eqn:RSP}, the functions representing angular momenta,
\[
	L^{i,j}(q,p) = q_ip_j-q_jp_i,\quad i,j\in \{1,\dots ,n\} 
\] 
are constants of motion. However, at most $2n-3$ of the $L^{i,j}$ are functionally independent (See Proposition \ref{Prp:Find}). So together with the Hamiltonian itself, this provides $2n-2$ functionally independent constants of motion for \eqref{Eqn:RSP}. A Hamiltonian with $2n-2$ functionally independent integrals is called \textit{quasimaximally superintegrable}. The functions $L^{i,j}$ are linear in momenta, and so are induced by Killing vectors for the flat Euclidean metric, see Section \ref{Sec:Poly}. The Hamiltonians \eqref{Eqn:I:HO} and \eqref{Eqn:I:Kep} have additional functionally independent constants of motion given by the functions
\[
	f_j^{HO}(q,p) = p_j^2+\lambda q_j^2\quad \text{and} \quad f_j^{KC}(q,p) = \sum_{i=1}^{n} ((q_ip_j-q_jp_i)p_i)-\frac{\lambda q_j}{2r},
\] 
respectively, so they have $2n-1$ functionally independent constants of motion, making them maximally superintegrable. These constants of motion are quadratic in momenta, being induced by rank two Killing tensors, or ``hidden symmetries''. A Hamiltonian is second-order maximally superintegrable if it has $2n-1$ functionally independent constants of motion of degree at most two in the momenta. Superintegrable systems with constants quadratic in the momenta are particularly important, see the series of papers \cite{kalnins_second-order_2005,kalnins_second_2005,kalnins_second_2005-1,kalnins_second_2006,kalnins_second-order_2006} for a full overview of second-order superintegrability in dimensions 2 and 3. The focus on quadratic constants of motion is largely due to their relation to separability of the Hamilton-Jacobi equations \cite{kalnins_separation_2018,boyer_stackel-equivalent_1986}.

It is a classical result of Bertrand \cite{bertrand_theoreme_1873} that the only Hamiltonians of the form \eqref{Eqn:RSP} on $\mathbb{R}^3$ for which all bounded orbits are closed are \eqref{Eqn:I:HO} and \eqref{Eqn:I:Kep}. These are the only second-order maximally superintegrable Hamiltonians of the form \eqref{Eqn:RSP} on $\mathbb R^3$. This result suggests the possibility for a classification of maximally superintegrable systems. 

In 1990 Evans \cite{evans_superintegrability_1990} classified all classical Hamiltonians defined on $\mathbb{R}{}^{3}$ with the Euclidean metric that are second-order maximally superintegrable. Evans' paper was the conclusion to a program to classify all potentials which are multiseparable in three dimensions, initiated in 1967 \cite{makarov_systematic_1967}, which followed a classification in two dimensions,\cite{fris_higher_1965, winternitz_symmetry_1967}. In \cite{evans_superintegrability_1990}, Evans proved that in three dimensions quadratic superintegrability and multiseparability are equivalent, allowing him to classify maximally superintegrable systems with quadratic constants of motion in three dimensional Euclidean space.

Recently there has been significant progress towards a classification of second-order maximally superintegrable systems in dimensions $n\ge 4$ \cite{kress_algebraic_2023,kress_algebraic_2024}. In this paper we prove existence and uniqueness of second-order maximally superintegrable systems within a general class of Hamiltonians on conformally flat manifolds of arbitrary dimension, proving a natural generalisation of Bertrand's theorem to a much larger class of potentials parametrised by Einstein rescaling factors. To prove these results, we reformulate the central notions of superintegrable systems, constants of motion and functional independence, using methods derived from conformal geometry and tractor calculus. We believe that this may enable a new direction in the classification of second-order maximally superintegrable systems in arbitrary dimension.

In $\mathbb{R}{}^{n}$, the usual Euclidean metric $g_e$ is conformally related to another flat metric, given by $r^{-4}g_e$. In our analysis, the important feature of Hamiltonians of the form \eqref{Eqn:RSP} is that the potential is a function of $r{}^{2} $, which is an \textit{Einstein rescaling factor}, meaning it is the conformal factor relating the base Einstein metric $g$ to another Einstein metric. Taking this more general point of view, we prove several results. 

We prove quasimaximal superintegrability, analogous to that of \eqref{Eqn:RSP}.

\begin{theorem}[Theorem \ref{Thm:QMS}]\label{Thm:Main1}
		Let $(M,\boldsymbol{c})$ be conformally flat and simply connected. Let ${I}{}_{A} = D_A \sigma$ and $J_A = D_A \tau$ for $\sigma,\tau\in \Gamma(\mathcal{E}[1])$ be linearly independent parallel tractors. Let $V\colon \mathbb{R}\to \mathbb{R} $ be smooth and consider the Hamiltonian
	\begin{equation}\label{Eqn:Main1}
	\mathcal{H}(q,p) = g^{-1}_q(p,p)+V\left( \frac{\tau}{\sigma} \right) ,
	\end{equation} 
	where $g = \sigma^{-2}\boldsymbol{g}$. Then $\mathcal{H}$ has $2n-2$ functionally independent constants of motion, $f_{0} = \mathcal{H}$, and $f_{1},\dots ,f_{2n-3}$ induced by vectors fields which are Killing for $g$ and $\hat{g} = \tau^{-2}\boldsymbol{g}$.
\end{theorem}

In Section \ref{Sec:Ex}, we prove existence of second-order maximally superintegrable potentials in the family \eqref{Eqn:Main1}, and in Section \ref{Sec:Uni}, we prove uniqueness of these systems, completing our generalisation of Bertrand's theorem.

\begin{theorem}[Corollary \ref{Cor:U}]\label{Thm:MainU}
	Let $(M,\boldsymbol{c})$ be a conformally flat, simply connected manifold of dimension $n\ge 3$. Let $\sigma,\tau\in \Gamma(\mathcal{E}[1])$ be almost Einstein scales and let ${I}{}_{A}$, ${J}{}_{A}$ be the corresponding parallel tractors. Let $s = \tau\sigma^{-1}$ and let $\tilde{\sigma} = f(s)\sigma$ for some smooth $f\colon \mathbb{R}\to \mathbb{R} $ such that $\tilde{\sigma}$ is non-zero on an open dense subset of $M$. Suppose that the subspace spanned by $I$ and $J$ is non-degenerate. Then, without loss of generality, we consider two cases (see Theorem \ref{Thm:Ex1}). 

	Suppose that $I$ and $J$ are null and that ${I}{}_{A}{J}{}^{A}=1$. Then the Hamiltonian
	\begin{equation}
		\mathcal{H}(q,p) = g^{-1}_q(p,p)+\frac{1}{f(s)^2} = g^{-1}_q(p,p)+V(s),
	\end{equation}
	where $g = \sigma^{-2}\boldsymbol{g}$ and $V(s) = f(s)^{-2}$, is second-order maximally superintegrable if and only if $f$ solves either
	\[
		v(s) = \frac{f(s)^4}{2}V''(s) = 0, \quad \text{in which case}\quad V(s) = as+b
	\]
	or
	\[
		u(s)-sv(s) = -\frac{f(s)^4}{2}(2sV''(s)+3 V'(s)) = 0, \quad \text{in which case}\quad V(s) = \frac{a}{\sqrt{|s|} }+b
	\] 

	Let $\iota=I_AI^A$ and $\jmath  = J_AJ^A$. Suppose that $I$ and $J$ are orthogonal and $\iota,\jmath \in \{-1,+1\} $. Then the Hamiltonian
	\begin{equation}
		\mathcal{H}(q,p) = g^{-1}_q(p,p)+\frac{1}{f(s)^2} = g^{-1}_q(p,p)+V(s),
	\end{equation}
	is second-order maximally superintegrable if and only if $f$ solves either
	\[
		u(s) = -\frac{f(s)^4}{2}(sV''(s)+3V'(s)) = 0, \quad \text{in which case}\quad V(s) = \frac{a}{s^2}+b
	\] 
	or
	\[
		v(s)-\iota\jmath su(s)=\frac{f(s)^4}{2}((\iota\jmath s^2+1)V''(s)+3\iota\jmath sV'(s))=0,\quad \text{in which case}\quad V(s) = a \frac{s}{\sqrt{|\iota\jmath s^2+1|} }+b.
	\] 
\end{theorem}

Our approach to constants of motion is through conformal Killing tensors. While classically superintegrable systems are primarily understood using the Poisson structure, our shift towards conformal Killing tensors allows us to apply the invariant tractor calculus, and related constructions to superintegrable systems.

The \textit{Killing equation} for rank $\ell$ symmetric covectors on a Riemannian manifold is
\begin{equation}\label{Eqn:Kill}
{{\nabla}^g}{}_{\mspace{-6mu}(a}{r}{}_{b_{1}\dots b_\ell)} = 0,
\end{equation}
where $\nabla^g$ is the Levi-Civita connection of $g$. The leading term of any polynomial constant of motion must come from a Killing tensor for the metric $g$ (see Section \ref{Sec:CM}). When the Hamiltonian has nonzero potential, $V$, the full requirement on a constant of motion can be reformulated in the language of conformal geometry. Quadratic and linear constants of motion for a Hamiltonian of the form $\mathcal{H} = g^{-1}+V$ are in one to one correspondence with conformal Killing vectors and rank two conformal Killing tensors for which the conformally related metrics $g$ and $Vg$ are both Killing scales, see Section \ref{Sec:CKT}. 

Conformal Killing tensors of ranks one and two and their Killing scales were the subject of a recent paper, \cite{gover_conformal_nodate} (and extensions in \cite{goodwin_killing_2026}), which uses the conformally invariant tractor calculus to prolong the conformal Killing equation. In addition they derive algebraic formulae which govern whether certain metrics in the conformal class are Killing scales.

Currently these methods require $(M,\boldsymbol{c})$ to be conformally flat, as the full prolongation connection and splitting operator for second order Killing tensors is not known in the non-flat case. The conformally flat case is the setting usually considered for second-order superintegrable systems. When considering so-called ``abundant'' nondegenerate quadratically superintegrable systems, one is forced to be on a conformally flat manifold \cite[Corollary 5.10]{kress_algebraic_2023}. Therefore restricting to the conformally flat case is natural, and fits with the current framework for superintegrable systems. Certainly superintegrable systems have been constructed in conformally curved spaces \cite{kalnins_superintegrability_2013} and with further developments in tractor prolongations for higher order Killing tensors, our methods may extend to those settings also. 

Our constructions and calculations are mostly local, so we require our manifolds to have trivial conformal holonomy, which is stronger than conformally flat by Ambrose-Singer, to show global existence of quadratic constants of motion. However, when proving uniqueness, we construct PDEs that must be satisfied locally, providing obstructions to superintegrability independent of global topology. We note that if $(M,\boldsymbol{c})$ is conformally flat and simply connected, then it admits a conformal immersion (the developing map) into the M\"obius space of signature $(P,Q)$, which in the Riemannian case is the conformal sphere \cite{kuiper_conformally-flat_1949}. In particular, $(M,\boldsymbol c)$ then has trivial conformal holonomy.

The structure of our paper is as follows. In Section \ref{Sec:Bcg}, we introduce the necessary background on Hamiltonian mechanics and conformal geometry. In Section \ref{Sec:CM}, we give a full derivation of the relationship between constants of motion and conformal Killing tensors. These sections keep our presentation self contained and in the language and notation found in conformal geometry. In Section \ref{Sec:Tct}, we give the standard construction of the conformal tractor bundle, and its invariant connection and calculus. We also summarise the results we require from the sources \cite{gover_conformal_2021,goodwin_killing_2026}. In Sections \ref{Sec:FI} and \ref{Sec:Ein}, we give a formulation of functional independence in terms of tractors and prove the first main result, Theorem \ref{Thm:Main1}. In Section \ref{Sec:Ex}, we show that there exist second-order maximally superintegrable systems contained in the families derived in Theorem \ref{Thm:Main1}. In Section \ref{Sec:Uni}, we show that the second-order maximally superintegrable systems constructed in Section \ref{Sec:Ex} are unique, completing the proof of Theorem \ref{Thm:MainU}. In Section \ref{Sec:Ham}, we discuss linear families of Hamiltonians, showing how extra potential terms can be added to the Hamiltonians previously derived. Finally in Section \ref{Sec:FM}, we show how classical Hamiltonians, written in coordinates, can be extracted from our tractor construction, and we calculate two examples on $\mathbb{R}{}^{n}$ with flat and hyperbolic metrics.

The superintegrable systems we construct include the isotropic harmonic oscillator,
\[
	|p|^2+\lambda r^2+\frac{a_{1}}{q_{1}^2}+\dots + \frac{a_n}{q_n^2},
\] 
the Kepler system
\[
	|p|^2+\frac{\lambda}{r}+\frac{a_{2}}{q_{2}^2}+\dots +\frac{a_n}{q_n^2},
\] 
and a higher dimensional generalisation of Evans' third potential \cite{evans_superintegrability_1990}
\[
q_{1}^2|p|^2+\frac{\lambda q_{2}}{\sqrt{q_{2}^2+q_{1}^2} } +a_{3}\frac{q_{1}^2}{q_{3}^2}+\dots +a_n \frac{q_{1}^2}{q_n^2}.
\] 

This reveals the importance of Einstein scales in writing down superintegrable systems, as well as classifying them in higher dimensions.

\section{Background}\label{Sec:Bcg}
In this section we give an introduction to Hamiltonian mechanics, with a view towards polynomial superintegrability on pseudo Riemannian manifolds, and conformal geometry.

Throughout this article $M$ will be a connected smooth manifold of dimension $n\ge 3$.
\subsection{Symplectic Geometry on \texorpdfstring{$T^*M$}{T*M}}:w
The cotangent bundle $T^*M$ is a smooth manifold of dimension $2n$. In the context of classical mechanics, $T^*M$ is often called the \textit{phase space}, a point of $M$ is called a \textit{position} and is typically denoted $q\in M$. A cotangent vector is called a \textit{(generalised) momenta} and is typically denoted by $p\in T_{q} ^*M$. Given a chart $\phi\colon U\to \mathbb{R}^n $, we denote the local coordinates by $(q_{1},\dots ,q_n) = \phi(q)$, which induce local coordinates $(q_{1},\dots ,q_n,p_{1},\dots ,p_n)$ on $T^*M$, where $p = \sum_{i=1}^{n} p_idq_i$.

There is a canonical symplectic form, $\omega\in \Gamma(\bigwedge^2(T^*(T^*M)))$ on $T^*M$ \cite{lee_introduction_2012}, which in local coordinates is given by
\[
\omega = \sum_{i=1}^{n} dq_i\wedge dp_i. 
\] 
Let $f_{1},f_{2}$ be smooth functions defined on open subsets of $T^*M$. The canonical symplectic form $\omega$ induces a Poisson bracket, $\{f_{1},f_{2}\} := \omega^{-1}(df_{1},df_{2})$ defined on the common domain of $f_{1},f_{2}$. In local coordinates, this is given by
\begin{equation}\label{Eqn:PB}
\left\{f_{1},f_{2}\right\}= \sum_{i=1}^{n} \left( \frac{\partial f_{1}}{\partial p_i} \frac{\partial f_{2}}{\partial q_i}-\frac{\partial f_{2}}{\partial p_i} \frac{\partial f_{1}}{\partial q_i}\right).
\end{equation}

The \textit{Hamiltonian vector field} of a real valued smooth function $\mathcal{H}$ defined on an open dense subset of $T^*M$ is $X_\mathcal{H} = \omega^{-1}(d\mathcal{H})$. By construction, the bracket of a Hamiltonian and a smooth function $f$ is $\left\{\mathcal{H},f\right\} =  X_\mathcal{H}(f)$, the derivative of $f$ along the flow of $X_\mathcal{H}$.

A nominated function $\mathcal{H}$ may be called a \textit{Hamiltonian}. Then Hamiltonian mechanics studies the dynamics generated by the flow of $X_\mathcal{H}$, governed by Hamilton's equations,
\[
\frac{dq}{dt} = \frac{\partial\mathcal{H}}{\partial p}\quad \text{and} \quad \frac{dp}{dt} = -\frac{\partial\mathcal{H}}{\partial q},
\] 
the trajectories and long term behaviour. 
A \textit{constant of motion} for a Hamiltonian $\mathcal{H}$ is a function, $f$, defined on an open dense subset of $T^*M$ such that, $\left\{\mathcal{H},f\right\} = X_\mathcal{H}(f) = 0$. Constants of motion provide first integrals for solving the equations of motion given by Hamilton's equations.

\subsection{Polynomial functions}\label{Sec:Poly}
Let $v\in \Gamma(S^\ell(TM))$ be a symmetric contravariant tensor field. Define a smooth function $\tilde{v}\colon T^*M\to \mathbb{R} $ by
\begin{equation}\label{Eqn:vqp}
	\tilde{v}(q,p) = v_q(p,\dots, p) = \sum_{i_{1}\dots i_\ell\in \{1,\dots ,n\} }v^{i_{1}\cdots i_\ell}(q)\,p_{i_1}\cdots p_{i_\ell},
\end{equation} 
for each $q\in M$ and $p\in T^*_qM$. Or simply $\tilde{v}(q,p) = v(q)$ when $\ell=0$. In local coordinates on $M$, the function $\tilde{v}$ is a homogeneous polynomial of degree $\ell$ in the momenta $p_{1},\dots ,p_n$ with coefficients smooth in $q$. Given a homogeneous polynomial function, one can then define a symmetric tensor $v$ for which $\tilde{v}$ is the given polynomial. 

Hence, we have a linear isomorphism between homogeneous polynomial functions of degree $\ell$ and symmetric contravariant tensors of degree $\ell$. Which then extends to an algebra isomorphism between the symmetric algebra $S(TM)$ with the symmetric product and the algebra of polynomial functions on $T^*M$ with pointwise multiplication.

Let $f,g$ be smooth functions on $M$ and let $v,u\in \Gamma(TM)$. Suppose that in local coordinates $u = \sum_{i=1}^{n} u^i\frac{\partial}{\partial q^i}$ and $v = \sum_{i=1}^{n} v^i\frac{\partial}{\partial q^i}$. Then we have $\{\tilde{f},\tilde{g}\}=0$, 
\begin{equation}\label{Eqn:vfbr}
\{\tilde{v},\tilde{f}\} = \sum_{i=1}^{n} v^i \frac{\partial f}{\partial q^i} = \widetilde{v(f)},\quad \text{and} \quad
\end{equation}
\begin{equation}\label{Eqn:vubr}
	\{\tilde{v},\tilde{u}\} = \sum_{i=1}^{n} \sum_{j=1}^{n} \left( v^i\frac{\partial u^j}{\partial q^i} p_j - u^i\frac{\partial v^j}{\partial q^i } p_j\right) = \widetilde{[v,u]},
\end{equation} 
where $[,]$ denotes the Lie bracket of vector fields on $M$. So there is an isomorphism of Lie algebras between $\Gamma(TM)$ equipped with the Lie bracket, and homogeneous degree one polynomial functions on $T^*M$ equipped with the Poisson bracket.

In classical mechanics, it is common to allow for singularities of functions on $M$ (e.g. \eqref{Eqn:I:Kep}). To avoid repetition, we will always allow for smooth functions, including the functional coefficients of polynomial functions, and tensor fields to be defined on an open dense subset of $M$, and include this assumption in notation such as $\Gamma(TM)$ and $C^\infty(M)$.

\subsection{The symmetric Schouten concomitant}

A Poisson algebra is a triple $(A,\boldsymbol{\cdot} ,\left\{,\right\})$ such that $(A,\boldsymbol{\cdot} )$ is an associative real algebra, $(A,\left\{,\right\})$ is a Lie algebra and the following \textit{Leibniz property} holds
\[
	\left\{x\boldsymbol{\cdot} y,z\right\} = x\boldsymbol{\cdot} \left\{y,z\right\}+\left\{x,z\right\}\boldsymbol{\cdot} y
\] 

The space of polynomial functions on $T^*M$ carries a Poisson algebra structure, where the product is pointwise multiplication and the bracket is the Poisson bracket. The bracket preserves being polynomial by \eqref{Eqn:PB}.

The space $S(TM)$ also has a Poisson algebra structure, the product is given by the symmetric product, which we will write $v\vee u$, and the bracket is the (symmetric) \textit{Schouten concomitant}, characterised by
\begin{enumerate}
	\item If $f,g\in C^\infty(M)$, then $\left\{f,g\right\}=0$,
	\item If $u\in \Gamma(TM)$ and $f\in C^\infty(M)$ is a function, then $\{u,f\} = u(f) = df(u)$,
	\item If $u,v\in \Gamma(TM)$, then $\{u,v\}$ is the Lie bracket of $u$ and $v$,
	\item If $u,v,w\in \Gamma(S(TM))$, then 
		\[
			\{u,v\vee w\} = \{u,v\}\vee w + \{u,w\}\vee v.
		\] 
\end{enumerate}

If $M$ is equipped with a torsion free affine connection, then the \textit{(symmetric) Schouten concomitant} of $v\in \Gamma(S^\ell_{1}TM)$ and $w\in \Gamma(S^\ell_{2}TM)$ is given by (see \cite{nijenhuis_jacobi-type_1955})
\begin{equation}\label{Eqn:SCC}
\{v,w\}^{a_1\cdots a_\ell}= \ell_{1}v^{c(a_{\ell_{2}+1}\cdots a_\ell}\nabla_{c}w^{a_1\cdots a_{\ell_{2}})}-\ell_{2} w^{c(a_{\ell_{1}+1}\cdots a_\ell}\nabla_{c}v^{a_1\cdots a_{\ell_{1}})}\in \Gamma(S^\ell(TM)),
\end{equation}
where $\ell=\ell_{1}+\ell_{2}-1$.

\begin{proposition}\label{Prp:Br}
	The Poisson algebra $(\Gamma(S(TM)),\vee,\left\{,\right\})$, is isomorphic to the Poisson algebra of polynomial functions on $T^*M$, equipped with pointwise multiplication and the Poisson bracket. The isomorphism is given by taking the associated polynomial, as in \eqref{Eqn:vqp}. 
\end{proposition}
\begin{proof}
	The calculations in \eqref{Eqn:vfbr} and \eqref{Eqn:vubr} show that the brackets are preserved for sections of $\mathcal{E}\oplus TM$, where $\mathcal{E} = M\times\mathbb{R}$. Sections of $\mathcal{E}\oplus TM$ generate $\Gamma(S(TM))$ under the operation $\vee$ and linear combinations, so $v\mapsto \tilde{v}$ extends to a Poisson algebra isomorphism.
\end{proof}In light of this isomorphism, we will drop the notation $\tilde{v}$ and use the same symbol $v$ to refer to the tensor and the function, when there is no room for confusion. 

\subsection{Conformal Geometry}\label{Sec:CfG}

Our method for studying constants of motion and superintegrability is through conformal geometry.
Two metrics, $g$ and $\hat{g}$ on $M$ are \textit{conformally related} if $\hat{g} = \Omega^2g$ for some smooth positive function $\Omega$.  An equivalence class $\boldsymbol{c} = [g]$ of conformally related metrics is called a \textit{conformal class}, and the pair $(M,\boldsymbol{c})$ is a \textit{conformal manifold}. Our reference for the following section is \cite{curry_introduction_2018}.

Let $(M,\boldsymbol{c})$ be a conformal manifold. Write $\nabla{}^{g} $ for the Levi-Civita connection of $g\in\boldsymbol{c}$. Suppose that $g,\hat{g}\in \boldsymbol{c}$ and $\Omega>0$ is such that $\hat{g} = \Omega^2 g$. The conformal transformation laws for $\omega\in\Gamma(\mathcal{E}_a)$ and $v^a\in\Gamma(\mathcal{E}^a) $ are
\begin{equation}\label{Eqn:CT}
	\nabla_a{}^{\mspace{-11mu}\hat{g}}\omega_b = \nabla_a{}^{\mspace{-11mu}{g}}\omega_b -\Upsilon_a\omega_b-\omega_a\Upsilon_b+\Upsilon^c\omega_cg_{ab},\quad \nabla_a{}^{\mspace{-11mu}\hat{g}}v^b = \nabla_a{}^{\mspace{-11mu}{g}} v^b+\Upsilon_av^b-v_a\Upsilon^b+\Upsilon_cv^c\delta{}_{a}{}^{b},
\end{equation}
where $\Upsilon_a = \Omega^{-1}\nabla_a \Omega$. Transformations for any higher rank tensor bundle can be deduced from the above. In particular for an $n$-form, $\omega\in \Gamma(\Lambda^n T^*M)$,
\[
	\nabla_a{}^{\mspace{-11mu}\hat{g}}\omega = \nabla_a{}^{\mspace{-11mu}g}\omega -n\Upsilon_a \omega.
\] 
So we call the oriented line bundle $\mathcal{E}[-2n] = (\Lambda^n T^*M)^2$ the \textit{conformal density bundle} of weight $-2n$, and oriented roots $\mathcal{E}[w] = \mathcal{E}[-2n]^{\frac{w}{-2n}}$ the conformal density bundles of weight $w$ for any $w\in \mathbb{R}$. Since the density bundles are oriented, they have a ray subbundle of positive elements, denoted $\mathcal{E}_+[w]$. Then
\begin{equation}\label{Eqn:TfW}
	\nabla_a{}^{\mspace{-11mu}\hat{g}}\sigma = \nabla_a{}^{\mspace{-11mu}g}\sigma +w\Upsilon_a \sigma
\end{equation} 
for a density $\sigma\in \Gamma(\mathcal{E}[w])$. For any vector bundle $\mathcal{B}$, we define $\mathcal{B}[w] = \mathcal{B}\otimes\mathcal{E}[w]$. Tensoring weights changes the transformation laws. For a section $k_a\in \Gamma(\mathcal{E}_a)$,
\[
	{\nabla{}^{\hat{g}}}{}_{\mspace{-8mu}(a}k{}_{b)}={\nabla^{g}}{}_{\mspace{-8mu}(a}k{}_{b)} -2\Upsilon_{(a}k_{b)} +\Upsilon^c k_c g_{ab},
\] 
whereas for a section $k_a\in \Gamma(\mathcal{E}_a[2])$,
\[
	{\nabla{}^{\hat{g}}}{}_{\mspace{-8mu}(a}k{}_{b)}={\nabla^{g}}{}_{\mspace{-8mu}(a}k{}_{b)} +\Upsilon^c k_c g_{ab},
\] 
which shows that the \textit{conformal Killing operator}, $k_b\mapsto \nabla{}_{(a} k{}_{b)_0} $ is \textit{conformally invariant}, that is it does not depend on the metric $g\in \boldsymbol{c}$, so we omit the $g$ from the connection.

Any metric $g\in \boldsymbol{c}$ determines the (pseudo) Riemannian density $\sigma^g\in \Gamma(\mathcal{E}[-2n])$ (see \cite[Proposition 2.44]{lee_introduction_2018}), which in an orientable neighbourhood is the square of the volume form. 
Then we call $\sigma = (\sigma^g)^{-1 / 2n}\in  \mathcal{E}_+[1]$ the \textit{scale density} corresponding to $g\in \boldsymbol{c}$. The scale density $\sigma$ corresponding to $g$ is parallel for the Levi-Civita connection of $g$. 

If $\hat{g} = \Omega^2 g$, then the scale densities are related by $\Omega\hat{\sigma} = \sigma$. So it follows that the \textit{conformal metric} $\boldsymbol{g} := \sigma^2 g=\hat{\sigma}^2\hat{g}$ is conformally invariant. Moreover,
\[
	\nabla_a^g\boldsymbol{g} = (\nabla^g_a\sigma)g+(\nabla^g_ag)\sigma=0
\]
for any $g\in \boldsymbol{c}$. So we will henceforth use $\boldsymbol{g}_{ab}\in \mathcal{E}_{(ab)}[2]$ and $\boldsymbol{g}^{ab}\in \Gamma(\mathcal{E}^{(ab)}[-2])$ to raise and lower indices. 

A metric $g\in \boldsymbol{c}$ corresponds bijectively to nowhere vanishing scale densities $\sigma\in \mathcal{E}_+[1]$, and are related by $g = \sigma^{-2}\boldsymbol{g}$. We will also consider scale densities with zeros $\sigma\in \mathcal{E}[1]$, with the assumption that $M\setminus\mathcal{Z}(\sigma)$ is dense in $M$. In this case, the metric $g = \sigma^{-2}\boldsymbol{g}$ is defined on $M\setminus\mathcal{Z}(\sigma)$. Such a density is called a \textit{scale}.

The conformal Killing operator is defined for $k_{a_1\cdots a_\ell}\in\Gamma(\mathcal{E}_{(a_1\dots a_\ell)_0}[2\ell])$ by
\begin{equation}\label{Eqn:CKE}
	\nabla{}_{(a_0}k{}_{a_1\dots a_\ell)_0}=0
\end{equation}
and is conformally invariant. A solution of \eqref{Eqn:CKE} is called a \textit{conformal Killing tensor}. Rank one conformal Killing tensors are called \textit{conformal Killing vectors}. Conformal Killing vectors are infinitesimal conformal transformations, much like Killing vectors in Riemannian geometry. 

If $r_{a_1\cdots a_\ell}\in\Gamma(\mathcal{E}_{(a_1\dots a_\ell)}[2\ell])$, then raising every index using the conformal metric gives an unweighted $r^{a_1\cdots a_\ell}\in\Gamma(\mathcal{E}^{(a_1\dots a_\ell)})$, which then corresponds to a polynomial function on $T^*M$, see Section \ref{Sec:Poly}. Where there is no room for confusion, we will use the same symbol for the weighted tensor, unweighted tensor and polynomial function on $T^*M$.

Let $r\in \Gamma(S^\ell(TM))$. Then $r$ is a \textit{Killing tensor} for a metric $g\in \boldsymbol{c}$ if and only if
\[
	\nabla{}^{g}{}_{(a} r_{a_{1}\dots a_\ell)}=0.
\] 
If $r$ is a Killing tensor for $g\in \boldsymbol{c}$, then the trace-free part of $r$,
\begin{equation*}
	{k}{}_{a_{1}\dots a_\ell} = {r}{}_{(a_{1}\dots a_\ell)_0}
\end{equation*}
is a conformal Killing tensor. It is natural to ask when a conformal Killing tensor is the trace free part of a Killing tensor for some metric in $\boldsymbol{c}$. The paper \cite{gover_conformal_nodate} studies this question for rank two conformal Killing tensors. Following that paper we make the following definition.
\begin{definition}
Let ${k}{}_{a_{1}\dots a_\ell}\in \Gamma(\mathcal{E}_{(a_{1}\dots a_\ell)_0}[2\ell])$ ($\ell\ge 2$) be a conformal Killing tensor of rank $\ell$, let $\sigma\in \Gamma(\mathcal{E}[1])$ be a scale and let $g = \sigma^{-2}\boldsymbol{g}$ be the metric defined on $M\setminus\mathcal{Z}(\sigma)$.
	\begin{tabel}
		\item\label{Def:KS:1} $\sigma$ is a \textit{strong Killing scale} for $k$ if $k$ is a Killing tensor for $g$,
		\item\label{Def:KS:2} $\sigma$ is a \textit{Killing scale} for $k$ if there exists a Killing tensor $r$ for $g$ such that $k$ is the trace-free part of $r$.
		\item\label{Def:KS:3} Let $k_a\in \Gamma(\mathcal{E}_a[2])$ be a conformal Killing vector. $\sigma$ is a \textit{Killing scale} for $k$ if $k$ is a Killing vector for $g$.
	\end{tabel}
	
\end{definition}

In this article we will primarily be concerned with rank one conformal Killing vectors and rank two conformal Killing tensors.

\section{Constants of Motion and Conformal Geometry}\label{Sec:CM}
\subsection{Constants of Motion}
Let $g$ be a pseudo-Riemannian metric on $M$ with inverse $g^{-1}\in \Gamma(S^2(TM))$. We will consider Hamiltonians of the form
\begin{equation}\label{Eqn:Ham}
\mathcal{H}(q,p) = g^{-1}_q(p,p)+V(q),
\end{equation} 
	where $V\colon M\setminus S\to \mathbb{R}$ is a smooth function, called the \textit{potential}. The singularity set $S$ must be contained in a finite union of properly embedded submanifolds of dimension less than $n$.

A function $f\colon T^*M\to \mathbb{R} $ is a \textit{constant of motion}, or \textit{first integral} for $\mathcal{H}$ if $\{\mathcal{H},f\}=0$. The order of a polynomial constant of motion is its polynomial degree. \textit{Linear} and \textit{quadratic} constants of motion refer to order one and order two constants of motion. 

The space of polynomial constants of motion is a Poisson subalgebra of $S(TM)$, which follows from the Leibniz property and the Jacobi identity. 

Suppose that the potential in $\mathcal{H}$ is zero, so $\mathcal{H} = g^{-1}$. Then the trajectories of $X_\mathcal{H}$ are geodesics of the base manifold \cite[\S 45.D]{arnold_mathematical_1989}. Let $r\in \Gamma(S^\ell(TM))$, then we have
\begin{equation}\label{Eqn:KILL} 
	\{g^{-1},r\} = 2\nabla^{(a}r^{a_{1}\dots a_\ell)},
\end{equation} 
by \eqref{Eqn:SCC}, where $\nabla$ is the Levi-Civita connection of $g$. This is the Killing operator on sections of $S^\ell(TM)$.
\begin{proposition}\label{Prp:Kill}
	Let $g$ be a pseudo-Riemannian metric on $M$. If $r\in \Gamma(S^\ell(TM))$, then
	\[
	\{g^{-1},r\}=0
\] if and only if $r$ is a rank $\ell$ Killing tensor for $g$.
\end{proposition}

When $\mathcal{H}$ has a non-zero potential, there is an additional term in the bracket,
\begin{equation}\label{Eqn:Hr}
\{\mathcal{H},r\} = \{g^{-1},r\}+\{V,r\}.
\end{equation} 

The terms on the right are homogeneous of orders $\ell+1$ and $\ell-1$, respectively. So $\{\mathcal{H},r\}=0$ if and only if both $\{g^{-1},r\}$ and $\{V,r\}$ are zero. The first requirement is that $r$ is a Killing tensor for $g$ and the second condition can be transformed into the Killing equation for $Vg$. Hence, we have
\begin{proposition}\label{Prp:PBKT}
	Let $(M,\boldsymbol{c})$ be a conformal manifold and let $\mathcal{H}(q,p) = g^{-1}+V$ for some $g\in \boldsymbol{c}$ and $V$ a potential which is nonzero on an open dense set. If $\ell\in \mathbb{N}$ and $r\in \Gamma(S^\ell(TM))$, then
	\begin{equation}\label{Eqn:Kill2}
	\left\{\mathcal{H},r\right\}=0
	\end{equation} 
	if and only if $r$ is a rank $\ell$ Killing tensor for $g$ and the metric $Vg$
\end{proposition}
\begin{proof}
	Where $V$ is defined and nonzero we have
	\begin{equation}\label{Eqn:Prp:CM1}
		\left\{g^{-1},r\right\}=\{VV^{-1}g^{-1},r\} = V\{(Vg)^{-1},r\}+(Vg)^{-1}\{V,r\}.
	\end{equation} 
by the Leibniz rule for $\left\{,\right\}$. So we have the following equivalences,
\[
	\{g^{-1}+V,r\} = 0 \iff \left\{g^{-1},r\right\} = 0 = \left\{V,r\right\} \iff  \left\{g^{-1},r\right\} = 0 = \left\{(Vg)^{-1},r\right\}.
\] 
By Proposition \ref{Prp:Kill}, the last equality is equivalent to $r$ being Killing for $g$ and $Vg$. \end{proof}

\subsection{Quadratic Constants}

For $\ell = 2$, restricting to constants of motion of the form \eqref{Eqn:Kill2} is too restrictive since it requires $r$ to be a Killing tensor for both $g$ and $Vg$. A natural weakening of this condition is to ask that $Vg$ is a Killing scale for the conformal Killing tensor $k{}_{ab} = r{}_{(ab)_0} $, whence the actual Killing tensor for $Vg$ may differ from $r$ by a trace term.In this section we relate conformal Killing tensors and Killing scales to second order constants of motion for a Hamiltonian of the form $\mathcal{H} = g^{-1}+V$. Note that a section ${k}{}_{ab}\in \Gamma({\mathcal{E}}{}_{(ab)}[4])$ has the correct weight to give an unweighted tensor field after raising indices, ${k}{}^{ab}\in \Gamma(S^2(TM))$. Therefore we will raise indices and use Proposition \ref{Prp:Br} to relate ${k}{}_{ab}\in \Gamma({\mathcal{E}}{}_{(ab)}[4])$ to a second order function on $T^*M$.

A second order constant of motion is of the form $f(q,p) = r^2_q(p,p) + r^1_q(p)+U(q)$. Expanding $\left\{\mathcal{H},f\right\}=0$ and grouping terms by order gives
\[
\left\{\mathcal{H},f\right\} = \left\{g^{-1},r^2\right\}+\left\{g^{-1},r^1\right\}+(\left\{V,r^2 \right\}+\left\{g^{-1} ,U\right\})+\left\{V,r^1\right\}=0
\] 
Summing the order two and order zero terms gives $\left\{\mathcal{H},r^1\right\}=\left\{g^{-1},r^1\right\}+\left\{V,r^1\right\}=0$, so $r^1$ is a linear constant of motion. So we consider $r^1$ and $r^2+U$ independently. Therefore we use the term \textit{quadratic constant of motion} to refer to constants of motion of the form $f(q,p) = r_q(p,p) +U(q)$

\begin{theorem}\label{Thm:QCns}
	Let $(M,\boldsymbol{c})$ be a conformal manifold and let $\mathcal{H}(q,p) = g_q^{-1}(p,p)+V(q)$ for some $g\in \boldsymbol{c}$ and $V\in C^\infty(M)$, such that $V$ is non-zero on an open dense set. If $k_{ab}\in \Gamma(\mathcal{E}_{(ab)_0}[4])$, then there exist $U,T\in C^\infty(M)$ such that 
	\begin{equation}\label{Eqn:SOC}
		\{\mathcal{H},k+Tg^{-1}+U\}=0
	,\end{equation} 
	if and only if $k$ is a conformal Killing tensor and both $g$ and $\hat{g}=Vg$ are Killing scales for $k$. 

	If \eqref{Eqn:SOC} holds, then $r{}_{ab}  = k{}_{ab} +T\sigma^{4}g{}_{ab} = k{}_{ab} +T\sigma^2\boldsymbol{g}_{ab}$ is a Killing tensor for $g$ and
	\[
		\hat{r}{}_{ab}  = k{}_{ab} +(VT-U)\hat{\sigma}^{4}\hat{g}{}_{ab} = k{}_{ab} +(VT-U)\hat{\sigma}^{2}\boldsymbol{g}_{ab} = k{}_{ab} +\frac{VT-U}{V}\sigma^2\boldsymbol{g}_{ab}=r{}_{ab} -\frac{U}{V}\sigma^{4}g{}_{ab}
	\]
	is a Killing tensor for $\hat{g}$.
\end{theorem}
\begin{proof}
	We first expand the bracket,
	\[
		\{\mathcal{H},k+Tg^{-1}+U\} = \{g^{-1}+V,k+Tg^{-1}+U\} = \{g^{-1},r\}+\{V,r\}+\{g^{-1},U\}
	,\]
	where ${r} = {k} +T{g}^{-1} $. Taking terms of each order gives
	\begin{equation}\label{Eqn:QCns:1}
	\{\mathcal{H},k+Tg^{-1}+U\} = 0\iff \{g^{-1},r\}=0 = \{V,r\}+\{g^{-1},U\}
	\end{equation}
	We then calculate
	\begin{align*}
		\{g^{-1},r\}
		&=\{V(Vg)^{-1},r\} \\
		&= V\{\hat{g}^{-1},r\}+\hat{g}^{-1}\{V,r\}\\
		&= V(\{\hat{g}^{-1},r\}-\{\hat{g}^{-1},U \hat{g}^{-1}\})+\hat{g}^{-1}\{V,r\}+V\{\hat{g}^{-1},U \hat{g}^{-1}\}\\
		&= V(\{\hat{g}^{-1},r-U \hat{g}^{-1}\})+\hat{g}^{-1}\{V,r\}+V\hat{g}^{-1}\{V^{-1}{g}^{-1},U\}\\
		&= V(\{\hat{g}^{-1},r-U \hat{g}^{-1}\})+\hat{g}^{-1}\left(\{V,r\}+\{{g}^{-1},U\}\right).\\
	\end{align*}So on the open dense set where $\hat{g}$ is non-singular, and $V$ is non-zero, we have that
\begin{equation}\label{Eqn:QCns:2}
\{g^{-1},r\}=0 = \{V,r\}+\{g^{-1},U\}\iff \{g^{-1},r\}=0=\left\{\hat{g}^{-1},r-U \hat{g}^{-1}\right\}
\end{equation} 
Therefore the existence of $T$ and $U$ imply the existence of Killing tensors $r = k+Tg^{-1}$ for $g$ and $\hat{r} = r-U \hat{g}^{-1} = k+(VT-U)\hat{g}^{-1}$ for $\hat{g}$, so both $g$ and $\hat{g}$ are Killing scales.

Conversely, if $g$ and $\hat{g}$ are Killing scales for $k$, with corresponding Killing tensors $r = k+Tg^{-1}$ and $\hat{r} = k+\hat{T}\hat{g}^{-1}$, then taking $U = VT-\hat{T}$ gives
\[
\left\{g^{-1},r\right\}=0 = \left\{\hat{g}^{-1},k+\hat{T}\hat{g}^{-1}\right\} = \left\{\hat{g}^{-1},r-U\hat{g}^{-1}\right\}
\] 
So $\left\{\mathcal{H},k+Tg^{-1}+U\right\} = 0$ by \eqref{Eqn:QCns:1} and \eqref{Eqn:QCns:2}.
\end{proof}

At a fixed energy, $\mathcal{H} = g^{-1}+V = E$, the trajectories of $\mathcal{H}$ are, up to reparametrisation, geodesics of the Jacobi metric $g_J = (E-V)g$. So the metric $Vg$ could be called the Jacobi metric with zero energy (up to sign). This move from $\mathcal{H}$ to the Hamiltonian $\hat{\mathcal{H}} = (Vg)^{-1}$ is called the \textit{St\"ackel transform} \cite{boyer_stackel-equivalent_1986}.
\begin{corollary}\label{Cor:QCns}
	Let $(M,\boldsymbol{c})$ be a conformal manifold, let $\sigma\in \Gamma(\mathcal{E}[1])$ be a scale and let $V$ be a potential and let $\mathcal{H} = g^{-1}+V$, where $g = \sigma^{-2}\boldsymbol{g}$. Then there is an $\mathbb{R}$-linear isomorphism between the space of quadratic constants of motion for $\mathcal{H}$, modulo constant functions and constant multiples of $\mathcal{H}$, and rank two conformal Killing tensors for which $g$ and $Vg$ are both Killing scales.
\end{corollary}
\begin{proof}
	Let $f(q,p) = r_q(p,p)+U(q)$ be a constant of motion for $\mathcal{H}$. The order three term in $\left\{\mathcal{H},f\right\}$ is $\left\{g^{-1},r\right\}$, so $r$ is a Killing tensor for $g$. Then ${k}{}_{ab} = {r}{}_{(ab)_0}$ is a conformal Killing tensor and there exist functions $T,U$ such that $k+Tg^{-1}+U$ is a constant of motion for $\mathcal{H}$. We call this linear map from constants of motion to conformal Killing tensors $\Phi$.

	In the proof of Theorem \ref{Thm:QCns}, we showed that if $r = k+Tg^{-1}$ and $\hat{r} = k+\hat{T}\hat{g}^{-1}$ are Killing for $g$ and $\hat{g}$, then $k+Tg^{-1}+(VT-\hat{T})$ is a quadratic constant of motion for $\mathcal{H}$. This shows that $\Phi$ is surjective.

	The kernel of $\Phi$ consists of constants of motion of the form $Tg^{-1}+U$. Then
	\[
		0=\left\{\mathcal{H},Tg^{-1}+U\right\} = g^{-1}\left\{g^{-1},T\right\}+(T\left\{V,g^{-1}\right\}+\left\{g^{-1},U\right\}).
	\] 
	Since $\left\{g^{-1},T\right\} = 0\iff {d}{}_{a}T=0$ and $M$ is connected, we conclude that $T = C_{1}\in \mathbb{R}$. This gives
	\[
	0 = \left\{\mathcal{H},Tg^{-1}+U\right\} = \left\{\mathcal{H},T\mathcal{H}+(U-TV)\right\} = \left\{\mathcal{H},(U-TV)\right\} = \left\{g^{-1},(U-TV)\right\}
	\] 
	and we similarly conclude that $U-TV = C_{2}\in \mathbb{R}$. Hence,
	\[
	Tg^{-1}+U = T\mathcal{H}+\tilde{U} = C_{1}\mathcal{H}+C_{2}
	\] 
	So the kernel of $\Phi$ is spanned by constant functions and $\mathcal{H}$.
\end{proof}

\section{Tractors}\label{Sec:Tct}

Here we give a construction of the standard conformal tractor bundle for a conformal manifold $(M,\boldsymbol{c})$. See \cite{bailey_thomass_1994,cap_tractor_2000,curry_introduction_2018}. 
\subsection{Conformal Tractors}
We first need the conformal decomposition of curvature. Let $g\in \boldsymbol{c}$, and let $R{}_{ab}{}^{c}{}_{d} $ be its Riemann curvature tensor, given by $R{}_{ab}{}^{c}{}_{d} v^d = 2\nabla{}_{[a} \nabla{}_{b]} v^c$. The curvature tensor has symmetries $R{}_{abcd}  = -R{}_{bacd} =-R{}_{abdc} $ and the Bianchi identity $R{}_{[abc]d} =0$. We can decompose $R$ into trace and trace-free parts,
\[
	R{}_{abcd}  = W{}_{abcd} +2\boldsymbol{g}{}_{c[a} \mathsf{P}{}_{b]d} +2\boldsymbol{g}{}_{d[b} \mathsf{P}{}_{a]c}.
\] 
Since we lowered the $c$ index using the conformal metric, the above identity is between weight two tensors. The trace-free part $W{}_{abcd} $ is called the \textit{Weyl Curvature}, it is conformally invariant, and has the same symmetries as $R{}_{abcd} $ with the addition of $W{}_{ab}{}^{a}{}_{d} =0$. The tensor $\mathsf{P}{}_{bd} $ is called the \textit{conformal Schouten curvature} and is symmetric. We write $\mathsf{J} = \mathsf{P}{}_{b}{}^{b} $ and we have the identities
\[
	 \operatorname{Ric}{}_{bd}  = R{}_{ab}{}^{a}{}_{d}  = (n-2)\mathsf{P}{}_{bd} +\mathsf{J}\boldsymbol{g}{}_{bd}\quad \text{and} \quad \operatorname{Sc} = R{}_{ab}{}^{a}{}_{d}g{}^{bd} =2(n-1)\sigma^2 \mathsf{J}.
\] 
Under a change of metric $\hat{g} = \Omega^2g$, the Weyl curvature does not change, and the Schouten transforms by
\[
\hat{\mathsf{P}}{}_{ab}  = \mathsf{P}{}_{ab}  -\nabla_a\Upsilon_b +\Upsilon_a\Upsilon_b-\tfrac{1}{2}g{}_{ab} \Upsilon_c\Upsilon^c.
\] 
Consider the Jet bundles, $\mathcal{J}^k(\mathcal{E}[1])$ for the conformal $1$-density bundle. Through Taylor series, we can obtain the jet exact sequences for all $\mathcal{J}^k(\mathcal{E}[1])$. The sequence for $\mathcal{J}^1(\mathcal{E}[1])$ is 
\begin{equation}\label{Eqn:ES1}
	0\longrightarrow \mathcal{E}_a[1]\longrightarrow \mathcal{J}^1(\mathcal{E}[1])\longrightarrow \mathcal{E}[1]\longrightarrow 0.	
\end{equation} 
We write $\mathcal{J}^1(\mathcal{E}[1]) = \mathcal{E}[1]\loplus\mathcal{E}_a[1]$ to reflect the filtration of $\mathcal{J}^1(\mathcal{E}[1])$ given by the above sequence.
As discussed in Section \ref{Sec:CfG}, a choice of metric $g\in \boldsymbol{c}$ gives a connection on $\mathcal{E}[1]$. A connection on $\mathcal{E}[1]$ is equivalent to a splitting of the above sequence, where $j^1_q(\sigma)\in \mathcal{J}^1(\mathcal{E}[1])_q$ maps to $\nabla_a\sigma\in  \mathcal{E}_a[1]_q$. The exact sequence at $\mathcal{J}^2(\mathcal{E}[1])$ is
\begin{equation}\label{Eqn:ES2}
	0\longrightarrow{\mathcal{E}}{}_{(ab)}[1]\longrightarrow\mathcal{J}^2(\mathcal{E}[1])\longrightarrow\mathcal{J}^1(\mathcal{E}[1])\longrightarrow 0
\end{equation} 
The conformal metric gives an invariant decomposition of $\mathcal{E}{}_{(ab)} [1]$ into trace and trace-free parts, $\mathcal{E}{}_{(ab)} [1] = \mathcal{E}{}_{(ab)_0} [1]\oplus \mathcal{E}[-1]\otimes\boldsymbol{g}{}_{ab}$. The \textit{conformal cotractor bundle} is then the quotient 
\[
	{\mathcal{E}}{}_{A} = \mathcal{T}^* = \frac{\mathcal{J}^2(\mathcal{E}[1])}{{\mathcal{E}}{}_{(ab)_0}[1]}.
\]
This definition of the tractor bundle immediately implies the existence of a second order operator given by the 2nd jet operator $j^2(q\mapsto \sigma) = (q\mapsto j^2_q(\sigma))$ composed with the projection onto the above quotient. The definition of the tractor bundle is characterised by an exact sequence
\[
	0\longrightarrow \mathcal{E}{}_{(ab)_0}[1]  \longrightarrow \mathcal{J}^2(\mathcal{E}[1]) \longrightarrow \mathcal{E}{}_{A}  \longrightarrow 0,
\] 
which is invariantly split on the left by the \textit{almost Einstein operator},
\[
	AE(\sigma)_{ab} = \nabla{}_{(a} \nabla{}_{b)_0} \sigma +\mathsf{P}{}_{(ab)_0}\sigma.
\] 
This identifies $\mathcal{E}{}_{A} $ with the kernel of $AE$ inside $\mathcal{J}^2(\mathcal{E}[1])$, see \cite{bailey_thomass_1994}. 
Taking quotients in the exact sequence \eqref{Eqn:ES2} gives the tractor exact sequence,
\[
	0\longrightarrow \mathcal{E}[-1] \longrightarrow {\mathcal{E}}_{A} \longrightarrow \mathcal{J}^1(\mathcal{E}[1]) \longrightarrow0 
\] 
which implies the invariant filtration $\mathcal{E}_A = \mathcal{J}^1(\mathcal{E}[1])\loplus \mathcal{E}[-1] = \mathcal{E}[1]\loplus\mathcal{E}_a[1]\loplus\mathcal{E}[-1]$. A choice of metric $g\in \boldsymbol{c}$ induces an isomorphism $\mathcal{E}_A\cong \mathcal{E}[1]\oplus \mathcal{E}_a[1]\oplus \mathcal{E}[-1]$ of vector bundles, given by 
\begin{equation}\label{Eqn:j2q}
j^2_q(\sigma)\mapsto (\sigma,\nabla_a\sigma,-\frac{1}{n}(\Delta\sigma+\mathsf{J}\sigma)).
\end{equation} 
When working in a scale $g\in \boldsymbol{c}$, we usually write a tractor $U_A\in \Gamma(\mathcal{E}_A)$ in terms of its components, $U_A \overset{g}{=}  (\sigma,\mu_a,\rho)$. Under a change of scale, $\hat{g} = \Omega^2 g$, the components transform as
\begin{equation}\label{Eqn:TctC}
	U_A \overset{g}{=} (\sigma,\mu_a,\rho) \overset{\hat{g}}{=} (\sigma,\mu_a+\sigma\Upsilon_a, \rho-\Upsilon^b\mu_b-\tfrac{1}{2}\sigma\Upsilon^b\Upsilon_b),
\end{equation}
where $\Upsilon_a = \Omega^{-1}\nabla_a\Omega$. There is a canonical inclusion, $\mathcal{E}[-1]\to \mathcal{E}{}_{A} $, which we represent by the canonical tractor $X_A\in \Gamma(\mathcal{E}_A[1])$, which in any scale is given by $X_A = (0,0,1)$. There is also a canonical projection, $\mathcal{E}_A \to \mathcal{J}^1(\mathcal{E}[1])\to \mathcal{E}[1]$, which is represented by a cotractor, $X^A\in \Gamma(\mathcal{E}^A[1])$, and we have $X^AU_A = \sigma$. 

\begin{figure}[tbh]	
\begin{tikzcd}
            &  &  & 0 \arrow[d]   &   \\
            & & & {\mathcal{E}_a[1]} \arrow[d] \arrow[ld, "{Z}{}^{a}{}_{A}"', dashed, bend right, shift right] &   \\
0 \arrow[r] & {\mathcal{E}[-1]} \arrow[r, "X_A"] & \mathcal{E}_A \arrow[r] \arrow[ru, "{Z}{}^{A}{}_{a}"', dashed] \arrow[rd, "X^A"] \arrow[l, "Y^A", dashed, bend left, shift left] & {\mathcal{J}^1(\mathcal{E}[1])} \arrow[r] \arrow[d]  & 0 \\
	    &  &  & {\mathcal{E}[1]} \arrow[d] \arrow[lu, "Y_A", dashed, bend left, shift left]   &   \\
            &  &  & 0   &  
\end{tikzcd}
\caption{\label{Fig:T} The tractor projectors}
\end{figure}
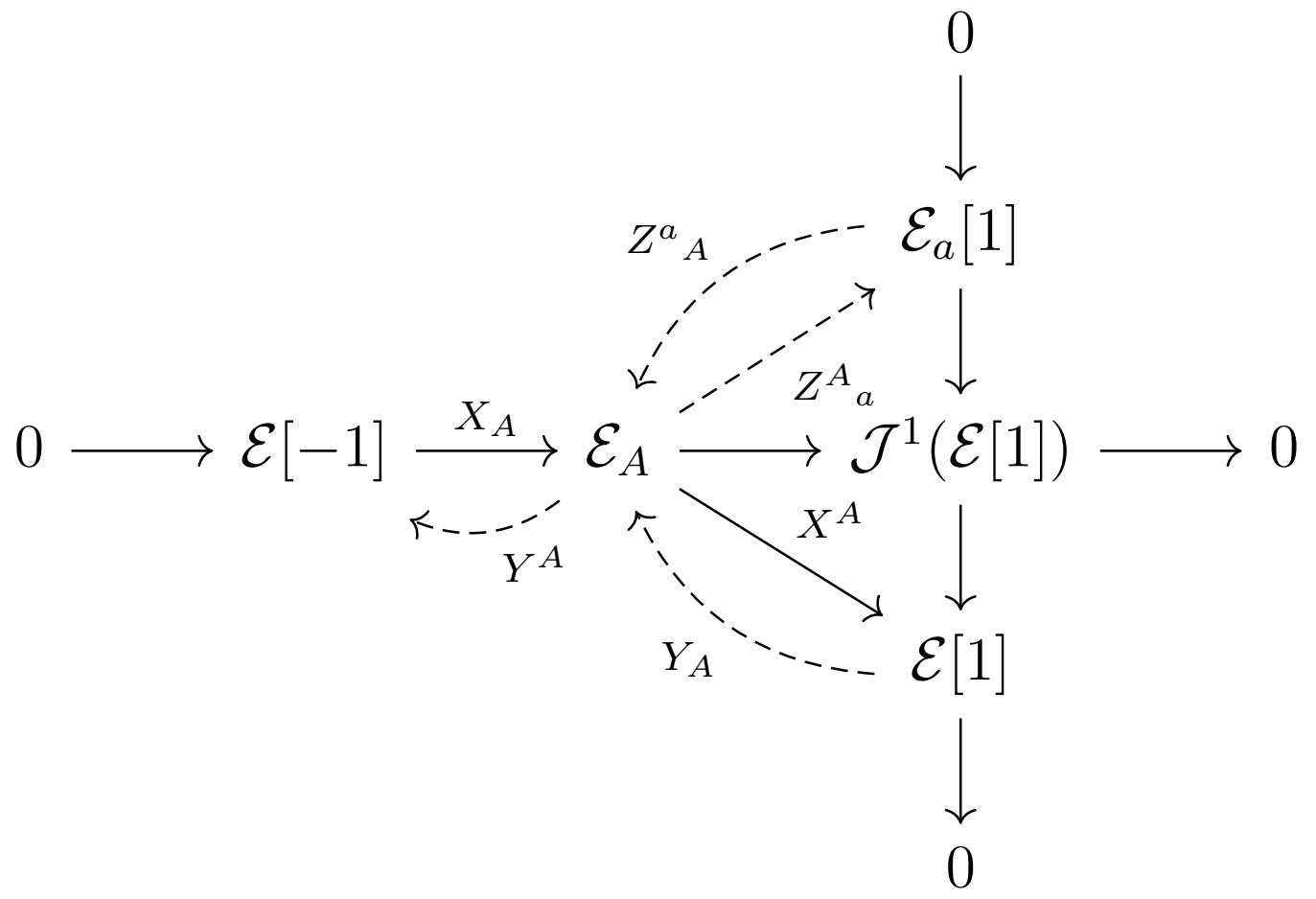

Given a metric $g\in \boldsymbol{c} $, there are two other projectors, 
\[
	Z{}^{A}{}_{a} U{}_{A}  = \mu_a \quad \text{and} \quad Y^AU_A = \rho
\]
that are not invariant. These projectors come from the scale dependent splitting \eqref{Eqn:j2q}, and we can read off the transformation laws from \eqref{Eqn:TctC},
\[
\hat{X}^A=X^A, \quad  {\hat{Z}}{}^{A}{}_{a} = {Z}{}^{A}{}_{a}+{\Upsilon}{}_{a}{X}{}^{A}, \quad {\hat{Y}}{}^{A} = {Y}{}^{A} - {\Upsilon}{}^{a}{Z}{}^{A}{}_{a}-\tfrac{1}{2}{\Upsilon}{}_{a}{\Upsilon}{}^{a}X^A.
\] 
The tractor bundle has an unweighted invariant non-degenerate bilinear form,
\[
{h}{}^{AB} = {X}{}^{A}{Y}{}^{B}+{Y}{}^{A}{X}{}^{B}+{Z}{}^{A}{}^{a}{Z}{}^{B}{}_{a}
\] 
called the \textit{tractor metric}. The tractor metric allows us to raise and lower tractor indices, relating the projectors ${X}{}^{A}$, ${Y}{}^{A}$, ${Z}{}^{A}{}_{a}$ to the corresponding inclusion maps ${X}{}_{A}$, ${Y}{}_{A}$ ${Z}{}_{A}{}^{a}$. All the projectors and the splittings are illustrated in Figure \ref{Fig:T}. 

The tractor bundle has a natural invariant calculus, given by the conformally invariant tractor connection, which is defined by

\[
\nabla_a U_A = \nabla_a\left( \sigma,\mu_b,\rho \right) = (\nabla_a \sigma-\mu_a,\nabla_a\mu_b + \mathsf{P}{}_{ab}\sigma+\boldsymbol{g}{}_{ab} \rho, \nabla_a\rho -\mathsf{P}{}_{a}{}^{b} \mu{}_{b}  ),
\] 
or equivalently in terms of the projectors,
\begin{equation}\label{Eqn:TcCon}
{\nabla}{}_{a}X^A = {Z}{}^{A}{}_{a},\quad {\nabla}{}_{a}{Z}{}^{A}{}_{b} = -{\mathsf{P}}{}_{ab}{X}{}^{A}-{\boldsymbol{g}}{}_{ab}Y^A,\quad {\nabla}{}_{a}{Y}{}^{A} = {\mathsf{P}}{}_{a}{}^{b}{Z}{}^{A}{}_{b},
\end{equation} 
where in \eqref{Eqn:TcCon} all connections are the tractor connection paired with a connection on ${\mathcal{E}}{}_{a}$ and $\mathcal{E}[w]$ from the conformal class. The tractor metric and its inverse are parallel, so the formulae in \eqref{Eqn:TcCon} hold with indices raised.

The tractor bundle and its associated calculus are particularly suitable for prolongation of invariant PDEs. The most fundamental example is the conformally invariant almost Einstein equation \cite{bailey_thomass_1994}
\[
AE(\sigma)_{ab} = {\nabla}{}_{(a}{\nabla}{}_{b)_0}\sigma+{\mathsf{P}}{}_{(ab)_0}\sigma=0,
\] 
for $\sigma\in \Gamma(\mathcal{E}[1])$. Solutions of the almost Einstein equation are in one-to-one correspondence with parallel unweighted standard tractors; ${I}{}_{A}\in \Gamma(\mathcal{E}_A)$ such that ${\nabla}{}_{a}{I}{}_{A}=0$. The corresponding solution is $\sigma = {X}{}^{A}{I}{}_{A}$ and the inverse operation from solutions to tractors is given by the second order tractor $D$ operator, $D\colon \Gamma(\mathcal{E}[1])\to \Gamma(\mathcal{E}_A) $ given in a scale by
\begin{equation}\label{Eqn:TcD}
	{D}{}_{A}\sigma =  {Y}{}_{A}\sigma+{Z}{}^{a}{}_{A}{\nabla}{}_{a}\sigma-{X}{}_{A} \frac{\Delta\sigma +\mathsf{J}\sigma}{n},
\end{equation}
which is easily checked to be conformally invariant. In fact, this is the invariant operator defined earlier in this section, induced by the 2nd jet operator $j{}^{2} \colon \Gamma(\mathcal{E}[1])\to \Gamma(J^2(\mathcal{E}[1])) $ composed with the projection onto the tractor bundle.

Let $\sigma\in \Gamma(\mathcal{E}[1])$ be a scale density. Then the tractor $I_A = D_A\sigma$ is called the \textit{scale tractor} of $\sigma$. On the open dense $M\setminus\mathcal{Z}(\sigma)$, working in the splitting induced by $g = \sigma^{-2}\boldsymbol{g}$, we have that
\[
I_A = (\sigma,0,-\tfrac{1}{n}(\mathsf{J}^{g}\sigma)).
\] 

The tractor $D$ is induced by the invariant \textit{unnormalised tractor} $\tilde{D}$ operator \cite{bailey_thomass_1994} which is defined for any weighted conformal tractor bundle, with any number of tractor indices and any conformal weight,
\[
\tilde{D}_A V = (n+2w-2)wY_AV + (n+2w-2)Z{}_{A}{}^{a} \nabla_a V-X_A(\Delta V+w\mathsf{J}V),
\] 
where we have suppressed any tractor indices of $V$.
The normalised tractor $D$ is then defined by $D = \left( (n+2w-2)w \right)^{-1} \tilde{D}$ for any $w\not\in \{0,\frac{2-n}{2}\} $, and by $D = (n+2w-2)^{-1} \tilde{D}$ when $w=0$. We will only use the tractor $D$ \eqref{Eqn:TcD} and the tractor $D\colon \Gamma(\mathcal{E}{}_{A} )\to \Gamma(\mathcal{E}{}_{AB} [-1])$ by
\[
D_A U_B = Z{}_{A}{}^{a} \nabla_a I_B - \frac{1}{n-2}X_A\Delta I_B .
\] 
When applying this tractor $D$ to a scale tractor $I_A = D_A \sigma$, and working in the splitting induced by $g = \sigma^{-2}\boldsymbol{g}$, we have
\[
	D_A I_B = Z{}_{A}{}^{a} Z{}_{B}{}^{b} \mathsf{P}_{(ab)_{0}}-\frac{1}{n}\nabla_a\mathsf{J} Z{}_{A}{}^{a} X{}_{B} -\frac{1}{n}\nabla_b \mathsf{J} Z{}_{B}{}^{b} X_A+\frac{1}{n-2}\left( \frac{1}{n}\Delta\mathsf{J} + \mathsf{P}^{(ab)_0}\mathsf{P}_{(ab)_0} \right) 
\] 

When considering tractor bundles with skew symmetries, it is useful to define new projectors. We let
\[
	\mathbb{X}{}^{AB}{}_{b} = X{}^{[A} Z{}^{B]}{}_{b},\quad \mathbb{W}{}^{AB} = X{}^{[A} Y{}^{B]} ,\quad \mathbb{Z}{}^{AB}{}_{ab} =Z{}^{[A}{}_{a} Z{}^{B]}{}_{b} ,\quad \mathbb{Y}{}^{AB}{}_{b}  = Y{}^{[A} Z{}^{B]}{}_{b}.
\] 
The tractor $\mathbb{X}$ is canonical, that is, it does not depend on the scale we are calculating in, since the symmetries kill the transformation of $Z$. We will not need the general transformation laws, although they are simple to derive. Instead we will need the connection applied to these projectors,

\begin{equation}\label{Eqn:XZWY}
\nabla_a \mathbb{X}{}^{AB}{}_{b} = \mathbb{Z}{}^{AB}{}_{ab} -\boldsymbol{g}{}_{ab} \mathbb{W}{}^{AB} ,\quad \nabla_a \mathbb{W}{}^{AB}  = \mathsf{P}{}_{a}{}^{b} \mathbb{X}{}^{AB}{}_{b} -\mathbb{Y}{}^{AB}{}_{a},
\end{equation} 

\[
	\nabla_a \mathbb{Z}{}^{AB}{}_{bc} = 2\boldsymbol{g}{}_{a[c} \mathbb{Y}{}^{AB}{}_{b]} +2\mathsf{P}{}_{a[c} \mathbb{X}{}^{AB}{}_{b]},\quad \nabla_a \mathbb{Y}{}^{AB}{}_{b} = \mathsf{P}{}_{ab} \mathbb{W}{}^{AB} -\mathsf{P}{}_{a}{}^{c} \mathbb{Z}{}^{AB}{}_{bc}  
\] 
\subsection{BGG Machinery}\label{Sec:BGG}
Bernstein-Gelfand-Gelfand (BGG) sequences in parabolic geometry \cite{cap_bernstein-gelfand-gelfand_2001, calderbank_differential_2001} are invariant sequences of differential operators associated to any irreducible tractor bundle. In the context of conformal geometry, a tractor bundle refers to a $\mathfrak{g} = \mathfrak{so}(p+1,q+1)$ invariant subbundle of a tensor power of the standard tractor bundle, equipped with the signature $(p+1,q+1)$ bilinear form $h{}_{AB}$. In particular the irreducible tractor bundles are the trace-free tractors with symmetries given by a Young diagram. 

Associated to a tractor bundle $\mathcal{V}$ is the \textit{twisted de Rham} sequence,
\[
	0\longrightarrow \mathcal{V} \overset{d^\nabla}{\longrightarrow} T^*M\otimes \mathcal{V} \overset{d^\nabla}{\longrightarrow} \Lambda^2T^*M\otimes\mathcal{V} \overset{d^\nabla}{\longrightarrow} \dots,
\] 
where $d^{\nabla}(V_{a_{1}\dots a_k\; A}) = \nabla{}_{[a} (V_{a_{1}\dots a_k]\; A_{1}\dots A_j})$, where $A_{1}\dots A_j$ are the tractor indices of $\mathcal{V}$ and $\nabla$ is the tractor connection paired with any torsion free affine connection on $M$. When the tractor connection is flat, the twisted de Rham sequence is a complex. The Kostant codifferential is an algebraic operator acting in the opposite direction,
\[
0\overset{\partial^*}{\longleftarrow} \mathcal{V} \overset{\partial^*}{\longleftarrow} T^*M\otimes \mathcal{V} \overset{\partial^*}{\longleftarrow} \Lambda^2T^*M\otimes\mathcal{V} \overset{\partial^*}{\longleftarrow} \dots.
\] 
For conformal tractor bundles, we have
\[
	\partial^*(V_{a_{1}\dots a_k\; A_{1}\dots A_j}) = \sum_{i=1}^{j} 2\mathbb{X}_{A_i}{}^B{}^{a_1}V_{a_{1}\dots a_k A_{1}\dots A_{i-1} B A_{i+1}\dots A_j},
\]
where $\mathbb{X}{}_{AB}{}^{b}  = X{}_{[A} Z{}_{B]}{}^{b} $ is an invariant bundle inclusion, $\mathbb{X}\colon T^*M\hookrightarrow \mathcal{E}_{[AB]}$. The homology bundles are defined by $H_k(\mathcal{V}) = \operatorname{ker}(\partial^*)/\operatorname{im}(\partial^*)$, where both are subbundles of $\Lambda^kT^*M\otimes\mathcal{V}$. We denote the projection onto homology by $\Pi\colon  \operatorname{ker}(\partial^*)\to H_k(\mathcal{V})$ and note that for $k=0$, this gives a map $\Pi\colon \mathcal{V}\to H_0(\mathcal{V}) $. For irreducible parabolic geometries, including conformal, the homology bundles $H_k(\mathcal{V})$ are irreducible $\mathfrak{g}_0$ bundles. For conformal geometry, $\mathfrak{g}_0 = \mathfrak{co}(P,Q) = \mathfrak{so}(P,Q)\oplus \mathbb{R}$, so the homology bundles are trace-free irreducible weighted tensor bundles. 

The Bernstein-Gelfand-Gelfand machinery gives two important families of operators. The first are the BGG splitting operators. 
\begin{theorem}[\!\!\cite{cap_bernstein-gelfand-gelfand_2001} Theorem 2.5 and Lemma 2.7]\label{Thm:BGGS}
	There exist invariant splitting operators, $L\colon \Gamma(H_k(\mathcal{V}))\to \Gamma(\operatorname{ker}(\partial^*))$ such that $\Pi\circ L = \operatorname{id}_{H_k(\mathcal{V})}$ and $\partial^*\circ d^{\nabla}\circ L=0$. The subspace $\operatorname{im}(L)\subset \Gamma(\Lambda^k(T^*M)\otimes\mathcal{V})$ consists of sections $s$ such that $\partial^*(s)=0$ and $\partial^*\circ d^{\nabla}(s)=0$.
\end{theorem}
The above conditions are sufficient to calculate the splitting operators. Take the $\mathcal{V} = \mathcal{E}_A$, standard conformal cotractor bundle. We can compute the action of $\partial^*$ on $\mathcal{E}{}_{aA} $ in a splitting,
\begin{equation}\label{Eqn:dsT}
	\partial^*(\sigma_aY_A+\mu_{ab} Z{}^{b}{}_{A} +\rho_aX_A) = 2\mathbb{X}{}_{A}{}^{Ba} (\sigma_aY_B+\mu_{ab}Z^b_B+\rho_aX_B) = -\sigma{}_{b} Z{}_{A}{}^{b}  +\mu{}_{a}{}^{a} X_A.
\end{equation} 
So suppose $\sigma\in \Gamma(H_0(\mathcal{E}_A)) = \mathcal{E}[1]$ and that $L(\sigma)_A = (\sigma,\mu_a,\rho)$. Then $d^{\nabla}\circ L(\sigma) = \nabla_a(\sigma,\mu_b,\rho) = (\nabla_a\sigma-\mu_a,\nabla_a\mu_b+ \sigma\mathsf{P}{}_{ab} +\rho\boldsymbol{g}{}_{ab} , \nabla_a\rho- \mathsf{P}{}_{a}{}^{b} \mu_b) $. So the condition $\partial^*\circ d^{\nabla}\circ L=0$ becomes $\nabla_a\sigma-\mu_a = 0$ and $\nabla_a\mu^a+\mathsf{J}\sigma+n\rho=0$. Hence, $\mu_a = \nabla_a\sigma$ and $\rho=-\frac{\Delta+\mathsf{J}}{n}\sigma$. This shows that the first example of a splitting operator is the (normalised) tractor $D$ operator \eqref{Eqn:TcD}. 

The splitting operators give rise to the second family of operators, the BGG differential operators $\mathscr{D}_k:= \Pi\circ d^{\nabla}\circ L \colon \Gamma(H_k(\mathcal{V})) \to \Gamma(H_{k+1}(\mathcal{V}))$ and the BGG sequence,
\[
	0\longrightarrow \Gamma(H_0(\mathcal{V}))\overset{\mathscr{D}_0}{\longrightarrow} \Gamma(H_1(\mathcal{V}))\overset{\mathscr{D}_1}{\longrightarrow} \Gamma(H_2(\mathcal{V}))\overset{\mathscr{D}_2}{\longrightarrow}\cdots 
\]
which in the flat case is a complex. In general, the squares in the following diagram do not commute, however in the flat case they do.
\[
\begin{array}{ccccccccccc}
	0
&\longrightarrow
&\Gamma(\mathcal{V})
& \xrightarrow{\; d^\nabla \;}
& \Gamma\!\big(\Lambda^1(T^*M)\otimes \mathcal{V}\big)
& \xrightarrow{\; d^\nabla \;}
& \Gamma\!\big(\Lambda^2(T^*M)\otimes \mathcal{V}\big)
& \xrightarrow{\; d^\nabla \;}
& \cdots
\\[0.2em]
&&\big\uparrow\, L_0
&
&
\big\uparrow\, L_1
&
&
\big\uparrow\, L_2
&
\\[0.2em]
0
&\longrightarrow
&\Gamma\!\big({H}_0(\mathcal{V})\big)
& \xrightarrow{\; \mathscr{D}_0 \;}
& \Gamma\!\big({H}_1(\mathcal{V})\big)
& \xrightarrow{\; \mathscr{D}_1 \;}
& \Gamma\!\big({H}_2(\mathcal{V})\big)
& \xrightarrow{\; \mathscr{D}_2 \;}
& \cdots
\end{array}
\]
\begin{theorem}\label{Thm:BGGL}
	Let $(M,\boldsymbol{c})$ be a conformally flat manifold and let $\mathcal{V}$ be an irreducible tractor bundle. Then there is a one-to-one correspondence between tensors in the kernel of the $k$-th BGG operator $\mathscr{D}_k \colon \Gamma(H_k(\mathcal{V}))\to \Gamma(H_{k+1}(\mathcal{V})) $ and sections of $\Lambda^k(T^*M)\otimes\mathcal{V}$ which are in the kernel of $\partial^*$ and $d^{\nabla}$. 

	In particular, for $k=0$ there is a one-to-one correspondence between solutions of $\mathscr{D}_0 s=0$ and parallel tractors in $\mathcal{V}$.
\end{theorem}
For $\mathcal{V} = \mathcal{E}{}_{A} $, the image of $\partial^*$ can be determined from \eqref{Eqn:dsT} to be the subspace spanned by $Z$ and $X$. So the projection onto homology is given by $\sigma Y_A+\mu_a Z{}^{a}{}_{A} +\rho X_A\mapsto \sigma$, equivalently $U_A\mapsto X^AU_A$. More generally for an irreducible tractor bundle with certain symmetries, the projection $\Pi\colon \mathcal{V}\to H_0(\mathcal{V}) $ is given by first applying as many $X$ projectors as is allowed by the symmetries and then applying $Z$ projectors to the remaining indices. For a tractor $V\in \Gamma(\mathcal{V})$, the tensor $\Pi(V)$ is called the \textit{top slot} of $V$, the projection is always conformally invariant. The following Lemma is frequently useful to move from tensor to tractor equations.

\begin{lemma}\label{Lem:BGG}
	Let $\mathcal{V}$ be an irreducible tractor bundle and let $V\in\Gamma(\mathcal{V})$ be parallel. If $\Pi(V)=0$, then $V=0$.
\end{lemma}
\begin{proof}
	By Theorem \ref{Thm:BGGS}, $V = L\circ \Pi(V) = 0$, since $d^{\nabla}(V) = 0$ and $ \partial^*(V) = 0$.
\end{proof}
For any particular tractor bundle, one can prove the above lemma by brute force by differentiating and showing that all slots of $V$ vanish. For example, if $V{}_{A}\in \Gamma(\mathcal{E}_A)$ is parallel and $X^AV_A=0$, then, in any scale,
\[
	\nabla_a(X^AV_A) = Z{}^{A}{}_{a}  V_A = 0 \implies \nabla_b( Z{}^{A}{}_{a} V_A) = -\mathsf{P}{}_{ab}  X^AV_A-\boldsymbol{g}_{ab} Y^AV_A = -\boldsymbol{g}_{ab} Y^AV_A=0,
\] 
so $X^AV_A = 0$, $Z{}^{A}{}_{a} V_A=0$ and $Y^AV_A=0$, so $V=0$. The calculations get more tedious as the number of tractor indices increases, making the BGG machinery a convenient route to take. 

Aside from the above, we are genuinely interested in many BGG operators in this article. We are primarily interested in $\mathscr{D}_0$ for various tractor bundles. We will refer to $\mathscr{D}_0$ as the \textit{BGG operator} associated to $\mathcal{V}$.

The conformal Killing operator \eqref{Eqn:CKE} is a BGG operator for all ranks $\ell$. The rank one conformal Killing vector operator is the BGG operator for the tractor bundle of antisymmetric tractors, also called the adjoint tractor bundle, $\mathcal{A} = \mathcal{E}{}_{[AB]} $ and the rank two conformal Killing operator is the BGG operator for the tractor bundle consisting of degree four trace-free tractors with curvature symmetries, 
\[
	\mathcal{W} = \{W\in {\mathcal{E}}{}_{ABCD}:W{}_{ABCD} = - W{}_{BACD} = -W{}_{ABDC},\; {W}{}_{AB}{}^{A}{}_{D}=0,\; {W}{}_{[ABC]D}=0\},
\]
represented by the trace-free Young diagram $\begin{picture}(16, 6.5)(0,1.5)
	\put(-0.2, 10){\line(1,0){10.4}}
	\put(0, 5){\line(1,0){10}}
	\put(-0.2, 0){\line(1,0){10.4}}
	\put(0, 10){\line(0,-1){10}}
	\put(5, 10){\line(0,-1){10}}
	\put(10, 10){\line(0,-1){10}}
	\put(10, 0){${}_{\,0}$}
\end{picture}(\mathcal{E}_A) $. Theorem \ref{Thm:BGGL} applied to these tractor bundles give the following,
\begin{theorem}\label{Thm:W}
	Let $(M,\boldsymbol{c})$ be conformally flat. There is a one-to-one correspondence between
	\begin{tabel}
		\item \label{Thm:R:1} conformal Killing vectors $k_{a}\in \Gamma(\mathcal{E}_a[2])$,
		\item \label{Thm:R:3} parallel antisymmetric conformal tractors $K_{AB}$.
	\end{tabel}
	There is a one-to-one correspondence between
		\begin{tabel}
		\item \label{Thm:W:1} rank $2$ conformal Killing tensors $k_{ab}\in \Gamma(\mathcal{E}{}_{(ab)_0}[4])$,
		\item \label{Thm:W:3} parallel conformal tractors $W_{ABCD}$ with Weyl tensor symmetries.
	\end{tabel}
\end{theorem}
The tractor $W_{ABCD}$ is called the \textit{full prolongation} of ${k}{}_{ab}$.

The maps \ref{Thm:R:3}$\to $\ref{Thm:R:1} are given by taking the top slot, that is $k_b = X{}^{A} Z{}^{B}{}_{b} K{}_{AB} $ and ${k}{}_{bd} = {X}{}^{A}{Z}{}^{B}{}_{b}{X}{}^{C}{Z}{}^{D}{}_{d}{W}{}_{ABCD}$

The splitting operators, $L$ which give the maps \ref{Thm:R:1}$\to $\ref{Thm:R:3} are harder to construct. For rank one conformal Killing tensors, the direct approach we applied to the standard tractor bundle is still reasonably simple, however this approach for rank two tensors is far less tractable. For rank one and two conformal Killing tensors on conformally flat manifolds, the splitting operators were invariantly derived in \cite{gover_conformal_nodate}. In that paper, and extended in \cite{goodwin_killing_2026}, the authors derived algebraic conditions between full prolongations of conformal Killing tensors and scale tractors that characterise Killing scales, which is very important to our current application.

The results of those papers are summarised below.
\subsection{Prolongation of Killing Tensors}\label{Sec:CKT}

We will need the following results from \cite{gover_conformal_nodate} and \cite{goodwin_killing_2026}. The first proposition can also be found in \cite{gover_detecting_2014}

\begin{proposition}\label{Prp:ConPro}
	Let $I_A = D_A\sigma$ be a parallel scale tractor. Then consider the pseudo-Riemannian manifold $(M\setminus\mathcal{Z}(\sigma),g)$, where $g = \sigma^{-2}\boldsymbol{g}$. There is a connection preserving isomorphism between the projective tractor bundle on $M\setminus\mathcal{Z}(\sigma)$ and the subbundle of the conformal tractor bundle perpendicular to $I_A$ on $M\setminus\mathcal{Z}(\sigma)$. The projective projectors on the projective bundle agree with the projectors $X^A$ and $Z{}^{A}{}_{a} $ on $I^{\perp}$.
\end{proposition}
\begin{theorem}[\cite{gover_conformal_nodate} Theorem 4.18, \cite{goodwin_killing_2026} for null $I$]\label{Thm:KT}
	Let $(M,\boldsymbol{c})$ be conformally flat and let $I_A = D_A\sigma$ be a parallel scale tractor. There is a one-to-one correspondence between 
	\begin{tabel}
		\item\label{Prp:ConPro:1:1} Sections $k_a\in \Gamma(\mathcal{E}_a[2])$ that are Killing for the metric $g=\sigma^{-2}\boldsymbol{g}$ and
		\item\label{Prp:ConPro:1:2} parallel antisymmetric conformal tractors $K_{AB}$ that satisfy $I^AK_{AB}=0$.
	\end{tabel}
	There is a one-to-one correspondence between 
	\begin{tabel}
		\item\label{Prp:ConPro:2:1} Sections $r{}_{bd} \in \Gamma(\mathcal{E}_{(bd)}[4])$ that are Killing for the metric $g=\sigma^{-2}\boldsymbol{g}$ and
		\item\label{Prp:ConPro:2:2} parallel conformal tractors $R_{ABCD}$ with Riemann curvature symmetries that satisfy $I^AR_{ABCD}=0$.
	\end{tabel}
\end{theorem}
Later in this article, we will need the following Proposition, which says that the parallel tractor $R{}_{ABCD} $ is independent of the Einstein scale.
\begin{proposition}\label{Prp:Run}
	Let $(M,\boldsymbol{c})$ be conformally flat and let $I_A = D_A\sigma$ and $J_A = D_A\tau$ be parallel scale tractors. Let $R{}_{ABCD}$ be a parallel tractor with curvature symmetries such that $I^A R{}_{ABCD} =0$. Let $r{}_{bd} = X^A Z{}^{B}{}_{b} X^C Z{}^{D}{}_{d} R{}_{ABCD} $ and suppose that $r{}_{bd}$ is a Killing tensor for the metric $\tau^{-2}\boldsymbol{g}$. Then $J^A R{}_{ABCD}=0$.
\end{proposition} 
\begin{proof}
	By Theorem \ref{Thm:KT}, there is a parallel tractor $\hat{R}{}_{ABCD} $ which is parallel, $J^A \hat{R}{}_{ABCD} =0$ and $X^A Z{}^{B}{}_{b} X^C Z{}^{D}{}_{d} \hat{R}_{ABCD} = r{}_{bd} $. So $X^A Z{}^{B}{}_{b} X^C Z{}^{D}{}_{d} (R{}_{ABCD} -\hat{R}{}_{ABCD} )=0$. The conformal tractor bundle with curvature symmetries splits as
	\begin{equation}\label{Eqn:Run}
	\begin{picture}(16, 11.0)(0, 5.0)
	\put(-0.2, 16){\line(1,0){16.4}}
	\put(-0.2, 8){\line(1,0){16.4}}
	\put(-0.2, 0){\line(1,0){16.4}}
	\put(0, 16){\line(0,-1){16}}
	\put(8, 16){\line(0,-1){16}}
	\put(16, 16){\line(0,-1){16}}
\end{picture}(\mathcal{E}_A) = \begin{picture}(22, 11.0)(0, 5.0)
	\put(-0.2, 16){\line(1,0){16.4}}
	\put(-0.2, 8){\line(1,0){16.4}}
	\put(-0.2, 0){\line(1,0){16.4}}
	\put(0, 16){\line(0,-1){16}}
	\put(8, 16){\line(0,-1){16}}
	\put(16, 16){\line(0,-1){16}}
	\put(16, 0){${}_{\,0}$}
\end{picture}(\mathcal{E}_A)\oplus \begin{picture}(22, 7.0)(0, 1.0)
	\put(-0.2, 8){\line(1,0){16.4}}
	\put(-0.2, 0){\line(1,0){16.4}}
	\put(0, 8){\line(0,-1){8}}
	\put(8, 8){\line(0,-1){8}}
	\put(16, 8){\line(0,-1){8}}
	\put(16, 0){${}_{\,0}$}
\end{picture}(\mathcal{E}_A)\oplus \mathcal{E}.
	\end{equation} 
	The projector $X^A Z{}^{B}{}_{b} X^C Z{}^{D}{}_{d} $ directly extracts the top slot of the trace-free part. Therefore, by Lemma \ref{Lem:BGG}, $R-\hat{R}$ has zero projection onto the first direct summand. The top slot of the second summand is given by
	\[
		X^A X^C h{}^{BD} (R-\hat{R}){}_{ABCD} = X^A Z{}^{B}{}_{b} X^C Z{}^{D}_d g{}^{bd}  (R-\hat{R})_{ABCD} = r{}_{b}{}^{b} -r{}_{b}{}^{b} =0.
	\] 
	So, by Lemma \ref{Lem:BGG}, $R-\hat{R}$ has zero projection onto the second direct summand of \eqref{Eqn:Run}. The final summand is spanned by the tractor $(h\owedge h)_{ABCD} = 4h{}_{A[C} h{}_{D]B} $. So, $R{}_{ABCD} -\hat{R}{}_{ABCD} = \lambda h{}_{A[C} h{}_{D]B} $. However, since $I^A R{}_{ABCD}  = 0 = J^A \hat{R}{}_{ABCD} $,
	\[
		0 = I^AJ^B\left( R{}_{ABCD} -\hat{R}{}_{ABCD}  \right) = \lambda I^A J^B h{}_{A[C} h{}_{D]B}  = \lambda I{}_{[C}  J{}_{D]}.
	\] 
	Hence, either $I$ is proportional to $J$ and the conclusion is trivial, or $\lambda=0$, so $R = \hat{R}$ and $J^A R{}_{ABCD} =0$.
\end{proof}The analogous result for Killing vectors is clear since $\Lambda^2(\mathcal{E}_A)$ is already an irreducible tractor bundle, so $K{}_{AB} $ is determined by its top slot, $k_b = X^A Z{}^{B}{}_{b} K{}_{AB} $.

We call the tractors $K{}_{AB} $ and $R{}_{ABCD} $ \textit{full prolongations} analogously to prolongations of conformal Killing tensors. 
\begin{theorem}[\cite{gover_conformal_nodate} Theorem 1.3, \cite{goodwin_killing_2026} for null $I$]\label{Thm:IWI}
Let $(M,\boldsymbol{c})$ be conformally flat and let $k_{ab}$ be a conformal Killing tensor. Let $W_{ABCD}$ be the full prolongation of ${k}{}_{ab}$. Then there is a one-to-one correspondence between
\begin{tabel}
\item \label{Thm:IWI:1} scales $\sigma$ such that the metric $g_{ab}=\sigma^{-2}\boldsymbol{g}_{ab}$ is Einstein, and $k_{ab}$ is the tracefree part of a Killing tensor for $g_{ab}$,
 \item  \label{Thm:IWI:3} parallel tractors $I_C$ that satisfy
\begin{equation} \label{Eqn:IWI}
    I^AW_{AB[CD}I_{E]}=0.
\end{equation}
\item \label{Thm:IWI:2} parallel tractors ${I}{}_{A}$ for which there exists a parallel ${R}{}_{ABCD}$ with Riemann curvature symmetries satisfying
	\begin{equation}
		{I}{}^{A}{R}{}_{ABCD}=0\quad \text{and}\quad {W}{}_{ABCD} = \operatorname{trace-free}({R}){}_{ABCD}
	\end{equation}	
\end{tabel}
\end{theorem}

The generalisation to non-Einstein scales is proved in \cite{goodwin_killing_2026}. A scale $\sigma\in \Gamma(\mathcal{E}[1])$ is a \textit{local Killing scale} for $k{}_{ab} $ if for all $q\in M\setminus\mathcal{Z}(\sigma)$, there is an open neighbourhood $U$ of $q$ on which $g = \sigma^{-2}\boldsymbol{g}$ is a Killing scale for $k{}_{ab} $.
\begin{theorem}\label{Thm:KS2}
	Let $(M,\boldsymbol{c})$ be conformally flat. Let $k_{ab}$ be a conformal Killing tensor. Let $W_{ABCD}$ be the full prolongation of ${k}{}_{ab}$. Let $\sigma\in \Gamma(\mathcal{E}[1])$ and let $g_{ab}=\sigma^{-2}\boldsymbol{g}_{ab}$ on $M\setminus\mathcal{Z}(\sigma)$. Then $\sigma$ is a local Killing scale for $k_{ab}$ if and only if $I_A=D_A\sigma$ satisfies
      \begin{equation}\label{Eqn:3II}
3I^AX^B Z{}^{C}{}_{c}  Z{}^{D}{}_{d}  X^E W_{AB[CD}I_{E]}
=\sigma X^B Z{}^{C}{}_{c}  Z{}^{D}{}_{d}  X^EW_{AB[CD}D_{E]}I^A,
\end{equation}
\end{theorem}
Due to the characterisation in Theorem \ref{Thm:QCns}, we will be interested in conformal Killing tensors which are Killing in multiple Killing scales. We will make frequent use of the above Theorems to determine which scales are Killing scales. The general method is as follows. Start with an Einstein scale, $\sigma$ with parallel scale tractor $I_A=  D_A\sigma$. Killing tensors for $\sigma$ correspond to parallel tractors $R{}_{ABCD} $ with curvature symmetries that are orthogonal to $I$. Given such an $R$, take its trace-free part,
\begin{equation}\label{Eqn:tfW}
	W{}_{ABCD} = R{}_{ABCD} -(h\owedge P)_{ABCD} = R{}_{ABCD} - h{}_{AC} P{}_{BD} - h{}_{BD} P{}_{AC} +h{}_{AD} P{}_{BC} +h{}_{BC} P{}_{AD}, 
\end{equation} 
where
\[
	P{}_{BD}  = \frac{1}{n}\left(R{}_{AB}{}^{A}{}_{D} - \frac{1}{2(n+1)}R{}_{AC}{}^{AC} h{}_{BD}\right).
\]
Then one uses Theorem \ref{Thm:IWI} or Theorem \ref{Thm:KS2} to test if another scale is a Killing scale for $W$.

In particular, if the second scale is Einstein, corresponding to the parallel tractor ${J}{}_{A}$, then the following result shows that we don't need to compute $W$.

\begin{theorem}\label{Thm:IRI}
	Let $(M,\boldsymbol c)$ be conformally flat, let ${J}{}_{A} = D{}_{A}{\tau}$ be a parallel tractor and let $R_{ABCD}$ be a parallel tractor with Riemann curvature symmetries. Suppose that
	\[
		{J}{}^{A}{R}{}_{AB[CD}{J}{}_{E]}=0.
	\] 
	Then ${\tau}$ is a Killing scale for $k = \operatorname{top}(W)$, where $W = \operatorname{trace-free}(R)$.
\end{theorem}
\begin{proof}
	We have $W = \operatorname{trace-free}(R)$ as given by \eqref{Eqn:tfW}. So 
	\[
		J^A W{}_{ABCD}  = J^A R{}_{ABCD} -2J{}^{A} P{}_{A[C} h{}_{D]B} -2J_{[C} P{}_{D]B}
	\]
	So, since $J^A R{}_{AB[CD} J{}_{E]} =0$,
	\begin{align*}
		J^A W{}_{AB[CD} J{}_{E]}  = 2J{}^{A} P{}_{A[C} J_E h{}_{D]B}  = \frac{2}{n}J^A R{}_{AB'[C}{}^{B'} J_E h{}_{D]B} 
	\end{align*}
	However, also by assumption, we have
	\[
		0 = \operatorname{tr}{}_{BD} \left( J^A R{}_{AB[CD} J{}_{E]}  \right) = \frac{1}{3}\left( J^A R{}_{ABC}{}^{B} J_E +J^A R{}_{AB}{}^{B}{}_{E} J{}_{C} \right) = \frac{2}{3}J^A R{}_{AB[C}{}^{B} J{}_{E]} 
	\]
	So $J^A W{}_{AB[CD} J{}_{E]} =0$ and $\tau$ is a Killing scale by Theorem \ref{Thm:IWI}.
\end{proof}

\begin{corollary}\label{Cor:IwU}
Let $(M,\boldsymbol{c})$ be conformally flat and let ${I}{}_{A}$ be a parallel tractor and let ${R}{}_{ABCD} = {(I^2\owedge U)}{}_{ABCD}$ for some parallel $U\in \Gamma({\mathcal{E}}{}_{(AB)})$. Then $\sigma = {X}{}^{A}{I}{}_{A}$ is a Killing scale for the conformal Killing tensor ${k}{}_{(ab)} = {r}{}_{(ab)_0}$, where $r = \operatorname{top}(R)$. 
\end{corollary}
The following Lemma was proved for non-null $I$ in \cite{gover_conformal_nodate}. The proof in the null case is similar, but since we will need it, we will include the proof, including null $I$, for completeness.

\begin{lemma}\label{Lem:GABCD}
	Let $I_A\in \mathcal{E}_A|_q$ be a non-zero tractor at a point $q\in M$. Then the space of tractors with curvature symmetries which are orthogonal to $I$ and are have no trace-free part is spanned by
	\[
		(h\owedge (I^2-\tfrac{\iota}{2}h))_{ABCD},
	\] 
	where $\iota = I_AI^A$.
\end{lemma}
\begin{proof}
	Fix $I_A\in \mathcal{E}_A|_q$. Let ${P}{}_{BD} \in S^2(\mathcal{E}_A|_q)$ and suppose that $I^A(h\owedge P)_{ABCD}=0$. Then
	\begin{equation}\label{Eqn:IhP}
		I^A(h\owedge P)_{ABCD} = 2{P}{}_{B[D} I{}_{C]} +2I{}^{A} P{}_{A[C} h{}_{D]B} = 0
	\end{equation} 
	So $I{}^{A} P{}_{A[C} I{}_{E} h{}_{D]B} =0$, which implies $I^A P{}_{A[C} I{}_{E]} =0$, since $h$ is non-degenerate. Which then gives $I^A P{}_{AC}  = \lambda_{1}I_C$ for some $\lambda_{1}\in \mathbb{R} $. Substituting into \eqref{Eqn:IhP} gives
	\[
		P{}_{B[D} I{}_{C]} + \lambda_{1}I{}_{[C} h{}_{D]B} =0,
	\] 
	so $P{}_{BD} -\lambda_{1}h{}_{BD}  = \lambda_{2}I_BI_D$ for some $\lambda_{2}\in \mathbb{R} $. Then $I^A P{}_{AC}  = (\lambda_{2}\iota-\lambda_{1})I_C = \lambda_{1}I_C$. So $\lambda_{1} = -\frac{\iota}{2}\lambda_{2}$ and $P{}_{BD}  = \lambda_{2}(I_BI_D-\frac{\iota}{2} h{}_{BD}) $.
\end{proof}

Our main approach will be to first construct $2n-2$ functionally independent Killing tensors which are Killing for a whole family of scales, given by functional combinations of two Einstein scales (see Proposition \ref{Prop:fKS}), and then only use \eqref{Eqn:3II} for the final constant.

\section{Functional Independence}\label{Sec:FI}
A set of functions $f_1,\dots ,f_m \colon S\to \mathbb{R} $ defined on a manifold $S$ is \textit{functionally independent} if the differentials $df_{1},\dots ,df_m$ are linearly independent in $T^*S$ on an open dense subset of $S$. Clearly if $m>\dim(S)$, then the functions $f_{1},\dots ,f_m$ are functionally dependent.

Let $f_{1},\dots ,f_m$ be constants of motion for a Hamiltonian $\mathcal{H}$, all defined on an open subset of $T^*M$. Since $df_{1},\dots ,df_m$ are all orthogonal to $d \mathcal{H}$, with respect to the symplectic form $\omega$ on $T^*M$, the maximal number of functionally independent constants of motion for $\mathcal{H}$ is $2n-1$. The Hamiltonian is \textit{maximally superintegrable} if there are $2n-1$ functionally independent constants of motion for $\mathcal{H}$.

The differential of a function $f\colon T^*M\to \mathbb{R} $ is a section of $T^*(T^*M)$. Let $(q,p)\in T^*M$. Then there are canonical dual exact sequences of vector spaces,

\[\begin{tikzcd}
	0 & {T^*_qM} & {T_{(q,p)}(T^*M)} & {T_qM} & 0, \\
	0 & {T^*_qM} & {T^*_{(q,p)}(T^*M)} & {T_qM} & 0,
	\arrow[from=1-1, to=1-2]
	\arrow["\varphi", from=1-2, to=1-3]
	\arrow["{d\pi}", from=1-3, to=1-4]
	\arrow[from=1-4, to=1-5]
	\arrow[from=2-1, to=2-2]
	\arrow["{\pi^*}", from=2-2, to=2-3]
	\arrow["\varphi^*", from=2-3, to=2-4]
	\arrow[from=2-4, to=2-5]
\end{tikzcd}\]where the first map of the first sequence is the identification of the vertical bundle $V{}_{(q,p)} = \ker(d\pi)\subset T{}_{(q,p)} (T^*M)$ with $T^*_qM$. To construct this map explicitly, take $\omega_q\in T_q^*M$ and consider the curve in $T^*M$ originating from $(q,p)$ given by $(q,p+t\omega_q)$, the initial velocity of this curve is the image of $\omega{}_{q} $ in $T{}_{(q,p)} (T^*M)$. 

A splitting of the first sequence identifies $T{}_q M$ with a subspace $H{}_{(q,p)} \subset T{}_{(q,p)} (T^*M)$, usually called the \textit{horizontal subspace}, giving a direct sum decomposition $T{}_{(q,p)} (T^*M) = V{}_{(q,p)} \oplus H{}_{(q,p)} $. A smooth collection of splittings for all $(q,p)\in T^*M$ gives a horizontal subbundle $H\subset T(T^*M)$ such that $T(T^*M) = H\oplus V$, such a collection of splittings is equivalent to a connection on $T^*M$ (see, for example, \cite{kolar_natural_1993, lang_fundamentals_1999, kobayashi_foundations_1963}). Given a connection $\nabla$, the horizontal subbundle is generated by the initial velocities of all parallel transports of $\nabla$, which are all curves in $T^*M$. 

A splitting of the sequence for $T(T^*M)$ is equivalent to a dual splitting of the sequence for $T^*(T^*M)$, which we will denote at $(q,p)\in T^*M$ by
\begin{equation}\label{Eqn:TTM}
T^*_{(q,p)}(T^*M)\overset{\nabla}{=} T_qM\oplus T_q^*M.
\end{equation} 
This setting is where our interest lies since $df\in \Gamma(T^*(T^*M))$ for constants of motion $f$. So we explicitly compute for polynomial functions $f$, the horizontal and vertical components of $df{}_{(q,p)} $ and the corresponding element of $T_q^*M\oplus T_qM$. 

\begin{proposition}\label{Prp:df}
	Let $\nabla$ be an affine connection and let $f\colon T^*M\to \mathbb{R} $ be induced by a vector field, $f(q,p) = k_q(p)$. Then under the isomorphism $T^*_{(q,p)}(T^*M))\overset{\nabla}{=} T_q^*M\oplus T_qM$, we have
	\[
	df \overset{\nabla}{=} (k^a,p_a\nabla_bk^a)
	\] 
\end{proposition}
\begin{proof}
	The canonical component $\phi^*(df_{(q,p)})$ acts on $\omega_q\in T^*_qM$ by
	\[
		\langle \phi^*(df_{(q,p)}), \omega{}_{q}  \rangle = \langle df_{(q,p)},\phi(\omega)\rangle = \left.\frac{d}{dt}\right|_{t=0}k_q(p+t\omega_q) = k_q(\omega_q).
	\] 
	So $\varphi^*(df) = k^a$. Let $v\in T_qM$ and let $\hat{v}\in T_{(q,p)}(T^*M)$ be the horizontal lift given by $\nabla$, which is the initial velocity of a parallel transport curve $(p(t),q(t))$ with $q'(0)=v$. Then
	\[
		df(\hat{v}) = \left.\frac{d}{dt}\right|_{t=0}k_{q(t)}(p(t)) = \nabla_v\langle k, p \rangle = \langle \nabla_v k, p \rangle.  
	\] 
	Hence, $df \overset{\nabla}{=} (k^a,p_a\nabla_bk^a)$.
\end{proof}
By first considering simple sections of $S^\ell(TM)$, that is, $r_{a_{1}\cdots a_\ell} = k^1{}_{(a_{1}}\cdots k^\ell{}_{a_{\ell)}}$, one can easily extend the above proposition. If $r\in \Gamma(S^\ell(TM))$, then the differential of the induced function, also denoted $r$, on $T^*M$ maps to
\begin{equation}\label{Eqn:dr}
	dr_{(q,p)} \overset{\nabla}{=} (\ell{p}{}_{a_2}\dots {p}{}_{a_{\ell}}{r}{}^{b\,a_2\cdots a_{\ell}},{p}{}_{a_{1}}\dots {p}{}_{a_\ell}{\nabla}{}_{b}{r}{}^{a_{1}\cdots a_\ell})|_q\in T_qM \oplus T^*_qM
\end{equation} 

The first constant of motion taken for a Hamiltonian $\mathcal{H} = g^{-1}+V$ is usually $\mathcal{H}$ itself. Using the Levi-Civita connection of $g$ for the identification \eqref{Eqn:TTM} gives
\begin{equation}\label{Eqn:HJac}
	 (2{p}{}_{a} {g}{}^{ab}, {p}{}_{a}{p}{}_{c}{\nabla}{}_{b}{g}{}^{ac}+{\nabla}{}_{b}V) = (2{p}{}^{b},{\nabla}{}_{b}V)\in T_qM \oplus T^*_qM.
\end{equation} 

Let ${r}{}_{ab}$ be a Killing tensor for a metric $g\in \boldsymbol{c}$. By Theorems \ref{Thm:W} and \ref{Thm:KT} we can associate parallel tractors $R{}_{ABCD} $ and $W{}_{ABCD} $ to ${r}{}_{ab}$ and its trace-free part $k{}_{ab}  = {r}{}_{(ab)_0}$, respectively. Then using Theorem \ref{Thm:QCns} together with Theorems \ref{Thm:IWI} \ref{Thm:IWI:3} and \ref{Thm:KS2}, we can test for the existence of a function $U$ making $f = r+U$ a constant of motion for $\mathcal{H} = g^{-1}+V$. However, constructing the function $U$ and relating it to the tractors $R{}_{ABCD} $ and $W{}_{ABCD} $ is complicated. Therefore it becomes difficult to work with the entire differential $df = dr+dU$. However, using linear algebra we can show that it is sufficient to consider $dr$.

\begin{proposition}\label{Prp:Ind}
	Let $f_{1},\dots ,f_m$ be polynomial functions on $T^*M$ and let $r_{1},\dots ,r_m$ be contravariant tensor fields representing the highest order terms of $f_{1},\dots ,f_m$. If the functions induced by $r_{1},\dots ,r_m$ are functionally independent, then so are the functions $f_{1},\dots ,f_m$.
\end{proposition}
\begin{proof}
	Let $\nabla$ be an affine connection on $M$. For each $i\in \{1,\dots ,m\} $, let $\ell{}_{i} $ be the order of $f{}_{i}$, so $r^i\in \Gamma(S^{\ell_i}(TM))$. Then under the identification \eqref{Eqn:dr},
	\[
		dr_i \overset{\nabla}{=} (\ell_i\,{p}{}_{a_2}\dots {p}{}_{a_{\ell_i}}{r}{}^{b\,a_2\cdots a_{\ell_i}},{p}{}_{a_{1}}\dots {p}{}_{a_{\ell_i}}{\nabla}{}_{b}{r}{}^{a_{1}\cdots a_{\ell_i}})|_q
	\]
	and by assumption, $dr_1,\dots ,dr_m$ are linearly independent for almost all $(q,p)\in T^*M$. Then let $F\colon T^*M\to \mathbb{R}^m $ be the function $F(q,p) = (f_{1}(q,p),\dots ,f_m(q,p))$. Fix $(q,p)\in T^*M$ and consider the differential $dF$ restricted to $\{(q,\lambda p):\lambda\in \mathbb{R},\lambda\neq 0\} $. Then $dF$ is a matrix with rows
	\[
	dF|_{(q,\lambda p)}\overset{\nabla}{=} \left(\lambda^{\ell_i-1}\ell_i{p}{}_{a_2}\dots {p}{}_{a_{{\ell_i}}}{r}_i{}^{\mspace{-6mu}b\,a_2\cdots a_{{\ell_i}}}+l.o.ts,\lambda^{\ell_i}{p}{}_{a_{1}}\dots {p}{}_{a_{\ell_i}}{\nabla}{}_{b}{r}_i{}^{\mspace{-6mu}a_{1}\cdots a_{\ell_i}}+l.o.ts\right)_{1\le i\le m}	
	\] 
	where $l.o.ts$ stands for lower order terms, with respect to the power of $\lambda$. Let $\ell$ be the maximum polynomial order among the $f_{1},\dots ,f_m$. We now multiply the first $n$ columns of $dF$ by $\lambda$ and multiply each row by $\lambda^{\ell-\ell_i}$. This operation does not change the rank of the matrix. Which gives
	\begin{align*}
		\operatorname{rank}(dF|_{(q,\lambda p)}) 
		&= \operatorname{rank} (\lambda^\ell \ell_i{p}{}_{a_2}\dots {p}{}_{a_{{\ell_i}}}{r}_i{}^{\mspace{-6mu}b\,a_2\cdots a_{{\ell_i}}}+l.o.ts,\lambda^\ell{p}{}_{a_{1}}\dots {p}{}_{a_{\ell_i}}{\nabla}{}_{b}{r}_i{}^{\mspace{-6mu}a_{1}\cdots a_{\ell_i}}+l.o.ts)_{1\le i\le m}\\ 
		&=\operatorname{rank}(\lambda^\ell A_\ell+\lambda^{\ell-1}A_{\ell-1}+\dots \lambda A_{1}+A_0)
	\end{align*} 
	where $A_{\ell}$ is the matrix with rows 
	\[
		dr_i|_{(q,p)}\overset{\nabla}{=} (\ell_i\,{p}{}_{a_2}\dots {p}{}_{a_{\ell_i}}{r}{}^{b\,a_2\cdots a_{\ell_i}},{p}{}_{a_{1}}\dots {p}{}_{a_{\ell_i}}{\nabla}{}_{b}{r}{}^{a_{1}\cdots a_{\ell_i}})|_q,
	\] 
	so $A{}_{\ell} $ has rank $m$. We now need the following result from linear algebra. Let 
	\[
	A(\lambda) = \lambda^{\ell}A_\ell+\dots +\lambda A_1+A_0
	\] 
	be a matrix pencil. If $A_\ell$ has rank $R$, then $\operatorname{rank}(A(\lambda))\ge R$ for all but finitely many $\lambda\in \mathbb{R}$.

	Indeed, for nonzero $\lambda$, let $\alpha=\lambda^{-1}$ and consider
	\[
	B(\alpha) = \alpha^\ell A_0+\dots \alpha A_{\ell-1}+A_\ell.
	\] 
	Since $B(0) = A_\ell$ has rank $R$ and the rank of the pencil $B(\alpha)$ can only drop at finitely many $\alpha\in \mathbb{R}$, it follows that $\lambda^{\ell}B(\alpha) = A(\lambda)$ has rank at least $R$ for all but finitely many $\lambda\in \mathbb{R}$.

	By assumption, there is an open dense subset $U\subset T^*M$ on which $dr_{1},\dots ,dr_m$ are linearly independent. If $(q,p)\in U$, then $\operatorname{rank}(dF) = \operatorname{rank}(A_\ell) = m$ for all but finitely many values of $\lambda$. So, letting $V\subset T^*M$ be the open set on which $dF$ has rank $m$, we have shown that $U\subset \overline{V}$, since $U$ is dense, this implies $\overline{V}=T^*M$.
\end{proof}

Let $k_b\in \Gamma(\mathcal{E}_b)$ be a conformal Killing vector. Recall the BGG prolongation of $k_b$ from Section \ref{Sec:BGG}, $K{}_{AB} \in \Gamma(\mathcal{E}{}_{[AB]} )$ satisfying $\mathbb{X}{}^{AB}{}_{b} K{}_{AB} =k_b$ and $2\mathbb{X}{}_{[B}{}^{Cb} \nabla{}_{|b|} K{}_{A]C} =0$. One can verify by direct computation that these conditions, when $k_b$ conformal Killing, are equivalent to 
\begin{equation}\label{Eqn:XZWYK}
	\begin{gathered}
		\mathbb{X}{}^{AB}{}_{b} K{}_{AB} =k_b,\quad  \mathbb{Z}{}^{AB}{}_{ab} K{}_{AB} =\nabla{}_{[a} k{}_{b]} \\
	\mathbb{W}{}^{AB} K{}_{AB} =-\frac{1}{n}\nabla{}_{b} k{}^{b}, \quad \textrm{and}\quad \mathbb{Y}{}^{AB}{}_{b} K{}_{AB} =  \mathsf{P}{}_{b}{}^{a} k_a+\tfrac{1}{n}\nabla_b\nabla_ak^a ,
	\end{gathered}
\end{equation}
see similar calculations in \cite{gover_conformal_nodate,gover_conformal_2008-1}. Alternatively we write the above as
\[
	K{}_{AB}  = \begin{pmatrix} k_b\\ \nabla{}_{[a} k{}_{b]} ,\quad  -\frac{1}{n} \nabla_bk^b\\  \mathsf{P}{}_{b}{}^{a} k_a+\tfrac{1}{n}\nabla_b\nabla_ak^a  \end{pmatrix},
\] 
similarly to how we represent a tractor as in \eqref{Eqn:TctC}. We use the full prolongation to capture the identification \eqref{Eqn:TTM} for $dk$ in a conformally invariant way.
\begin{theorem}\label{Thm:PhiPsi1}
	Let $(M,\boldsymbol{c}) $ be a conformal manifold and let $g\in \boldsymbol{c}$. Let $k_a\in \Gamma(\mathcal{E}_a [2])$ be a conformal Killing vector, with full prolongation $K{}_{AB} $. Let $\nabla$ be the Levi-Civita connection of $g$ and let $(q,p)\in T^*M$. Denote the identification \eqref{Eqn:TTM} by $\Phi\colon T^*_{(q,p)}(T^*M)\to T_qM\oplus T^*_qM $. Define $\Psi\colon T_qM\oplus T^*_qM\to \mathcal{E}{}_{[AB]} |_q $ by
	\[
		\Psi(v^a,\omega_b) = \begin{pmatrix} 0 \\ 2v{}_{[a} p{}_{b]} ,\; v^bp_b \\-\omega_b\end{pmatrix},
	\] 
	where we use the isomorphism $\mathcal{E}{}_{[AB]} \cong \mathcal{E}_a[2]\oplus \mathcal{E}_{[ab]}[2]\oplus \mathcal{E}\oplus \mathcal{E}_b$ determined by $g\in \boldsymbol{c}$. The map $\Psi$ is injective whenever $p$ is not null, and	\begin{equation}\label{Eqn:dkKS}
		\Psi\circ \Phi(dk) = 4K{}_{[A}{}^{C} S{}_{B]C},
	\end{equation} 
	where $S{}_{AB}  = \mathbb{X}{}_{AB}{}^{b} p_b$.
\end{theorem}
\begin{proof}
	The kernel of $\Psi$ is spanned by $(p_a,0)$ when $p$ is null and is trivial otherwise. By Proposition \ref{Prp:df}, $\Phi(dk_{(q,p)}) = (k^a,p_a\nabla_b k^a)$. Then, by \eqref{Eqn:XZWYK},
	\[
		4K{}_{[A}{}^{C} S{}_{B]C} = \begin{pmatrix} 0 \\ 2k{}_{[a} p{}_{b]} ,\quad k^bp_b \\ -\frac{1}{n}\nabla_ak^a p_b- \nabla{}_{[b} k{}_{a]} p^a\end{pmatrix} = \begin{pmatrix} 0 \\ 2k{}_{[a} p{}_{b]} ,\quad k^bp_b \\ -p_a\nabla_bk^a \end{pmatrix},
	\] 
	since $\nabla{}_{(a} k{}_{b)_0}=0 $. Which gives \eqref{Eqn:dkKS}.
\end{proof}
The points $(q,p)\in T^*M$ such that $p$ is null for the conformal metric forms a nowhere dense subset of $T^*M$. Which gives the following corollary.
\begin{corollary}\label{Thm:Find}
	Let $k^1_b,\dots ,k^m_b\in \Gamma(\mathcal{E}_b[2])$ be conformal Killing vectors, with full prolongations $K{}^{1}{}_{AB} ,\dots ,K{}^{m}{}_{AB} $. Then the functions induced by the $k_i$ are functionally independent if and only if the set of points $(q,p)\in T^*M$ for which
	\[
		{K^1}{}_{[A}{}^{C}{S}{}_{B]C} ,\dots ,{K^m}{}_{[A}{}^{C}{S}{}_{B]C}\in {\mathcal{E}}{}_{[AB]}|_q \quad \textrm{are linearly independent}
	\]
	is dense in $T^*M$, where ${S}{}_{AB} = X{}_{[A} Z{}_{B]}{}^{b}  {p}{}_{b}$.
\end{corollary}Second order constants of motion for $\mathcal{H} = g^{-1}+V$ have leading term a Killing tensor for $g$. We can extend the above to Killing tensors for almost Einstein scales in the conformally flat setting using the prolongation from Theorem \ref{Thm:KT}.

\begin{theorem}\label{Thm:PhiPsi2}
	Let $(M,\boldsymbol{c}) $ be a conformal manifold and let $g\in \boldsymbol{c}$. Let $R{}_{ABCD} $ be a parallel tractor with curvature symmetries and let $r{}_{bd} = \mathbb{X}{}^{AB}{}_{b} \mathbb{X}{}^{CD}{}_{d} R{}_{ABCD} $. Fix $(q,p)\in T^*M$ and let $\Phi\colon T^*_{(q,p)}(T^*M)\to T_qM\oplus T^*_qM $ and $\Psi\colon T_qM\oplus T^*_qM\to \mathcal{E}{}_{[CD]} |_q $ be the linear maps as in Theorem \ref{Thm:PhiPsi1}. Then 
	\begin{equation}\label{Eqn:drSRS}
		\Psi\circ \Phi(dr) = 8S{}^{AB} R{}_{ABE[D} S{}_{C]}{}^{E} ,
	\end{equation} 
	where $S{}_{AB}  = \mathbb{X}{}_{AB}{}^{b} p_b$.
\end{theorem}
\begin{proof}
	For the higher rank tractor $R$, we calculate only the slots which we absolutely need, we have
	\[
		8S{}^{AB} R{}_{ABE[D} S_{C]}{}^E Z{}^{C}{}_{c} Z{}^{D}{}_{d} = -4p^bX^A Z{}^{B}{}_{b} X{}^{C} Z{}^{D}{}_{[d} p{}_{c]}R{}_{ABCD}  = 4p^b r{}_{b[c} p{}_{d]},
	\] 
	\[
		8S{}^{AB} R{}_{ABE[D} S_{C]}{}^E X{}^{C} Y{}^{D} = -2p^bX^A Z{}^{B}{}_{b} Z{}^{C}_c p{}^{c}X^DR{}_{ABCD}  = 2p^bp^d r{}_{bd}.
	\]
	For the $Y^C Z{}^{D}{}_{d} $ slot we need
	\[
		\nabla{}_{c} r{}_{bd}  = 2X^A  Z{}^{B}{}_{(b} Z{}^{C}_{|c|} Z{}^{D}{}_{d)} R{}_{ABCD} -2g{}_{c(d} X^A Z{}^{B}{}_{b)} X^C Y^D R{}_{ABCD},
	\] 
	Then we have
	\begin{align*}
		8S{}^{AB} R{}_{ABE[D} S_{C]}{}^E Y{}^{C} Z{}^{D}{}_{c} 
		&= 2p^bp^dX^A Z{}^{B}{}_{b} Z{}^{E}{}_d Z{}^{D}{}_{c} R{}_{ABED} + 2p_cp^b X^AZ{}^{B}{}_{b}X^CY^D R{}_{ABCD}  \\ 
		&= -(\nabla_c r{}_{bd}) p^bp^d,\textrm{ and}
		\end{align*}	
	\begin{equation}\label{Eqn:SRS}
		8{S}{}^{AB}{R}{}_{ABE[D}S_{C]}{}^E = \begin{pmatrix} 0\\ 4p^b r{}_{b[c} p{}_{d]},\quad 2p^bp^d r{}_{bd} \\ -(\nabla_a r{}_{bd}) p^bp^d \end{pmatrix} = \Psi\circ \Phi(dr)
	\end{equation} 
	by \eqref{Eqn:dr}.
\end{proof}
	There is a simpler way to show \eqref{Eqn:SRS} for parallel curvature tractors of the form
	\begin{equation}\label{Eqn:Find2}
		R^i{}_{ABCD}  = K{}^{1}{}_{AB} K{}^{2}{}_{CD} +K{}^{2}{}_{AB} K{}^{1}{}_{CD} -2K{}^{1}{}_{[AB} K{}^{2}{}_{CD]}
	\end{equation} 
	with $X{}^{A} Z{}^{B}{}_{b} X{}^{C} Z{}^{D}{}_{d} R^i{}_{ABCD}  = r{}_{bd}  = 2k{}^{1}{}_{(b} k{}^{2}{}_{d)} $, where $k{}^{i}{}_{b}  = X{}^{A} Z{}^{B}{}_{b} K{}^{i}{}_{AB} $. It is straightforward to calculate that
	\[
		S{}^{AB} R^i{}_{ABE[D} S_{C]}{}^E  =  k{}^{1}{}_{a} p^a K{}^{2}{}_{E[D} S_{C]}{}^E +k{}^{2}{}_{a} p^a K{}^{1}{}_{E[D} S_{C]}{}^E   ,
	\] 
	which immediately gives \eqref{Eqn:SRS}. On a conformally flat manifold, it is clear that locally tractors of the form \eqref{Eqn:Find2} span the space of curvature tractors. 

Combining the results in this section gives the framework which we will use to talk about functional independence of constants of motion.
\begin{corollary}\label{Cor:Find}
	Let $(M,\boldsymbol{c})$ be conformally flat. Let $\sigma\in \Gamma(\mathcal{E}[1])$ be an almost Einstein scale, let $g = \sigma^{-2}\boldsymbol{g}$ and let $\mathcal{H} = g^{-1}+V$ be a Hamiltonian defined on $M\setminus\mathcal{Z}(\sigma)$. 

	Let $f_{1},\dots ,f_{m_1}$ be linear constants of motion corresponding to conformal Killing vectors $k{}^{1}{}_{b},\dots ,k^{m_1}{}_b\in \Gamma(\mathcal{E}_b[2])$ with full prolongations $K{}^{1}{}_{AB} ,\dots ,K{}^{m_1}{}_{AB} $. 

	Let $f_{m_{1}+1},\dots ,f_{m_{1}+m_{2}}$ be quadratic constants of motion with quadratic part induced by Killing tensors for $g$, $r{}^{1}{}_{bd} ,\dots ,r^{m_{2}}{}_{bd}\in \Gamma(\mathcal{E}{}_{bd} [4])$ with full prolongations, as in Theorem \ref{Thm:KT}, $R{}^{1}{}_{ABCD} ,\dots ,R{}^{m_{2}}{}_{ABCD}$. 

	Then if the set of points $(q,p)\in T^*M$ for which the tractors 
	\begin{equation}\label{Eqn:KSSRS}
		K^1{}_{[C}{}^{E}{S}{}_{D]E} ,\dots ,K^{m_1}{}_{[C}{}^{E}{S}{}_{D]E}, \;2{S}{}^{AB}{R^1}{}_{ABE[D}S_{C]}{}^E,\dots ,2{S}{}^{AB}{R^{m_2}}{}_{ABE[D}S_{C]}{}^E
	\end{equation} 
	are linearly independent in ${\mathcal{E}}{}_{[CD]}|_q,$ is dense in $T^*M$, where $S{}_{AB}  = \mathbb{X}{}_{AB}{}^{b} p_b$, then the constants $f_{1},\dots ,f_{m_{1}+m_{2}}$ are functionally independent.
\end{corollary}

The following Lemma gives explicit conditions to check whether the tractors \eqref{Eqn:KSSRS} are linearly independent. We will only need the non-null case.
\begin{lemma}\label{Lem:sl}
	Let $(q,p)\in T^*M$ and let ${S}{}_{AB} = X_{[A}Z_{B]}{}^b p{}_{b} $. Suppose that $p$ is not null. Then ${T}{}_{AB}\in {\mathcal{E}}{}_{[AB]}|_q$ satisfies ${T}{}_{[A}{}^{C}{S}{}_{B]C}=0$ if and only if in all scales
	\begin{equation}\label{Eqn:sl1}
		{\mathbb{X}}{}^{AB}{}_{b}{T}{}_{AB}=0,\quad {\mathbb{W}}{}^{AB}{T}{}_{AB}=0,\quad \text{and} \quad {p}{}^{a}\mathbb{Z}{}^{AB}{}_{ab} {T}{}_{AB}=0
	\end{equation}
	The first two requirements can be replaced by $X^A T{}_{AB} =0$. In this case, the image of the map $T\mapsto T{}_{[A}{}^{C} S{}_{B]C} $ has dimension $2n$.

	Otherwise if $p$ is null, then ${T}{}_{[A}{}^{C}{S}{}_{B]C}=0$ if and only if in all scales
	\begin{equation}\label{Eqn:sl2}
		{\mathbb{X}}{}^{AB}{}_{b}{T}{}_{AB}=\lambda p_b,\quad {\mathbb{W}}{}^{AB}{T}{}_{AB}p_b+{p}{}^{a}\mathbb{Z}{}^{AB}{}_{ab}  {T}{}_{AB}=0.
	\end{equation} 
\end{lemma}
\begin{proof}
	Given the following form of $T$ in some scale,
	\begin{equation}\label{Eqn:SdotT}
	T{}_{AB}  = \begin{pmatrix} t_b \\ \nu{}_{ab} ,\; \tau \\ * \end{pmatrix},
	\;\;\text{we have that}\;\;
		4{T}{}_{[A}{}^{C}{S}{}_{B]C} = \begin{pmatrix} 0 \\ 2t{}_{[a} p{}_{b]} ,\; t^ap_a \\ \tau p_b+ \nu{}_{ab}  p^a \end{pmatrix},
	\end{equation} 
	so $4{T}{}_{[A}{}^{C}{S}{}_{B]C}\in \operatorname{im}(\Psi{}_{(q,p)}) $, where $\Psi$ is as in Theorem \ref{Thm:PhiPsi1}. In fact, it is clear that elements of the above form span $\operatorname{im}(\Psi{}_{(q,p)} )$, which has dimension $2n$. Then
	\[
		\Psi{}_{(q,p)}^{-1}\left( 4{T}{}_{[A}{}^{C}{S}{}_{B]C} \right) = (t^a, -(\tau p_b+ \nu{}_{ab}  p^a)),
	\] 
	which vanishes if and only if $t=0$, $\tau=0$ and $p^a\nu_{ab}=0$, since $\nu{}_{ab}$ is skew and $p$ is not null.

	In the null case, we can see from \eqref{Eqn:SdotT} that ${T}{}_{[A}{}^{C}{S}{}_{B]C}=0$ if and only if $t{}_{[a} p{}_{b]} = 0$, $t^ap_a=0$ and $\tau p_b+\nu{}_{ab} p^a=0$. The first condition is equivalent to $ t_a = \lambda p_a\implies t^ap_a=0$ since $p$ is null and the second condition is just the second condition in \eqref{Eqn:sl2}.

\end{proof}

\section{Einstein Scales}\label{Sec:Ein}
Projectively flat pseudo-Riemannian manifolds (spaces of constant curvature) have maximal local symmetry: the isotropy group at each point is $O(P,Q)$, and on simply connected spaces the space of Killing vectors has the maximal dimension $n(n+1)/2$. On a conformally flat manifold, projectively flat metrics are equivalent to Einstein metrics. 

Let $(M,\boldsymbol{c})$ be conformally flat and let $\sigma\in \Gamma(\mathcal{E}[1])$ be an almost Einstein scale. These conditions are equivalent to $(M\setminus\mathcal{Z}(\sigma),g)$ being projectively flat. The space of Killing vectors for $g$ is in one to one correspondence with parallel sections of $\Lambda^2(I^{\perp})$, where $I_A = D_A\sigma$ by Theorem \ref{Thm:KT}. If $\tau$ is another almost Einstein scale, then the space of vectors which are Killing for both scales is in one to one correspondence with parallel sections of $\Lambda^2(I^{\perp}\cap J^{\perp})$, where $J_A = D_A \tau$. When working with such a pair of Einstein scales, we will write $\iota = I_AI^A$, $\jmath = J_AJ^A$ and $\kappa = I_AJ^A$. We assume that $I_A$ and $J_A$ are not collinear. Since they are parallel and $M$ is connected, collinearity at any point of $M$ is equivalent with collinearity on all of $M$. Then the rescaling factor $s = \tau^{-1}\sigma$ is regular on an open dense subset of $M$.

\begin{lemma}\label{Lem:Regs}
	Let $I{}_{A}$ and $J_A$ be parallel tractors which are not collinear. Let $\sigma = X^AI_A$ and $\tau = X^AJ_A$ and let $s = \tau\sigma^{-1}\in C^\infty(M\setminus\mathcal{Z}(\sigma))$. Then $ds$ is non-zero on an open dense subset of $M$. In particular, the zero set of $ds$ is contained in the set of points for which
	\begin{equation}\label{Eqn:ijk}
			\iota s(q)^2-2\kappa s(q)+\jmath = 0,
	\end{equation}
	where $\iota = I_AI^A$, $\jmath  = J_AJ^A$ and $\kappa = I_AJ^A$.
\end{lemma}
\begin{proof}
	Writing $\tau = s\sigma$, we have $J_A \overset{\sigma}{=}  (\sigma s,\sigma\nabla_as,-\frac{1}{n}\mathsf{J}^{\sigma}\sigma s-\frac{1}{n}\sigma(\Delta s))$. Suppose that $ds=0$ on an open set $U\subset M$. Then $\Delta s=0$ on $U$, so $J\overset{\sigma}{=} (s\sigma,0,-\frac{1}{n}\mathsf{J}^{\sigma}s\sigma) = sI_A$ on $U$. Then since $I$ and $J$ are parallel, they are collinear on $M$, contradicting our assumption.

	Now let $q\in M\setminus\mathcal{Z}(\sigma)$ be such that $ds|_q=0$. Then in the scale $\sigma$, $Z{}^{A}{}_{a} J{}_{A} =0$ so there is a $\rho\in \mathcal{E}[-1]_q$ such that $J_A|_q = s(q)I_A|_q+\rho X_A|_q$. Then at $q$,
	\[
		\kappa = I^AJ_A = s(q)\iota+\rho\sigma \quad \text{and} \quad \jmath = J^AJ_A = s(q)\kappa +\rho \sigma s(q).
	\] 
	Substituting the first expression for $\rho\sigma$ into the second equation gives \eqref{Eqn:ijk}.
\end{proof}

It turns out that these Killing vectors are shared with a much larger smooth class of metrics.
\begin{proposition}\label{Prop:fKS}
	Let $r_{a_{1}\cdots a_\ell}\in \mathcal{E}_{(a_{1}\cdots a_\ell)}[2\ell]$ and let $\sigma, \tau \in \Gamma(\mathcal{E}[1])$. Let $g = \sigma^{-2}\boldsymbol{g}$ and let $\hat{g} = \tau^{-2}\boldsymbol{g}$ and suppose that $r_{a_{1}\cdots a_\ell}$ is a Killing tensor for both $g$ and $\hat{g}$. Let $s = \tau\sigma^{-1}\in C^\infty(M\setminus\mathcal{Z}(\sigma))$ and let $f$ be a smooth real valued function defined on the range of $s$. Then $r_{a_{1}\cdots a_\ell}$ is a Killing tensor for $\tilde{g} = \tilde{\sigma}^{-2}g$, where $\tilde{\sigma} = f(s)\sigma$, defined on $M\setminus\mathcal{Z}(\sigma\tau f(s))$.
\end{proposition}
\begin{proof}
	Let $\nabla, \hat{\nabla}$ and $\tilde{\nabla}$ be the Levi-Civita connections of $g,\hat g$ and $\tilde g$, respectively. By the conformal transformation laws,
	\[
		\hat{\nabla}_{(a}r_{a_{1}\cdots a_\ell)} = \nabla_{(a}r_{a_{1}\cdots a_\ell)}+\ell\Upsilon^dr_{d(a_{2}\cdots a_\ell}\boldsymbol{g}_{a\,a_{1})},
	\]
	where $\Upsilon_a = s\nabla_a(s^{-1}) = -s^{-1} \nabla_a s$. So, if ${r}{}_{ab}$ is Killing for $g$ and $\hat{g}$, then
	\[
		0 = \hat{\nabla}_{(a}r_{a_{1}\cdots a_\ell)} = \nabla_{(a}r_{a_{1}\cdots a_\ell)}\implies \Upsilon^dr_{d(a_{2}\cdots a_\ell}\boldsymbol{g}_{a\,a_{1})}=0
	\] 
	The conformal transformation laws between $\nabla$ and $\tilde{\nabla}$ use
	\[{\tilde{\Upsilon}}{}_{a} = f(s){\nabla}{}_{a}(f(s)^{-1}) = -\frac{f(s)f'(s)}{f(s)^2}\nabla_a s=s \frac{f'(s)}{f(s)}{\Upsilon}{}_{a},\]
	Then
\[\tilde{\nabla}_{(a}r_{a_{1}\cdots a_\ell)} = \nabla_{(a}r_{a_{1}\cdots a_\ell)}+\ell\tilde{\Upsilon}^dr_{d(a_{2}\cdots a_\ell}\boldsymbol{g}_{a\,a_{1})} = \ell s \frac{f'(s)}{f(s)}\Upsilon^dr_{d(a_{2}\cdots a_\ell}\boldsymbol{g}_{a\,a_{1})}
=0.\]
\end{proof}
Consider the conformal manifold associated to $\mathbb{R}^n$ with the flat metric, which we will discuss in detail in section \ref{Sec:FM}. Let $\sigma\in \mathcal{E}[1]$ represent the Euclidean metric. Then $\tau = -\frac{r^2}{2}\sigma$ is another flat scale, and the space $I^{\perp}\cap J^{\perp}$ is spanned by parallel tractors $I^i_A  = (\sigma q_i,\sigma p_i,0)$. So the Killing tensors shared by $\sigma$ and $\tau$ are the top slots
\[
	2X^A Z{}^{B}{}_{b} I^i_{[A}I^j_{B]} = q_ip_j-q_jp_i.
\]
Then by the above Proposition, the angular momenta functions $L^{i,j}$ are constants of motion for any Hamiltonian of the form $|p|^2+V(r)$, and $2n-3$ of the $L^{i,j}$ are functionally independent. In the rest of this section, we generalise this reasoning to the larger class of Hamiltonians with functions parameterised by an Einstein rescaling function. In addition, we are careful to prove functional independence of $2n-3$ Killing vectors using the tractor framework developed in the previous section, proving quasimaximal superintegrability of the associated families of Hamiltonians.
To prove functional independence, we need two technical lemmas. Let $\mathcal{B}$ be a vector bundle equipped with a connection $\nabla^\mathcal{B}$. A vector subbundle $\mathcal{V}\subset \mathcal{B}$ is \textit{parallel} if $\nabla^\mathcal{B}{}_a s\in \mathcal{E}_a\otimes\mathcal{V}$ for all $s\in\Gamma(\mathcal{V})$.
\begin{lemma}\label{Lem:QM1}
	Let $(M,\boldsymbol{c})$ be a conformal manifold and $\mathcal{V}\subset \mathcal{E}^A$ be a parallel subbundle of the conformal tractor bundle. If there exists an open $U\subset M$ such that $X_q\in \mathcal{V}[1]_q$ for all $q\in U$, then $\mathcal{V}=\mathcal{E}^A$.
\end{lemma}

\begin{proof}
	Let $\nabla$ be the connection of some metric in the conformal class. Then on the open set $U$, 
	\[
		\nabla_a X^A = Z{}^{A}{}_{a} \in \mathcal{E}_a[1]\otimes\mathcal{V}
	\] 
	Then at all points of $U$,
	\[
		\nabla_b Z{}^{A}{}_{a} = -\mathsf{P}{}_{ab} X^A - \boldsymbol{g}{}_{ab} Y^A\in \Gamma(\mathcal{E}{}_{ab}[1]\otimes \mathcal{V} )
	\] 
	Since $X^A\in \mathcal{V}\left[1\right]$, we conclude that $Y^A\in \mathcal{V}[-1]$ at all points of $U$. Hence, at each point $q\in U$, $\mathcal{E}^A_q\subset \mathcal{V}_q $. Since $\mathcal{V}$ is a subbundle and has constant rank, we conclude that $\mathcal{V} = \mathcal{E}^A$.
\end{proof}
\begin{proposition}\label{Prp:Find}
	Let $1< m\le n+1$ and let $I^1,\dots ,I^m\in \Gamma(\mathcal{T^*})$ be linearly independent parallel tractors. Let ${K}^{ij}{}_{AB} = 2{I}{}^{i}{}_{[A}{I}{}^{j}{}_{B]}$ for all $1\le i<j\le m$. Then the $2m-3$ functions induced by the conformal Killing vectors given by the top slots of $K^{1,2},\dots ,K^{1,m}$ and $K^{2,m},\dots K^{m-1,m}$ are functionally independent. Moreover, this is the maximal number of possible functionally independent conformal Killing vectors given by the top slots of tractors in $\Gamma(\operatorname{span}(K^{i,j}))$.
\end{proposition}
\begin{proof}
	By Corollary \ref{Thm:Find}, it suffices to show that, for all $(q,p) $ in an open dense subset of $T^*M$, the tractors
	\begin{equation}\label{Eqn:Ind:0}
		{K}^{1,2}{}_{[A}{}^{C}{S}{}_{B]C},\dots ,{K}^{1,m}{}_{[A}{}^{C}{S}{}_{B]C},\dots ,{K}^{m-1,m}{}_{[A}{}^{C}{S}{}_{B]C},
	\end{equation}
	are linearly independent, where ${S}{}_{AB} = \mathbb{X}{}_{AB}{}^{b} {p}{}_{b}$. So suppose that $(q,p)\in T^*M$ and that ${\lambda}{}_{1,2},\dots ,{\lambda}{}_{1,m},\dots {\lambda}{}_{m-1,m}\in \mathbb{R}$ are such that
	\begin{equation}
		{\lambda}{}_{1,2}{K}^{1,2}{}_{[A}{}^{C}{S}{}_{B]C}+\dots +{\lambda}{}_{1,m}{K}^{1,m}{}_{[A}{}^{C}{S}{}_{B]C}+\dots +{\lambda}{}_{m-1,m}{K}^{m-1,m}{}_{[A}{}^{C}{S}{}_{B]C}=0.
	\end{equation}
	Let ${K} = {\lambda}{}_{1,2}{K}^{1,2}+\dots +{\lambda}{}_{1,m}{K}^{1,m}+\dots +{\lambda}{}_{m-1,m}{K}^{m-1,m}$. Now suppose that $p$ is not null. Then, by Lemma \ref{Lem:sl},
	\begin{equation}\label{Eqn:Ind:1}
		X^A K{}_{AB} =0,\quad \text{and} \quad {p}{}^{a}\mathbb{Z}{}^{AB}{}_{ab} {K}{}_{AB}=0
	\end{equation}
	Expanding the first condition gives
	\begin{align}\label{Eqn:Ind:1.5}
	0=X{}^{A}{K}{}_{AB} 
		&= (\lambda{}_{1,2} \sigma_{1}-\lambda{}_{2,m} \sigma_m) I{}^{2}{}_{B} +\cdots + \left( \lambda{}_{1,m-1} \sigma_{1}-\lambda{}_{m-1,m} \sigma{}_{m}  \right) I{}^{m-1}{}_{B}\notag \\ 
		&\mathrel{\phantom{=}} -\left( \sum_{i=2}^{m} \lambda{}_{1,i} \sigma_i \right) I{}^{1}{}_{B} +\left( \sum_{i=1}^{m-1} \lambda{}_{i,m} \sigma_i \right) I{}^{m}{}_{B} .
	\end{align} 
	Since the tractors $I_{1},\dots ,I{}^{m} $ are linearly independent, it follows that
	\begin{equation}\label{Eqn:Ind:2}
	\lambda{}_{1,i} \sigma_{1}-\lambda{}_{i,m} \sigma_m=0
	\end{equation} 
	for each $i\in \{2,\dots ,m-1\} $. We then calculate
	\begin{align*}
		0&=2I{}^{1\,B} K{}_{[A}{}^{C} S{}_{B]C}  \\ 
		 &= -I{}^{1\,B} K{}_{B}{}^{C} S{}_{AC} +\lambda{}_{1,2} \left( I{}^{1}{}_{B} I{}^{2}{}_{C} S{}^{BC}\right)I{}^{1}{}_{A}  +\cdots +\lambda{}_{1,m} \left( I{}^{1}{}_{B} I{}^{m}{}_{C} S{}^{BC}\right)I{}^{1}{}_{A} \\ 
		 &\mathrel{\phantom{=}}+ \lambda{}_{2,m} \left( I{}^{1}{}_{B} I{}^{m}{}_{C} S{}^{BC}  \right)I{}^{2}{}_{A}  +\cdots + \lambda{}_{m-1,m} \left( I{}^{1}{}_{B} I{}^{m}{}_{C} S{}^{BC}  \right) I{}^{m-1}{}_{A} \\
		 &\mathrel{\phantom{=}}-\lambda{}_{2,m} \left( I{}^{1}{}_{B} I{}^{2}{}_{C} S{}^{BC} \right) I{}^{m}{}_{A} -\cdots -\lambda{}_{m-1,m} \left( I{}^{1}{}_{B} I{}^{m-1}{}_{C} S{}^{BC}  \right) I{}^{m}{}_{A} 
	\end{align*}
	Since $X{}^{C} K{}_{AC} =0$, the first term is a multiple of $X_A$. Since $m\le n+1$, the parallel subbundle spanned by $I^1,\dots ,I^m$ does not contain $X$ for all $q$ in an open dense subset of $M$ by Lemma \ref{Lem:QM1}. Given such a $q$, we conclude that 
	\[
	\lambda{}_{i,m} I{}^{m}{}_{A}  I{}^{1}{}_{C} S{}^{AC} =0
	\] 
	for $i\in \{2,\dots ,m-1\}$. We have $2I{}^{m}{}_{A}  I{}^{1}{}_{C} S{}^{AC} =-K{}^{1,m}{}_{AC}S{}^{AC}  = p^c k^{1,m}_c$. Since $k{}^{1,m}{}_{c} $ is the top slot of a non-zero parallel tractor, it is non-zero on an open dense subset of $M$. This follows from Lemma \ref{Lem:BGG}, but can also quickly be derived from \eqref{Eqn:XZWY}. Suppose that $\mathbb{X}{}^{AB}{}_{b} K{}^{1,m}{}_{AB} =0$ on an open subset of $M$. Then
	\[
	0=\nabla_a(\mathbb{X}{}^{AB}{}_{b} K{}^{1,m}{}_{AB}) = \mathbb{Z}{}^{AB}{}_{ab} K{}^{1,m}{}_{AB} -\boldsymbol{g}{}_{ab} \mathbb{W}{}^{AB} K{}^{1,m}{}_{AB}.
	\] 
	So both $\mathbb{Z}{}^{AB}{}_{ab} K{}^{1,m}{}_{AB}$ and $\mathbb{W}{}^{AB} K{}^{1,m}{}_{AB}$ vanish. Applying the connection again to either of these gives $\mathbb{Y}{}^{AB}{}_b K{}^{1,m}{}_{AB}=0$.

	Hence, if $(q,p)\in T^*M$ are contained in the open dense subset of $T^*M$ satisfying
	\[
		p_ap^a\neq 0,\quad \sigma_1(q)\neq 0,\quad \sigma_m(q)\neq 0,\quad \text{and} \quad I{}^{m}{}_{A}  I{}^{1}{}_{C} S{}^{AC} \neq 0,
	\]
	then $\lambda_{i,m}=0$ for all $i\in \{2,\dots ,m-1\} $. Then since $\sigma_{1}\neq 0$, we have $\lambda_{1,i}=0$ for all $i\in \{2,\dots ,m-1\} $ by \eqref{Eqn:Ind:2}. This reduces \eqref{Eqn:Ind:1.5} to
	\[
	0=X^A K{}_{AB}  = \lambda_{1,m}(\sigma_1I{}^{m}{}_{B} - \sigma_m I{}^{1}{}_{B}).
	\] 
	Then by linear independence and $\sigma_{m}\neq 0$, we conclude that $\lambda{}_{1,m} =0$ and the tractors in \eqref{Eqn:Ind:0} are linearly independent on an open dense subset of $T^*M$.

	To show that $2m-3$ is maximal, we fix $q\in M$ such that $\sigma_{1}(q)\neq 0$, so also $X_q\not\in \operatorname{span}(I^1|_q,\dots ,I^m|_q)^{\perp}$. Consider the linear map $\phi_{p,q}\colon \operatorname{span}(K{}^{i,j}|_q)\to \mathcal{E}_A|_q\oplus \mathcal{E}_a|_q$ given by
	\[
	\phi_{p,q}({K}{}_{AB}) = (X{}^{A}{K}{}_{AB} ,p^b\mathbb{Z}^{AB}{}_{ab} {K}{}_{AB}),
	\]
	calculating in the scale $\sigma_{1}$. Then $X{}^{A} K{}_{AB} \in \operatorname{span}(I^{1},\dots ,I^m)\cap X^{\perp}$, which has dimension $m-1$ at $q$ by assumption. Then $p^b\mathbb{Z}^{AB}{}_{ab} {K}{}_{AB}\in \operatorname{span}(Z{}^{A}{}_{a} I{}^{2}{}_{A} ,\dots ,I{}^{m}{}_{A} Z{}^{A}{}_{a} )\cap p^{\perp}$, which at a fixed $q\in M$ has dimension $m-2$ for $p$ outside of a lower dimensional subspace of $T_q^*M$. Hence, the range of $\phi_{p,q} $ has dimension at most $2m-3$ for all $(q,p)$ in an open dense subset of $T^*M$.
\end{proof}
The tractors $S{}_{AB} = \mathbb{X}{}_{AB}{}^{b} p{}_{b} $ span an invariant subbundle of the adjoint tractor bundle. The product, $R{}_{C[A} S{}^{C}{}_{B]} $ is the same as the \textit{action} of the adjoint tractor $R$ on $S$. The adjoint tractor bundle is an associated bundle, $\mathcal{G}\times_P\mathfrak{g}$, where $\mathfrak{g} = \mathfrak{so}(p+1,q+1)$, and the action is induced by the action of $\mathfrak{g}$ on $\mathfrak{g}$, given by the Lie bracket. The product of two adjoint tractors, given by $R{}_{AB} T{}^{AB} $ is induced by the Killing form $B$ on $\mathfrak{g}$, for more details see \cite{cap_tractor_2000}. In particular it gives the, otherwise easy to show identity
\begin{equation}
	R{}_{C[B} S{}^{C}{}_{A]} T{}^{AB}  = -T{}_{C[B} S{}^{C}{}_{A]} R{}^{AB} 
\end{equation} 
for any adjoint tractors $R,S,T$, since $B([r,s],t) = B(r,[s,t]) = -B([t,s],r)$.

\begin{lemma}\label{Lem:gind}
	Let $g$ be a metric and let $k^1,\dots ,k^m$ be vector fields contained in a proper subbundle of $TM$. Then the functions induced by $g^{-1},k^1,\dots ,k^m$ are functionally independent if and only if the functions induced by $k^1,\dots ,k^m$ are functionally independent.
\end{lemma}
\begin{proof}
	Let $\nabla$ be the Levi-Civita connection of $g$. Then, calculating in the splitting $T^*(T^*M) \overset{\nabla}{=} T^*M\oplus TM$,
	\[
	dg^{-1}(q,p) \overset{\nabla}{=}  (2p_a,0)\quad \text{and} \quad dk^i \overset{\nabla}{=} (k_a^i,p^b{\nabla}{}_{a}{k}{}_{b}{}^{i}).
	\] 
	Suppose that $(q,p)\in T^*M$ is such that $p\not\in \operatorname{span}(k^{1}|_q,\dots ,k^m|_q)$ and that $\lambda,\lambda_{1},\dots ,\lambda_m\in \mathbb{R}$ are such that 
	\[
	\lambda dg^{-1}+\lambda_{1}dk^1+ \dots + \lambda_m dk^m=0 \implies \lambda p+\lambda_{1}k^1+\dots +\lambda_mk^m=0.
	\] 
	Then $\lambda=0$. Since $k^1,\dots ,k^m$ are contained in a proper subbundle of $TM$, there is an open dense subset of $T^*M$ such that $p\not\in \operatorname{span}(k^{1}|_q,\dots ,k^m|_q)$. Hence, linear independence of $dg^{-1},dk^1, \dots, dk^m$ on an open dense subset is equivalent to linear independence of $dk^1, \dots, dk^m$ on an open dense subset.
\end{proof}

We are now ready to state the main result of this section. In our assumptions, we ask for our manifold to be conformally flat and simply connected. This is to ensure that we can use parallel transport to extend any tractor at a point to a parallel tractor on all of $M$. Hence, one can replace ``Let $(M,\boldsymbol{c}) $ be conformally flat and simply connected'' by ``Let $(M,\boldsymbol{c})$ have trivial conformal holonomy''.

\begin{theorem}\label{Thm:QMS}
	Let $(M,\boldsymbol{c})$ be conformally flat and simply connected. Let ${I}{}_{A} = D_A \sigma$ and $J_A = D_A \tau$ for $\sigma,\tau\in \Gamma(\mathcal{E}[1])$ be linearly independent parallel tractors. Let $s = \tau\sigma^{-1}$ and let $V\colon \mathbb{R}\setminus A\to \mathbb{R} $ be a smooth function, where $A\subset \mathbb{R}$ is finite. Then the Hamiltonian
	\[
	\mathcal{H}(q,p) = g^{-1}_q(p,p)+V\left( \frac{\tau}{\sigma} \right) ,
	\] 
	where $g = \sigma^{-2}\boldsymbol{g}$, has $2n-2$ functionally independent constants of motion. The constants of motion are $f_{0} = \mathcal{H}$, and $f_{1},\dots ,f_{2n-3}$ induced by vectors fields which are Killing for $g$ and $\hat{g} = \tau^{-2}\boldsymbol{g}$.
\end{theorem}

\begin{proof}
	Let $I^1,\dots ,I^n\in \Gamma(\mathcal{E}_A)$ be parallel tractors spanning the $n$-dimensional subbundle $\operatorname{span}(I,J)^{\perp}\subset \mathcal{E}_A$. Let $K^{i,j}_{AB} = 2{I}{}^{i}{}_{[A}{I}{}^{j}{}_{B]}$ and $k^{i,j}_b = \mathbb{X}{}^{AB}{}_{b} {K}^{i,j}_{AB}$ for each $1\le i<j\le n$. Then ${I}{}^{A}{K}^{i,j}_{AB}=0=J^A{K}^{i,j}_{AB}$ and so by Theorem \ref{Thm:KT}, all the $k^{i,j}_b$ are Killing for both $g$ and $\hat{g}$. Then by Proposition \ref{Prop:fKS}, $k^{i,j}_b$ is also Killing for $\tilde{g} = V\left(\frac{\tau}{\sigma}\right)g$. Then by Proposition \ref{Prp:PBKT}, the function induced by ${k}{}^{i,j}{}_{b}$ is a constant of motion for $\mathcal{H}$.  

	By Proposition \ref{Prp:Find}, the $2n-3$ functions induced by $k^{1,2},\dots ,k^{1,n},\dots ,k^{n-1,n}$, are functionally independent. Observe that also $k^{i,j}\in \operatorname{span}(\mu^2,\dots ,\mu^n)$ for all $1\le i<j\le n$, where ${\mu}{}^{i}{}_{a} = {Z^1}{}^{A}{}_{a}{I}{}^{i}{}_{A}$ and $Z^1$ is the projector in the scale $\sigma_{1} = {X}{}^{A}{I}^{1}_{A}$. Hence, by Lemma \ref{Lem:gind}, the functions induced by $g^{-1},k^{1,2},\dots ,k^{1,n},\dots ,k^{n-1,n}$ are functionally independent.
	\end{proof}

Take $I_A$ and $J_A$ as in the above theorem. By Lemma \ref{Lem:gind}, The differentials $dk^{i,j}$ and $dg^{-1}$ span a $2n-2$ dimensional subspace of $T^*(T^*M)$, in the following lemma, we show explicitly which subspace of adjoint tractors they correspond to and give criteria which will be used in Section \ref{Sec:Ex} to demonstrate functional independence of the last constant of motion.
\begin{lemma}\label{Lem:Fin}
	Let $(M,\boldsymbol{c})$ have trivial conformal holonomy and let $I_A = D_A\sigma$ and $J_A = D_A\tau$ be parallel tractors which are not collinear. Let $R_{ABCD}$ be a parallel tractor with Riemann curvature symmetries which is orthogonal to $I$ and let $r{}_{bd} =\mathbb{X}{}^{AB}{}_{b} \mathbb{X}{}^{CD}{}_{d}  R{}_{ABCD}$. Then the following are equivalent
	\begin{tabel}
		\item\label{Lem:Fin:1} $r{}_{bd} $ is Killing in the scale $\tau = X^AJ_A$,
		\item\label{Lem:Fin:2} $r_{bd}\nabla^b\tau  =0$ on $M\setminus\mathcal Z(\sigma)$, where $\nabla=\nabla^\sigma$,
		\item\label{Lem:Fin:3} $J{}^{D} R{}_{ABCD} =0$ and 
		\item\label{Lem:Fin:4} the function induced by $r{}_{bd} $ is not functionally independent from the functions induced by $k_b = \mathbb{X}{}^{AB}{}_{b} K{}_{AB}$ for all parallel $K{}_{AB} \in \Gamma(\mathcal{E}_{[AB]})$ such that $I{}^{A} K{}_{AB} =J{}^{A} K{}_{AB} =0$.
	\end{tabel}
\end{lemma}
\begin{proof}
	Let $\nabla$ and $\hat{\nabla}$ be the Levi-Civita connections of $g = \sigma^{-2}\boldsymbol{g}$ and $\hat{g} = \tau^{-2}\boldsymbol{g}$, respectively, where they are defined. Let $\Upsilon_a = -\tau^{-1}\nabla{}_{a} \tau$ on $M\setminus(\mathcal{Z}(\tau)\cup \mathcal{Z}(\sigma))$. This is the same $\Upsilon$ from \eqref{Eqn:CT}. The conformal transformation law for $\nabla{}_{a} r{}_{bd} $ is 
	\[
		\hat{\nabla}_{a} r{}_{bd}  -\nabla{}_{a}  r{}_{bd} =2\Upsilon{}_{a} r{}_{bd} -\Upsilon{}_{b} r{}_{ad} -\Upsilon{}_{d} r{}_{ba} + \Upsilon{}^{c} r{}_{cd} \boldsymbol{g}{}_{ab} +\Upsilon{}^{c} r{}_{bc} \boldsymbol{g}{}_{ad}
	\] 
	so $\hat{\nabla}_{(a} r{}_{bd)}  -\nabla{}_{(a}  r{}_{bd)}=\Upsilon^c r{}_{c(d} \boldsymbol{g}{}_{ab)}$, which vanishes if and only if $\Upsilon{}^{b} r{}_{bd} =0$. Since $\tau\Upsilon_a = \nabla_a\tau$ on a dense open subset of $M$, this proves \ref{Lem:Fin:1}$\iff $\ref{Lem:Fin:2}. 

	Proposition \ref{Prp:Run} gives \ref{Lem:Fin:1} $\implies$ \ref{Lem:Fin:3}, and Theorem \ref{Thm:KT} gives the converse.
	
	Suppose that $J^A R{}_{ABCD} =0$. Then $S{}^{AB} R{}_{ABCD} \in \Gamma(\Lambda^2(I^{\perp}\cap J^{\perp}))$. So for all $(q,p)\in T^*M$, $S{}^{AB} R{}_{ABE[D} S{}^{E}{}_{C]} $ is always in the subspace spanned by $K{}_{E[D} S{}^{E}{}_{C]} $ for all parallel $K\in \Lambda^2(I^{\perp}\cap J^{\perp})$. Hence, $r{}_{bd}$ is not independent from the functions induced by the corresponding conformal Killing vectors by Theorems \ref{Thm:PhiPsi1} and \ref{Thm:PhiPsi2}.

	Now suppose that \ref{Lem:Fin:4} holds. Let $I^1,\dots ,I^n$ be parallel tractors forming a basis for $I^{\perp}\cap J^{\perp}$ on all of $M$. Let $K^1,\dots ,K{}^{2n-3} \in \Lambda^2(I^{\perp}\cap J^{\perp})$ be the prolongations of the maximal set of functionally independent Killing vectors constructed in Proposition \ref{Prp:Find}. Then
	\begin{equation}\label{Eqn:KJSJ}
		2K^i{}_{E[D} S{}^{E}{}_{C]} J{}^{D} (SJ)^{C} = -K^i{}_{EC} (SJ)^E (SJ)^C= 0
	\end{equation}
	for all $i\le 2n-3$.
	By the contrapositive of Corollary \ref{Cor:Find}, there is an open subset of $\hat{U}\subset T^*M$ on which $S{}^{AB} R{}_{ABE[D} S{}^{E}{}_{C]}$ lies in the span of $K^i{}_{E[D} S{}^{E}{}_{C]}$, where $S{}_{AB}  = \mathbb{X}{}_{AB}{}^{b}  p_b$ and $i\le 2n-3$. We then restrict to the non-empty open subset $U\subset \hat{U}$ consisting of $(q,p)\in \hat{U}$ such that $p_ap^a\neq 0$, $\tau(q)\neq 0$. Then
	\begin{align}\label{Eqn:SRJJ}
		8S{}^{AB} R{}_{ABE[D} S{}^{E}{}_{C]} J{}^{D} S{}^{C}{}_{A'} J{}^{A'} 
		&= 4S{}^{AB} R{}_{ABED} J^D S{}^{E}{}_{C} S{}^{C}{}_{A'}J{}^{A'} - 4S{}^{AB} R{}_{ABEC} (SJ){}^{E} (SJ)^{C}  \notag\\
		&= -p_ap^a \tau S{}^{AB} R{}_{ABCD} X^C J^D \\
		&= -p_ap^a \tau p^b r{}_{bd}\nabla^d \tau =0.\notag
	\end{align} 
	So $p^b r{}_{bd}\nabla^d \tau=0$ at $(q,p)$. Since $U$ is open, the intersection $U\cap T_q^*M$ is a spanning set for $T_q^*M$. So, by linearity, $p^br{}_{bd}\nabla^d\tau  = 0$ for all $p\in T^*_qM$ . So $r{}_{bd} \nabla^b \tau=0$ on $\pi(U)$, which is an open subset of $M$. Let $V \subset  \pi(U)$ be simply connected and apply \ref{Lem:Fin:2}$\implies$\ref{Lem:Fin:3} to $r$ as a Killing tensor on $(V,\boldsymbol{c}|_V)$. From which we conclude that $J{}^{D} R{}_{ABCD} =0$ on $V$. Since $J^D R{}_{ABCD} $ is parallel and $M$ is connected, it follows that $R{}_{ABCD} J^D=0$ on $M$. So we have shown that \ref{Lem:Fin:3} holds.\end{proof}

\begin{lemma}\label{Lem:Indfromg}
Let $(M,\boldsymbol{c})$ have trivial conformal holonomy and let $I_A = D_A\sigma$ and $J_A = D_A\tau$ be parallel tractors which are not collinear. Let $\iota = I_AI^A$ and $\kappa = I_AJ^A$ and suppose that $(\iota,\kappa)\neq (0,0)$. Let $R_{ABCD}$ be a parallel tractor with Riemann curvature symmetries which is orthogonal to $I$ and let $r{}_{bd} =X{}^{A} Z{}^{B}{}_{b} X{}^{C} Z{}^{D}{}_{d} R{}_{ABCD} $. Then the following are equivalent,
	\begin{tabel}
		\item\label{Lem:Fin:2:1} there exists $\lambda\in \mathbb{R}$ such that $r{}_{bd} -\lambda \sigma^2 \boldsymbol{g}_{bd}$ is Killing in the scale $\tau=X^AJ_A$,
		\item\label{Lem:Fin:2:2} there exists $\lambda\in \mathbb{R}$ such that $ r{}_{bd}\nabla^b \tau  = \lambda\sigma^2\nabla_d\tau$ on $M\setminus\mathcal Z(\sigma)$, where $\nabla=\nabla^\sigma$,
		\item\label{Lem:Fin:2:3} there is a constant $\lambda\in \mathbb{R}$ such that $J{}^{A} (R{}_{ABCD} - \lambda G{}_{ABCD} )=0$, where $G{}_{ABCD} = (h\owedge (I^2-\frac{\iota}{2}h))_{ABCD}$,
		\item\label{Lem:Fin:2:3b} there is a parallel $S{}_{ABCD} \in I^{\perp}\cap J^{\perp}$ with the same trace-free part as $R$ and
		\item\label{Lem:Fin:2:4} the function induced by $r{}_{bd} $ is not functionally independent from the functions induced by $\sigma^2\boldsymbol{g}$ and $k_b = X{}^{A} Z{}^{B}{}_{b}K{}_{AB}$ for all parallel $K{}_{AB} \in \Gamma(\mathcal{E}_{[AB]})$ such that $I{}^{A} K{}_{AB} =J{}^{A} K{}_{AB} =0$.
	\end{tabel}
\end{lemma}
\begin{proof}
	Observe that $G{}_{ABCD} $ is parallel, $I^A G{}_{ABCD} =0$ and $X{}^{A} Z{}^{B}{}_{b} X{}^{C} Z{}^{D}{}_{d} G{}_{ABCD}  = \sigma^2\boldsymbol{g}$. The equivalences \ref{Lem:Fin:2:1} $\iff $ \ref{Lem:Fin:2:2} $\iff $ \ref{Lem:Fin:2:3} follow from the same correspondences in Lemma \ref{Lem:Fin}.

	The implication \ref{Lem:Fin:2:3} $\implies$ \ref{Lem:Fin:2:3b} is clear. The converse follows from Lemma \ref{Lem:GABCD}, which states that $G{}_{ABCD} $ spans the space of curvature tractors which are pure trace and which are orthogonal to $I$, which then implies that $S-R = \lambda G$.
	
	The implication \ref{Lem:Fin:2:1} $\implies$ \ref{Lem:Fin:2:4} follows from \ref{Lem:Fin:1} $\implies$ \ref{Lem:Fin:4} of Lemma \ref{Lem:Fin}

	Let $K{}^{i,j}{}_{AB} = 2I{}^{i}{}_{[A} I{}^{j}{}_{B]} $ for each $1\le i<j\le n$ be a spanning set of parallel sections for $\Lambda^2(I^{\perp}\cap J^{\perp})$, where $I^1,\dots ,I^n$ is a parallel basis for $I^{\perp}\cap J^{\perp}$. Then let $k{}^{i,j}{}_{b} =  \mathbb{X}{}^{AB}{}_{b} K{}^{i,j}{}_{AB} $. 

	Suppose \ref{Lem:Fin:2:4} holds, that is, $r$ is functionally dependent on $k{}^{i,j} $ and $\sigma{}^{2}\boldsymbol{g}$. Then, by the contrapositive of Corollary \ref{Cor:Find}, there are functions $\lambda^{i,j}$ and $\lambda$ defined on an open subset $\hat{U}\subset T^*M$ such that
	\begin{equation}\label{Eqn:SRSKG}
		S{}^{AB}R{}_{ABE[D} S{}^{E}{}_{C]} = \sum_{i,j} \lambda{}^{i,j} (q,p)K{}_{E[D} S{}^{E}{}_{C]} +\lambda(q,p)S{}^{AB} G{}_{ABE[D} S{}^{E}{}_{C]},
	\end{equation} 
	where $S{}_{AB}  = \mathbb{X}{}_{AB}{}^{b} p{}_{b} $. We then restrict to the subset $U\subset \hat{U}$ consisting of $(q,p)\in \hat{U}$ such that $\sigma(q)\neq 0$, $\tau(q)\neq 0$, $\kappa\sigma-\iota\tau\neq 0$ and $\nabla_a\tau\neq 0$. Then $U$ is non-empty as $\sigma,\tau$ and $\kappa\sigma-\iota\tau$ are all scales, and $\nabla_a \tau = \sigma (ds)_a$ is non-zero on an open dense subset of $M$ by Lemma \ref{Lem:Regs}. Then contracting $J^C(SJ)^D$ into \eqref{Eqn:SRSKG} gives  
	\[
		\begin{gathered}
			8 S{}^{AB}R{}_{ABE[D} S{}^{E}{}_{C]} J^D (SJ)^C = 8\lambda(q,p)S{}^{AB} G{}_{ABE[D} S{}^{E}{}_{C]} J^D (SJ)^C,\\
			\implies-\tau p_ap^a p^b X^A Z{}^{B}{}_{b} X^C J^D R{}_{ABCD} = \lambda(q,p) \tau p_ap^ap^b X^AZ^B_b X^CJ^D G{}_{ABCD} 
		\end{gathered}
	\] 
	with the terms from $K^{i,j}$ vanishing by \eqref{Eqn:KJSJ}. Expanding the above in the scale $\sigma$ gives
	\[
		-p_ap^a \tau p^b r{}_{bd}\nabla^d \tau = -\lambda(q,p)p_ap^a\tau\sigma^2p^b \boldsymbol{g}{}_{bd}\nabla^d \tau\implies p^b r{}_{bd}\nabla^d \tau  = \lambda(q,p)\sigma^2p^d\nabla_d\tau.
	\] 
	We then contract $I^CJ^D$ into \eqref{Eqn:SRSKG}. Since the $K^{i,j}$ are orthogonal to $I,J$, this gives
	\[
		4S{}^{AB} R{}_{ABED}  S{}^{E}{}_{C} I^CJ^D = 4\lambda(q,p)S{}^{AB} G{}_{ABED}  S{}^{E}{}_{C} I^CJ^D.
	\] 
	We calculate in the scale $\sigma$, which gives $S{}^{E}{}_{C} I^C = -\frac{1}{2}\sigma Z{}^{E}{}_{e} p^e$, so
	\[
		\begin{aligned}
			\sigma p^bp^c X^A Z{}^{B}{}_{b} Z{}^{C}{}_{c} J^D R{}_{ABCD} 
		&= \sigma\lambda(q,p) p^bp^c X^A Z{}^{B}{}_{b} Z{}^{C}{}_{c} J^D G{}_{ABCD} \\
		&= -\lambda(q,p)\sigma\left( \sigma\kappa-\iota\tau \right) p_bp^b ,
		\end{aligned}
	\] 
	which we write as $p^bp^c B{}_{bc} = -\lambda(q,p) ( \sigma\kappa-\iota\tau ) p_cp^c$, with $B{}_{bc} = X^A Z{}^{B}{}_{(b} Z{}^{C}{}_{c)} J^D R{}_{ABCD}$. 

	Since $(\sigma\kappa-\iota\tau)	 p_cp^c\neq 0$ at $(q,p)$, we can eliminate $\lambda$ from the two identities we have derived, giving
	\begin{equation}\label{Eqn:Poly}
		-(\sigma\kappa-\iota\tau)p_cp^c p^b r{}_{bd} \nabla^d \tau = \sigma^2p^d\nabla_d\tau p^bp^c B{}_{bc}.
	\end{equation} 
	By varying over $p\in T_q^*M\cap U$, the above holds on an open subset of $T^*_qM$. Two polynomials are equal on an open subset of $\mathbb{R}^n$ if and only if they are equal on all of $\mathbb{R}^n$. Hence, \eqref{Eqn:Poly} holds on $T_q^*M$.

	If $p\in T_q^*M$ is such that $p^d\nabla_d \tau=0$ and $p_cp^c\neq 0$, then $p^b r{}_{bd} \nabla^d \tau=0$ by \eqref{Eqn:Poly}. Since $(\nabla_d \tau|_q)^{\perp}$ has dimension $n-1>n/2$, it cannot be a totally null subspace. Therefore, the set of non-null $p$ which are orthogonal to $\nabla_d \tau$ is dense in $(\nabla_d \tau|_q)^{\perp}$. So $p^b r{}_{bd} \nabla^d \tau=0$ for all $p\in(\nabla_d \tau|_q)^{\perp}$. From which it follows that there is a $c_q\in \mathbb{R}$ such that $r{}_{bd} \nabla{}^{d} \tau = c_q\sigma^2\nabla_b\tau$ at $q$. Varying over $q\in V:=\pi(U)$, this gives that $\lambda(q,p) = c_q$ for all $(q,p)\in U$. So, writing $\lambda(q) = c_q$, we have
	\[
		\begin{gathered}
			p^b X^A Z{}^{B}{}_{b} X^CJ^D R{}_{ABCD}  = \lambda(q)p^b X^A Z{}^{B}{}_{b}X^C J^D G{}_{ABCD}  \quad \text{and} \\
			p^bp^c X^A Z{}^{B}{}_{b} Z{}^{C}{}_{c} J^D R{}_{ABCD} = \lambda(q) p^bp^c X^A Z{}^{B}{}_{b} Z{}^{C}{}_{c} J^D G{}_{ABCD}
		\end{gathered}	
	\] 
	for all $(q,p)\in U$. For each fixed $q$, the above are polynomials defined on an open subset of $T_q^*M$. Hence, the identities hold for all $q\in \pi(U)$ and $p\in T_q^*M$. That is,
	\[
		\begin{gathered}
			 X^A Z{}^{B}{}_{b} X^CJ^D R{}_{ABCD}  = \lambda  X^A Z{}^{B}{}_{b}X^C J^D G{}_{ABCD}  \quad \text{and} \\
			 X^A Z{}^{B}{}_{(b} Z{}^{C}{}_{c)} J^D R{}_{ABCD} = \lambda X^A Z{}^{B}{}_{(b} Z{}^{C}{}_{c)} J^D G{}_{ABCD}
		\end{gathered}	
	\] 
	on $V = \pi(U)\subset M$ open. So, on $V$
	\[
		\begin{aligned}
			\nabla_{(c}\big(X^AZ^B{}_{b)}X^CJ^DR_{ABCD}\big)
			&= X^AZ^B{}_{(b}Z^C{}_{c)}J^DR_{ABCD} - \boldsymbol{g}_{bc}\,X^AY^BX^CJ^DR_{ABCD}\\
			&= X^AZ^B{}_{(b}Z^C{}_{c)}J^DR_{ABCD},
		\end{aligned}
	\]
	since $R\in I^{\perp}$, and so the same holds with $G$ in the place of $R$. Therefore,
	\[
	\begin{gathered}
		X^A Z{}^{B}{}_{(b} Z{}^{C}{}_{c)} J^D R{}_{ABCD} = \lambda  X^A Z{}^{B}{}_{(b}Z{}^{C}{}_{c)}  J^D G{}_{ABCD} + (\nabla{}_{(c} \lambda)X^A Z{}^{B}{}_{b)}X^C J^D G{}_{ABCD}\\
		\implies 0 = (\nabla{}_{(c} \lambda)X^A Z{}^{B}{}_{b)}X^C J^D G{}_{ABCD} = (\nabla{}_{(c} \lambda) (\nabla{}_{b)} \tau)
	\end{gathered} 
\]
Since $\nabla_b\tau\neq 0$ throughout $V$, we conclude that $\nabla_b\lambda=0$ on $V$. So, by restricting to a connected open subset if necessary, $\lambda$ is constant on $V$. So
\[
	\nabla^b\tau r{}_{bd}  = \lambda \sigma^2 \nabla_d\tau
\] 
on $V$, with $\lambda$ constant. So, we then apply \ref{Lem:Fin:2:2}$\implies$\ref{Lem:Fin:2:3} to $r{}_{bd}|_{V'} $ and $(V',\boldsymbol{c}|_{V'})$, where $V'$ is a simply connected subset of $V$, so that it has trivial conformal holonomy, to see that $J^A (R{}_{ABCD} -\lambda G{}_{ABCD}) =0$ on $V'$, whence $J^A (R{}_{ABCD} -\lambda G{}_{ABCD}) =0$, since it is parallel and $M$ is connected. Completing the proof of \ref{Lem:Fin:2:4}$\implies$\ref{Lem:Fin:2:3}.

\end{proof}

\begin{remark}
The hypothesis $(\iota,\kappa)\neq(0,0)$ cannot be removed. On $\mathbb{R}^n$, $n\ge3$, let $I$ be the Euclidean scale tractor and $J=I^1$, so that $\tau = q_1\sigma$, $\iota=\kappa=0$ and $\tau^{-2}\boldsymbol{g}$ is the upper half space metric. The Euclidean Killing tensor $r=p_1p_2$ has full prolongation orthogonal to $I$, and on $\{p_1\neq0\}$
\[
	d(p_1p_2) = \frac{p_2}{2p_1}\,d|p|^2 + \Big(p_1-\frac{p_2^2}{p_1}\Big)dp_2 - \frac{p_2}{p_1}\sum_{i\ge3}p_i\,dp_i,
\]
so $r$ is functionally dependent on $|p|^2$ and on the translations $p_2,\dots,p_n$, which are Killing for both scales. However $r_{bd}\nabla^b\tau\propto dq_2$ while $\sigma^2\nabla_d\tau\propto dq_1$, so \ref{Lem:Fin:2:2} fails. In this case $\kappa\sigma-\iota\tau\equiv0$, so the contraction with $I^CJ^D$ in the proof above carries no information about $\lambda$.
\end{remark}
A Hamiltonian $\mathcal{H}$ is \textit{locally second-order maximally superintegrable} if and only if for every $q\in M$ at which $\mathcal{H}$ is defined, there exists an open neighbourhood $U$ of $q$ such that $\mathcal{H}|_{T^*U}$ is second-order maximally superintegrable. 

Proposition \ref{Prp:Find} also shows that $2n-3$ is the maximal number of functionally independent linear constants of motion for the Hamiltonian of Theorem \ref{Thm:QMS} when $V$ is non-constant. (Indeed, a Killing vector for $g$ and $V(s)g$ is also Killing for $s^{-2}g = \tau^{-2}\boldsymbol g$ by Proposition \ref{Prop:fKS}.) This results in the following requirement for local maximal quadratic superintegrability,
\begin{theorem}\label{Thm:2Ein}
	Let $(M,\boldsymbol{c})$ have trivial conformal holonomy and let $I_A = D_A\sigma$ and $J_A = D_A\tau$ be parallel tractors which are not collinear. Let $\iota = I_AI^A$ and $\kappa = I_AJ^A$ and suppose that $(\iota,\kappa)\neq (0,0)$. Let $s = \tau\sigma^{-1}\in C^\infty (M\setminus\mathcal{Z}(\sigma))$, let $f\colon \mathbb{R}\to \mathbb{R}$ be smooth, let $\tilde{\sigma} = f(s)\sigma$ and let $\tilde{I}_A = D_A \tilde{\sigma}$. Let $V(s) = f(s)^{-2}$ and consider the Hamiltonian
	\[
	\mathcal{H}(q,p) = g^{-1}_q(p,p)+V(s) ,
	\] 
	where $g = \sigma^{-2}\boldsymbol{g}$. Then $\mathcal{H}$ is locally second-order maximally superintegrable if and only if there exists a $W_{ABCD}$ with Weyl symmetries which is the trace free part of some parallel $R{}_{ABCD}$ with curvature symmetries and such that $I^A R{}_{ABCD} =0$, which does not satisfy any of the equivalent properties in Lemma \ref{Lem:Indfromg} and
\begin{equation}\label{Eqn:3DI1}
  	3\tilde I^AX^B Z^C_c Z^D_d X^E W_{AB[CD}\tilde I_{E]}
=\tilde{\sigma} X^BZ^C_cZ^D_dX^EW_{AB[CD}\hat D_{E]}\tilde I^A,
\end{equation}
$\tilde{\sigma} = f(s)\sigma$ and $\tilde{I}_A = D_A(\tilde{\sigma})$.
\end{theorem}
\begin{proof}
	Suppose that $R{}_{ABCD} $ is a parallel tractor with curvature symmetries such that $I{}^{A} R{}_{ABCD} =0$ and suppose that $W = \operatorname{trace-free}(R)$ satisfies \eqref{Eqn:3DI1}. Then, by Theorem \ref{Thm:KS2}, $\tilde{\sigma}$ is a local Killing scale for $k{}_{bd}  = X{}^{A} Z{}^{B}{}_{b} X^C Z{}^{D}{}_{d} W{}_{ABCD} $, so by Theorem \ref{Thm:QCns}, there is a constant of motion for $\mathcal{H}$ with quadratic part $r{}_{bd} = X{}^{A} Z{}^{B}{}_{b} X{}^{C} Z{}^{D}{}_{d} R{}_{ABCD} $. By Proposition \ref{Prp:Find}, there are $2n-3$ functionally independent linear constants of motion induced by conformal Killing vectors $k_b = X^A Z{}^{B}{}_{b} K{}_{AB} $, where $K{}_{AB} $ is parallel and orthogonal to $I$ and $J$. Then by Lemma \ref{Lem:gind}, the function induced by $g^{-1}$ gives $2n-2$ functionally independent functions, and by Lemma \ref{Lem:Indfromg} \ref{Lem:Fin:2:4}, the function induced by $r{}_{bd} $ is functionally independent from $g^{-1}$. Finally, by Proposition \ref{Prp:Ind}, the $2n-3$ linear constants of motion, the Hamiltonian, with quadratic part $g^{-1}$ and the final constant of motion with quadratic part $r{}_{bd} $ are $2n-1$ functionally independent constants of motion and make $\mathcal{H}$ second-order maximally superintegrable.

	Quadratic constants of motion for $\mathcal{H}$, modulo constant multiples of $\mathcal{H}$ are in one-to-one correspondence with conformal Killing tensors for which $g$ and $\tilde{g}$ are both Killing scales. Suppose that $\mathcal{H}$ has $2n-2$ such constants, $f_{1},\dots ,f_{2n-2}$, with the addition of $\mathcal{H}$ making $\mathcal{H}$ second-order maximally superintegrable. The quadratic parts of these constants $r^1,\dots ,r^{2n-2}$ are Killing tensors for $g$. Then let $R{}^{1}{}_{ABCD} ,\dots ,R^{2n-2}_{ABCD}$ be the corresponding prolongations of the Killing tensors.

	Suppose for contradiction that all of the $R{}^{i} $ satisfy the conditions of Lemma \ref{Lem:Indfromg}. Then we can adjust each constant of motion by a constant multiple of $\mathcal{H}$ such that the new quadratic parts satisfy $ \hat{r}{}^{i}{}_{ab}\nabla^a \tau =0$. This does not change functional dependence. By Lemma \ref{Lem:Indfromg}, the $\hat{r}^i$ are all Killing for $\hat{g} = \tau^{-2}\boldsymbol{g}$, so by Proposition \ref{Prop:fKS} each is also Killing for $(V(s)g)$. Hence, by Proposition \ref{Prp:PBKT}, the function induced by each $\hat{r}^i$ is itself a constant of motion, and the functional part of each adjusted constant of motion is constant. 

	However, this means that functional independence of these constants reduces to that of the $\hat{r}^i$ and $g$, which cannot hold by Lemma \ref{Lem:Indfromg} \ref{Lem:Fin:2:4}, since every one of these functions is functionally dependent on $k_b = X^A Z{}^{B}{}_{b} K{}_{AB} $ for parallel $K\in \Gamma(\Lambda^2(I^{\perp}\cap J^{\perp}))$ and $g$, which can generate at most $2n-2$ functionally independent functions on $T^*M$, including $\mathcal{H}$, by Proposition \ref{Prp:Find}.\end{proof}
In the above proof we did not need to consider linear constants as they can be squared to obtain quadratic constants without changing functional dependence. If we had instead not included $\mathcal{H}$ in our list of constants, then $\mathcal{H}$ is functionally dependent on the $2n-1$ constants, since this is the maximal number. So we can replace one constant by $\mathcal{H}$ and continue the above proof.

In the following two sections, we will construct and classify potentials $V(s)$ in the above Theorem which are second-order superintegrable. By the end of those sections, we will no longer need to consider local constants of motion, since we will show that local second-order maximal superintegrability can only be achieved by specific choices of $V$, all of which we will construct global constants of motion for.

\section{Existence}\label{Sec:Ex}

In this section we apply the results obtained in the previous section to construct second-order maximally superintegrable Hamiltonians. In Section \ref{Sec:Uni} we show that these are unique among the families of Hamiltonians with potential a smooth function of an Einstein rescaling factor. In Section \ref{Sec:FM} we explicitly construct some examples of these Hamiltonians.

\begin{lemma}\label{Lem:Big}
	Let $(M,\boldsymbol{c})$ be conformally flat. Let $\sigma,\tau\in \Gamma(\mathcal{E}[1])$ be almost Einstein scales and let $I_A = D_A \sigma$ and $J_A = D_A \tau$ be the corresponding parallel tractors.
	Let $s = \tau\sigma^{-1}\in C^\infty (M\setminus\mathcal{Z}(\sigma))$, let $f\colon \mathbb{R}\to \mathbb{R}$ be smooth, let $\tilde{\sigma} = f(s)\sigma$ and let $\tilde{I}_A = D_A \tilde{\sigma}$.
	Suppose that $k_{bd}\in \Gamma(\mathcal{E}_{(bd)_{0}} [4])$ is a conformal Killing tensor for which $\sigma$ is a local Killing scale. Then $\tilde{\sigma}$ is a Killing scale if and only if
	\begin{equation}\label{Eqn:Big}
		{X}{}^{B}Z{}^{C}{}_{c} Z{}^{D}{}_{d}  {X}{}^{E}\left(  u(s)\left[  {I}{}^{A}{W}{}_{AB[CD}{J}{}_{E]}+{J}{}^{A}{W}{}_{AB[CD}{I}{}_{E]}\right] +v(s){J}{}^{A}{W}{}_{AB[CD}{J}{}_{E]}\right)=0,
	\end{equation}
	where $u(s) = (3(f(s)-sf'(s))f'(s)+sf(s)f''(s))$ and $v(s) = (3f'(s)^2-f(s)f''(s))$. Moreover, if $R_{ABCD}\in \Gamma(I^\perp)$ is such that $W$ is the trace-free part of $R$, then
	\begin{align*}
		{X}{}^{B}Z{}^{C}{}_{c} Z{}^{D}{}_{d}  {X}{}^{E}\left[  {I}{}^{A}{W}{}_{AB[CD}{J}{}_{E]}+{J}{}^{A}{W}{}_{AB[CD}{I}{}_{E]}\right] &= J^A{X}{}^{B}Z{}^{C}{}_{c} Z{}^{D}{}_{d}  {X}{}^{E}R{}_{AB[CD} I{}_{E]} \quad \text{and}\\
		J{}^{A}X{}^{B} Z{}^{C}{}_{c} Z{}^{D}{}_{d} X{}^{E} W{}_{AB[CD} J{}_{E]}  &= J{}^{A}X{}^{B} Z{}^{C}{}_{c} Z{}^{D}{}_{d} X{}^{E} R{}_{AB[CD} J{}_{E]}.
	\end{align*}
	So \eqref{Eqn:Big} is equivalent to
	\begin{equation}\label{Eqn:BigR}
	 	{X}{}^{B}Z{}^{C}{}_{c} Z{}^{D}{}_{d}  {X}{}^{E}\left(  u(s){J}{}^{A}{R}{}_{AB[CD}{I}{}_{E]}+v(s){J}{}^{A}{R}{}_{AB[CD}{J}{}_{E]}\right)=0.
	 \end{equation}
 \end{lemma}
\begin{proof}
	By Theorem \ref{Thm:KS2}, $\tilde{\sigma} $ is a Killing scale if and only if
	\begin{equation}\label{Eqn:3DI}
  	3\tilde I^AX^B Z{}^{C}{}_{c}  Z{}^{D}{}_{d}  X^E W_{AB[CD}\tilde I_{E]}
=\tilde{\sigma} X^B Z{}^{C}{}_{c} Z{}^{D}{}_{d} X^EW_{AB[CD}\hat D_{E]}\tilde I^A.
	\end{equation}
	First observe that adding a tractor involving $X$ or $h$ to either $\tilde{I}_A$ or $D_{A} \tilde{I}_{B} $ does not change the value of \eqref{Eqn:3DI}, as it is killed by the projectors. So first we calculate $\tilde{I}_A$ mod $X$, in the scale $\sigma$,
\[\tilde{I}_A = D_A\tilde\sigma \overset{\mod X_A}{=} \begin{pmatrix}
	f(s)\sigma \\ \nabla_a\tilde\sigma
\end{pmatrix} = \begin{pmatrix}
	f(s)\sigma\\ \sigma f'(s)\nabla_as 
\end{pmatrix} = (f(s)-sf'(s))\begin{pmatrix}
	\sigma\\0
\end{pmatrix}+f'(s)\begin{pmatrix}
	\sigma s\\ \sigma\nabla_a s
\end{pmatrix}\]
So $\tilde{I}_A = (f(s)-sf'(s))I_A+f'(s)J_A \mod X_A$.
Since $D_A(\phi X_B) = 0$ $\mod X,h$ for any density $\phi$,
\begin{align*}
	D_A\tilde{I}_B 
	&= Z{}_{A}{}^{a} \nabla_a\left((f(s)-sf'(s))I_A+f'(s)J_A\right) \mod X,h\\ 
	&= (-sf''(s))(\nabla_as) I_B+f''(s)(\nabla_as) J_B \\
	&= \sigma^{-1}f''(s)(J_A-sI_A)(J_B-sI_B)\mod X,
\end{align*} 
since $\sigma Z{}_{A}{}^{a} \nabla_as = J_A-sI_A$ mod $X$.
Hence, taking left hand side minus right hand side in \eqref{Eqn:3DI} becomes 
\begin{align}\label{Eqn:WBig}
	0={}X^B Z{}^{C}{}_{c} Z{}^{D}{}_{d}  X^E\Big(& (3(f(s)-sf'(s))^2-s^2f(s)f''(s))I^AW_{AB[CD}I_{E]}\notag\\
	& + (3(f(s)-sf'(s))f'(s)+sf(s)f''(s))I^AW_{AB[CD}J_{E]}\notag\\
	& + (3(f(s)-sf'(s))f'(s)+sf(s)f''(s))J^AW_{AB[CD}I_{E]}\notag\\
	& + (3f'(s)^2-f(s)f''(s))J^AW_{AB[CD}J_{E]}\Big).
\end{align}
By assumption, $\sigma$ is a Killing scale for $W$. So $I{}^{A} W{}_{AB[CD} I{}_{E]} =0$, by Theorem \ref{Thm:IWI}, completing the proof of \eqref{Eqn:Big}.

Now suppose that $W$ is the trace-free part of some $R\in I^\perp$, that is
\[W_{ABCD} = R_{ABCD}-(\mathsf{P}\owedge h)_{ABCD},\]
where $\mathsf{P}_{BD} = \tfrac{1}{n}(R{}_{AB}{}^{A}{}_{D} -rh_{BD})$, with $r = \frac{1}{2(n+1)}R{}_{AB}{}^{AB} $. Applying $J^A$ gives
\[
	J^AW_{ABCD} = J^AR_{ABCD}-2\mathsf{P}_{B[D}J_{C]}-2J^A\mathsf{P}_{A[C}h_{D]B},  \quad \text{and} \quad
\]
\begin{align*}
	  {I}{}^{A}{W}{}_{AB[CD}{J}{}_{E]}+{J}{}^{A}{W}{}_{AB[CD}{I}{}_{E]} 
	  &= -2\mathsf{P}_{B[D}I_{C}J_{E]}-2I^A\mathsf{P}_{A[C}J_{E}h_{D]B}+J^AR_{AB[CD}I_{E]}\\
	  &\mathrel{\phantom{=}} -2\mathsf{P}_{B[D}J_{C}I_{E]}-2J^A\mathsf{P}_{A[C}I_{E}h_{D]B}.
\end{align*} 
Since $X{}^{B} Z{}^{C}{}_{c} h{}_{BC} =0$, this gives 
\[
{X}{}^{B}Z{}^{C}{}_{c} Z{}^{D}{}_{d}  {X}{}^{E}\left[  {I}{}^{A}{W}{}_{AB[CD}{J}{}_{E]}+{J}{}^{A}{W}{}_{AB[CD}{I}{}_{E]}\right] = J^A{X}{}^{B}Z{}^{C}{}_{c} Z{}^{D}{}_{d} {X}{}^{E}R{}_{AB[CD} I{}_{E]}.
\] 
Similarly,
\[
	J{}^{A} W{}_{AB[CD} J{}_{E]}  = J^AR_{AB[CD}J_{E]}-2\mathsf{P}_{B[D}J_{C}J_{E]}-2J^A\mathsf{P}_{A[C}h_{D|B|}J_{E]}, \text{so}
\] 
\[
J{}^{A}X{}^{B} Z{}^{C}{}_{c} Z{}^{D}{}_{d} X{}^{E} W{}_{AB[CD} J{}_{E]}  = J{}^{A}X{}^{B} Z{}^{C}{}_{c} Z{}^{D}{}_{d} X{}^{E} R{}_{AB[CD} J{}_{E]},
\]
and \eqref{Eqn:BigR} follows.
\end{proof}

By Theorem \ref{Thm:QCns}, rank 2 Killing tensors which are Killing for both $\sigma$ and $\tilde{\sigma}$ correspond to constants of motion for $\mathcal{H} = g^{-1}+V(s)$, where $V(s) = f(s)^{-2}$. Then 
\begin{equation}\label{Eqn:hu}
	2v(s) = f(s)^4 V''(s) \quad \text{and} \quad -2u(s) = f(s)^4(sV''(s)+3V'(s)).  
\end{equation} 

We now need a Lemma to recover constants of motion from the parallel tractor $R{}_{ABCD} $.
\begin{lemma}\label{Lem:COM}
	Let $(M,\boldsymbol{c})$ be conformally flat. Let $R{}_{ABCD} $ be a parallel tractor with curvature symmetries. Let $I_A = D_A\sigma $ and $J_A = D_A\tau $ be parallel scale tractors and let $g = \sigma^{-2}\boldsymbol{g}$. Suppose that $I^A R{}_{ABCD} =0$ and that 
	\[
		J^A R{}_{ABCD}  = 2S{}_{B[D} T{}_{C]} 
	\] 
	for some parallel tractors $S\in \mathcal{E}{}_{(BD)} $ and $T\in \Gamma(\mathcal{E}_C)$. Let $W\in C^\infty(M) $ and let $V$ be a smooth real valued function defined outside of a finite singularity set on $\mathbb{R}$. Let $r{}_{bd}  = X^A Z{}^{B}{}_{b} X^C Z{}^{D}{}_{d} R{}_{ABCD} $. Then the function $r+W$ is a constant of motion for the Hamiltonian $\mathcal{H} = g^{-1}+V$ if and only if
	\[
	dW = V'(s)(\phi d\psi - \tfrac{1}{2}\psi d\phi),
	\] 
	where $\phi = \sigma^{-2}X^BX^D S{}_{BD} $ and $\psi = \sigma^{-1} X^C T_C$.
\end{lemma}
\begin{proof}
	By Theorem \ref{Thm:KT}, $r{}_{bd}  = X^A Z{}^{B}{}_{b}  X^C Z{}^{D}{}_{d} R{}_{ABCD} $ is a Killing tensor for $g$. So $r+W$ is a constant of motion for $\mathcal{H}$ if and only if
	\begin{equation}\label{Eqn:COM}
		0 = \left\{V,r\right\}+\left\{g^{-1},W\right\} = -r{}^{ab} (dV)_b + g^{ab} (dW)_b = -r{}^{ab} (dV)_b + \sigma^2 (dW)^a
	\end{equation} 
	Calculating in the scale $\sigma$,
	\[
		(d\psi)_c = \sigma^{-1}Z{}^{C}{}_{c} T_C\quad \text{and} \quad(d\phi)_d = 2\sigma^{-2}X{}^{B} Z{}^{D}{}_{d} S{}_{BD}.
	\] 
	Since $(dV)_a = V'(s)(ds)_a = \sigma^{-1}V'(s)\nabla_a \tau$, we have
	\[
		\begin{aligned}
			r{}_{bd} (dV)^d 
		&= \sigma^{-1}V'(s)X^A Z{}^{B}{}_{b} X^C J^D R{}_{ABCD} \\
		&= 2V'(s)\sigma^{-1} X{}^{B} Z{}^{C}{}_{b} X{}^{D} S{}_{B[D} T{}_{C]}  \\
		&=  V'(s)\sigma^2 (\phi(d\psi)_b-\tfrac{1}{2} \psi(d\phi)_b),\\
		\end{aligned}
	\] 
	so the result follows from \eqref{Eqn:COM}.
\end{proof}

We now restrict to the case where the subspace spanned by $I$ and $J$ is non-degenerate, or equivalently $\kappa^2-\iota\jmath \neq 0$. Since we are considering potentials of the form $V(s) = V(\tau\sigma^{-1})$, we can rescale $I$ and $J$ by constants without changing the Hamiltonians considered. Moreover, we can also replace $J$ by $J+cI$ for any $c\in \mathbb{R}$, as this replaces $s$ by $s+c$. Through these allowed transformations, we can reduce to two cases. 

First suppose that $\iota=I_AI^A = 0$. Then, by non-degeneracy, $\kappa\neq 0$. We replace $J$ by $J-\frac{\jmath}{2\kappa}I$, so that $J$ is null, we then multiply $J$ by a constant so that $\kappa = I_AJ^A=1$ and $\iota=\jmath =0$. 

Then if $\iota = I_AI^A\neq 0$, we replace $J$ by $J-\frac{\kappa}{\iota}I$, making $J$ and $I$ orthogonal. Then we rescale to obtain $\kappa=0$, $\iota,\jmath \in \{-1,1\}$.
\begin{theorem}\label{Thm:Ex1}
	Let $(M,\boldsymbol{c})$ have trivial conformal holonomy. Let $\sigma,\tau\in \Gamma(\mathcal{E}[1])$ be almost Einstein scales and let $I_A = D_A \sigma$ and $J_A = D_A \tau$ be the corresponding parallel tractors. Let $s = \tau\sigma^{-1}\in C^\infty (M\setminus\mathcal{Z}(\sigma))$. Let $\iota = I_AI^A$ and $\jmath  = J_AJ^A$. Suppose that the subspace spanned by $I$ and $J$ is non-degenerate, then without loss of generality either
	\begin{tabel}
		\item \label{Thm:Ex1:1:IJ}$I$ and $J$ are null and $I_AJ^A=1$, or
		\item \label{Thm:Ex1:2:IJ}$I$ and $J$ are orthogonal and $\iota,\jmath \in \{-1,+1\} $.
	\end{tabel}
	Let $f\colon \mathbb{R}\to \mathbb{R} $ be smooth such that $f\circ s$ is non-constant on $M\setminus\mathcal{Z}(\sigma)$ and let $V(s) = f(s)^{-2}$. If $f$ solves the ODE corresponding to the cases above,
	\begin{tabel}
		\item\label{Thm:Ex1:1} $v(s) = 3f'(s)^2-f(s)f''(s) = 0$, equivalently $V''(s) = 0$ with solutions
			\[
				V(s) = \lambda s+C.
			\]
		\item\label{Thm:Ex1:2} $u(s) = 3(f(s)-sf'(s))f'(s)+sf(s)f''(s)=0$, equivalently $sV''(s)+3V'(s)=0$ with solutions
			\[
				V(s) = \frac{\lambda}{s^2}+C.
			\]
	\end{tabel}
	Then the Hamiltonian
	\begin{equation}\label{Eqn:Ex1:Ham}
		\mathcal{H}(q,p) = g^{-1}_q(p,p)+V(s),
	\end{equation}
	where $g = \sigma^{-2}\boldsymbol{g}$ is second-order maximally superintegrable. 

	The space of Killing tensors $r{}_{bd} $ representing quadratic parts of constants of motion for $\mathcal{H}$ is characterised by the prolongation $R{}_{ABCD} \in I^{\perp}$ satisfying
	\begin{equation}\label{Eqn:Ex1:R}
		 {R}{}_{AB[CD}{\pi}{}_{E]}{}^{A}=0,
	\end{equation} 
	where $\pi$ is the natural parallel projection onto $I^\perp\cap J^{\perp}$. The space of curvature tractors satisfying this condition in each case, is characterised by
	\begin{tabel}
		\item\label{Thm:Ex1:1:R} 
			$R{}_{ABCD} = R^0{}_{ABCD} + (I^2\owedge S)_{ABCD}$, where $R^0 \in \begin{picture}(10, 10)(0,2)
	\put(0, 10){\line(1,0){10}}
	\put(0, 5){\line(1,0){10}}
	\put(0, 0){\line(1,0){10}}
	\put(0, 10){\line(0,-1){10}}
	\put(5, 10){\line(0,-1){10}}
	\put(10, 10){\line(0,-1){10}}
\end{picture}(I^{\perp}\cap J^{\perp})$ and $S \in S^2(I^{\perp}\cap J^{\perp})$. The corresponding constant of motion is then
\begin{equation}\label{Eqn:Ex1:1:CM}
	r_q(p,p)-\frac{\lambda}{2} \phi(q),
\end{equation}
where $r{}_{bd}  = X^A Z{}^{B}{}_{b} X^C Z{}^{D}{}_{d} R{}_{ABCD} $ and $\phi = \sigma^{-2}X^AX^B S{}_{AB} $.
		\item\label{Thm:Ex1:2:R} 
			$R{}_{ABCD} = R^0{}_{ABCD} + (J^2\owedge S)_{ABCD}$, where $R^0 \in \begin{picture}(10, 10)(0,2)
	\put(0, 10){\line(1,0){10}}
	\put(0, 5){\line(1,0){10}}
	\put(0, 0){\line(1,0){10}}
	\put(0, 10){\line(0,-1){10}}
	\put(5, 10){\line(0,-1){10}}
	\put(10, 10){\line(0,-1){10}}
\end{picture}(I^{\perp}\cap J^{\perp})$ and $S \in S^2(I^{\perp}\cap J^{\perp})$. The corresponding constant of motion is then
\begin{equation}\label{Eqn:Ex1:2:CM}
	r_q(p,p)+\lambda \jmath  \frac{\phi(s)}{s{}^{2} },
\end{equation}
where $r{}_{bd}  = X^A Z{}^{B}{}_{b} X^C Z{}^{D}{}_{d} R{}_{ABCD} $ and $\phi = \sigma^{-2}X^AX^B S{}_{AB} $.
	\end{tabel}

	The space of conformal Killing tensors given by $r_{(bd)_0}$, where $r{}_{bd}$ represents the quadratic part of a constant of motion for $\mathcal{H}$ is characterised by the prolongation $W{}_{ABCD} $ satisfying
	\begin{equation}\label{Eqn:Ex1:W}
		{W}{}_{AB[CD}{\pi}{}_{E]}{}^{A} = 0 .
	\end{equation}
	\end{theorem}
\begin{proof}
	Let $\tilde{\sigma} = f(s)\sigma$ and note that if $u=v=0$, then $V'(s)=0$, so $f(s)$ is constant, contrary to our assumption. We first consider \ref{Thm:Ex1:1:IJ}. In this case 
	\[
		{\pi}{}_{AB} = {\delta}{}_{AB}- {I}{}_{A} {J}{}_{B} -{J}{}_{A} {I}{}_{B}
	\] 
	So, under the assumption $v(s)=0$, the conditions \eqref{Eqn:Big} and \eqref{Eqn:BigR} become 
	\[
		{}-u(s){X}{}^{B}Z{}^{C}{}_{c} Z{}^{D}{}_{d}  {X}{}^{E}{W}{}_{AB[CD}{\pi}{}_{E]}{}^{A}=0 \quad \text{and} \quad -u(s)X^B Z{}^{C}{}_{c} Z{}^{D}{}_{d} X{}^{E} R{}_{AB[CD} \pi{}_{E]}{}^{A} =0
	\] 

	The case \ref{Thm:Ex1:2:IJ} is similar. Let $\iota = I_AI^A$ and $\jmath =J_AJ^A$. Then 
	\[
		{\pi}{}_{AB} = {\delta}{}_{AB}- \iota{I}{}_{A} {I}{}_{B} -\jmath {J}{}_{A} {J}{}_{B}.
	\] 
	So, under the assumption $u(s)=0$, the conditions \eqref{Eqn:Big} and \eqref{Eqn:BigR} become
	\[
		{}-v(s){X}{}^{B} Z{}^{C}{}_{c} Z{}^{D}{}_{d} {X}{}^{E}{W}{}_{AB[CD}{\pi}{}_{E]}{}^{A}=0 \quad \text{and} \quad -v(s) X{}^{B} Z{}^{C}{}_{c} Z{}^{D}{}_{d} X^E R{}_{AB[CD} \pi{}_{E]}{}^{A}=0 
	\]

	In both cases, any compatible $R{}_{ABCD} $ and $W{}_{ABCD} $ must satisfy
	\begin{equation}\label{Eqn:Ex1:WR2}
	{X}{}^{B}Z{}^{C}{}_{c} Z{}^{D}{}_{d}  {X}{}^{E}{W}{}_{AB[CD}{\pi}{}_{E]}{}^{A}=0= {X}{}^{B}Z{}^{C}{}_{c} Z{}^{D}{}_{d}  {X}{}^{E}{R}{}_{AB[CD}{\pi}{}_{E]}{}^{A},
	\end{equation}
	which we will now strengthen to \eqref{Eqn:Ex1:R} and \eqref{Eqn:Ex1:W}. The tractors ${W}{}_{AB[CD}{\pi}{}_{E]}{}^{A}$ and ${R}{}_{AB[CD}{\pi}{}_{E]}{}^{A}$ have symmetries given by the Young diagram \begin{picture}(8, 10)(0,3)
	\put(0, 12){\line(1,0){8}}
	\put(0, 8){\line(1,0){8}}
	\put(0, 4){\line(1,0){4}}
	\put(0, 0){\line(1,0){4}}
	\put(0, 12){\line(0,-1){12}}
	\put(4, 12){\line(0,-1){12}}
	\put(8, 12){\line(0,-1){4}}
	\end{picture}. Moreover, ${W}{}_{AB[CD}{\pi}{}_{E]}{}^{A}$ is trace-free. The projectors in \eqref{Eqn:Ex1:WR2} take the invariant top slot of the trace-free part of the parallel tractors, so by Lemma \ref{Lem:BGG}, ${W}{}_{AB[CD}{\pi}{}_{E]}{}^{A}=0$. If $\mathsf{P}$ is such that $W = R-(h\owedge \mathsf{P})$, then
	\[
		0 = {W}{}_{AB[CD}{\pi}{}_{E]}{}^{A} = R{}_{AB[CD} \pi{}_{E]}{}^{A} -2\mathsf{P}{}_{A[C} \pi{}^{A}{}_{E} h{}_{D]B}.
	\] 
	The invariant top slot of $R{}_{AB[CD} \pi{}_{E]}{}^{A}$ is then 
	\[
		X^E Z{}^{C}{}_{c} \mathsf{P}{}_{A[C} \pi{}^{A}{}_{E]} = 3 X^B Y^DX^E Z{}^{C}{}_{c} \mathsf{P}{}_{A[C} \pi{}^{A}{}_{E} h{}_{D]B} = \tfrac{3}{2\sigma}X^B I^DX^E Z{}^{C}{}_{c}R{}_{AB[CD} \pi{}_{E]}{}^{A}=0
	\] 
	Hence, $\mathsf{P}{}_{A[C} \pi{}^{A}{}_{E]} = 0$ by Lemma \ref{Lem:BGG}, so $R{}_{AB[CD} \pi{}_{E]}{}^{A} = 0$.
	
	For case \ref{Thm:Ex1:1:R}, the equation \eqref{Eqn:Ex1:R} becomes
	\[
		R{}_{AB[CD} \pi{}_{E]}{}^{A}  = -J^A R{}_{AB[CD} I{}_{E]} = 0,
	\] 
	which can be viewed as characterising the kernel of the composition of
	\[
		{}\lrcorner J\colon \begin{picture}(10, 10)(0,2)
			\put(0, 10){\line(1,0){10}}
			\put(0, 5){\line(1,0){10}}
			\put(0, 0){\line(1,0){10}}
			\put(0, 10){\line(0,-1){10}}
			\put(5, 10){\line(0,-1){10}}
			\put(10, 10){\line(0,-1){10}}
			\end{picture}(I^{\perp})\to \begin{picture}(10, 10)(0,2)
			\put(0, 10){\line(1,0){10}}
			\put(0, 5){\line(1,0){10}}
			\put(0, 0){\line(1,0){5}}
			\put(0, 10){\line(0,-1){10}}
			\put(5, 10){\line(0,-1){10}}
			\put(10, 10){\line(0,-1){5}}
		\end{picture}(I^{\perp})
		\quad \textrm{with} \quad \wedge I\colon \begin{picture}(10, 10)(0,2)
	\put(0, 10){\line(1,0){10}}
	\put(0, 5){\line(1,0){10}}
	\put(0, 0){\line(1,0){5}}
	\put(0, 10){\line(0,-1){10}}
	\put(5, 10){\line(0,-1){10}}
	\put(10, 10){\line(0,-1){5}}
\end{picture}(I^{\perp})\to \begin{picture}(10, 15)(0,4.5)
	\put(0, 15){\line(1,0){10}}
	\put(0, 10){\line(1,0){10}}
	\put(0, 5){\line(1,0){5}}
	\put(0, 0){\line(1,0){5}}
	\put(0, 15){\line(0,-1){15}}
	\put(5, 15){\line(0,-1){15}}
	\put(10, 15){\line(0,-1){5}}
\end{picture}(I^{\perp}).
	\]
	The kernel of ${}\wedge I$ is equal to $\{U{}_{B[D} I{}_{C]} :U\in S^2(I^{\perp})\} $. The preimage of this under the map $\lrcorner J$ consists of $R$ of the form
	\[
		R{}_{ABCD} = T{}_{ABCD} + (I^2\owedge U)_{ABCD},
	\] 
	such that $J^B U{}_{B[D} I{}_{C]} =0$. This implies $J^B U{}_{BD}  = \lambda I_D$, which gives $U  \in S^2(I^{\perp}\cap J^{\perp})\oplus \operatorname{span}(I^2)$. Since $(I^2\owedge I^2)=0$, it suffices to take $U\in S^2(I^{\perp}\cap J^{\perp})$.

	Now we apply Lemma \ref{Lem:COM} to find the constant of motion. Since $J^A R{}_{ABCD}  = 2U{}_{B[D} I{}_{C]} $, we set $\phi = \sigma^{-2}X^BX^D U{}_{BD} $ and $\psi = \sigma^{-1}X^CI_C = 1$. Since $V'(s) = \lambda$, $r+W$ is a constant of motion if and only if
	\[
		dW = \lambda(\phi d\psi-\tfrac{1}{2}\psi d\phi) = -\tfrac{\lambda}{2}d\phi.
	\] 
	So $r-\frac{\lambda}{2}\phi$ is a constant of motion and $\mathcal{H}$ is second-order maximally superintegrable.

	The case \ref{Thm:Ex1:2:R} is similar. First \eqref{Eqn:Ex1:R} becomes $J^A R{}_{AB[CD} J{}_{E]} =0$, so we consider the kernel of $\wedge J\circ \lrcorner J$. This kernel is the preimage $(\lrcorner J)^{-1} \{U{}_{B[D} J{}_{C]} :U\in S^2(I^{\perp})\}   $, which consists of $R$ of the form 
	\[
		R{}_{ABCD} = T{}_{ABCD} + (J^2\owedge U)_{ABCD},
	\] 
	such that $J^B U{}_{B[D} J{}_{C]} =0$. Eliminating $U$ of the form $J{}^{2}$ shows that it suffices to take $U\in S^2(I^{\perp}\cap J^{\perp})$. 

	To find the corresponding constant of motion, we observe that $J^A R{}_{ABCD} = \jmath 2U{}_{B[D} J{}_{C]} $. So we set $\phi = \sigma^{-2}X^BX^D U{}_{BD} $ and $\psi = \jmath \sigma^{-1} X^CJ_C = \jmath s$. Since $V'(s) = -2\lambda s^{-3}$, $r+W$ is a constant of motion if and only if
	\[
	dW = -\frac{2\lambda\jmath }{s^3}\left( \phi ds-\frac{1}{2}sd\phi \right) =   \lambda\jmath  d\left( \frac{\phi}{s^2} \right).
	\] 
	So $r+\lambda \jmath  \frac{\phi}{s{}^{2} }$ is a constant of motion and $\mathcal{H}$ is second-order maximally superintegrable.

\end{proof}
The above Theorem gives the harmonic oscillator in the case where $I$ and $J$ are Euclidean scales, whereas the following Theorem gives the Kepler-Coulomb System in that case.
\begin{theorem}\label{Thm:Ex2}
	Let $(M,\boldsymbol{c})$ have trivial conformal holonomy. Let $\sigma,\tau\in \Gamma(\mathcal{E}[1])$ be almost Einstein scales and let $I_A = D_A \sigma$ and $J_A = D_A \tau$ be the corresponding parallel tractors. Let $s = \tau\sigma^{-1}\in C^\infty (M\setminus\mathcal{Z}(\sigma))$. Let $\iota = I_AI^A$ and $\jmath  = J_AJ^A$. Suppose that the subspace spanned by $I$ and $J$ is non-degenerate, then without loss of generality either
	\begin{tabel}
		\item\label{Thm:Ex2:1:IJ} $I$ and $J$ are null and $I_AJ^A=1$, or
		\item\label{Thm:Ex2:2:IJ} $I$ and $J$ are orthogonal and $\iota,\jmath \in \{-1,+1\} $.
	\end{tabel}
	Let $f\colon \mathbb{R}\to \mathbb{R} $ be smooth such that $f\circ s$ is non-constant on $M\setminus\mathcal{Z}(\sigma)$ and let $V(s) = f(s)^{-2}$. If $f$ solves the ODE corresponding to the cases above;
	\begin{tabel}
	\item\label{Thm:Ex2:1} $u(s)-sv(s) = 3f(s)f'(s)-6sf'(s)^{2}+2sf(s)f''(s) = 0$, equivalently \\$2sV''(s)+3V'(s)=0$ with solutions
			\[
				V(s) = \frac{\lambda}{\sqrt{|s|}}+C,
			\] 
		\item\label{Thm:Ex2:2} $s\iota\jmath u(s)-v(s) = 3(\iota\jmath sf(s)-\iota\jmath s^2f'(s)-f'(s))f'(s)+(\iota\jmath s^{2}+1)f(s)f''(s)=0$, equivalently $(\iota\jmath s^2+1)V''(s)+3\iota\jmath sV'(s)=0$ with solutions
			\[
				V(s)=\frac{\lambda s}{\sqrt{|\iota\jmath s^2+1|}}+C,
			\] 
	\end{tabel}
	Then the Hamiltonian
	\begin{equation}
		\mathcal{H}(q,p) = g^{-1}_q(p,p)+V(s),
	\end{equation}
	where $g = \sigma^{-2}\boldsymbol{g}$ is second-order maximally superintegrable. 
	The space of Killing tensors $r{}_{bd} $ representing quadratic parts of constants of motion for $\mathcal{H}$ is characterised by the prolongation $R{}_{ABCD} \in I^{\perp}$ satisfying
	\begin{equation}\label{Eqn:Ex2:R}
		J{}^{A} R{}_{ABCD} = 2\pi{}_{B[D}U{}_{C]},
	\end{equation} 
	for some parallel $U\in I^{\perp}$, where $\pi$ is the natural parallel projection onto $I^\perp\cap J^{\perp}$. Then the constants of motion for \ref{Thm:Ex2:1:IJ} and \ref{Thm:Ex2:2:IJ} are
	\begin{equation}\label{Eqn:Ex2:1:CM}
		r_q(p,p)+\lambda \frac{\psi(q)}{\sqrt{|s|} },\quad \text{and} \quad
	\end{equation}
	\begin{equation}\label{Eqn:Ex2:2:CM}
			r_q(p,p)-\lambda \iota  \frac{\psi(q)}{\sqrt{|\iota\jmath s^2+1|}},
	\end{equation}
respectively, where $r{}_{bd}  = X^A Z{}^{B}{}_{b} X^C Z{}^{D}{}_{d} R{}_{ABCD} $ and $\psi = \sigma^{-1}X^A U{}_{A} $.
\end{theorem}\begin{proof}
	Let $\tilde{\sigma} = f(s)\sigma$ and note that if $u=v=0$, then $V'(s)=0$, so $f(s)$ is constant, contrary to our assumption. We first consider \ref{Thm:Ex2:1:IJ}. In this case 
	\[
		{\pi}{}_{AB} = {\delta}{}_{AB}-{I}{}_{A}{J}{}_{B} -{J}{}_{A} {I}{}_{B}.
	\]
	So, under the assumption $u(s)=sv(s)$, the condition \eqref{Eqn:BigR} becomes
	\begin{align*}
		0&= X^B Z{}^{C}{}_{c} Z{}^{D}{}_{d} X^E\left( sv(s) J{}^{A} R{}_{AB[CD} I{}_{E]} + v(s)J{}^{A} R{}_{AB[CD} J{}_{E]} \right) \\
		 &= v(s)X^B Z{}^{C}{}_{c} Z{}^{D}{}_{d} X{}^{E} J{}^{A} R{}_{AB[CD} (J+sI)_{E]} \\
		 &= -v(s)\sigma^{-1}J^AX^B Z{}^{C}{}_{c} Z{}^{D}{}_{d} X{}^{E} X{}^{F} R{}_{AB[CD} \pi{}_{E]F}.
	\end{align*}
		Taking $R = ((I\vee U)\owedge \pi)$ for some $U\in \Gamma(I^{\perp}\cap J^{\perp})$ gives $R\in I^{\perp}\setminus J^{\perp}$ and 
	\[
		J{}^{A}X^BX^F R{}_{AB[CD} \pi{}_{E]F} = X^B\pi{}_{B[D} U{}_{C} \pi{}_{E]F}X^F=0,
	\]
The case \ref{Thm:Ex2:2:IJ} is similar. In this case, since $\iota^2 = \jmath^2=1$, 
	\[
		{\pi}{}_{AB} = {\delta}{}_{AB}-\iota{I}{}_{A}{I}{}_{B} -\jmath {J}{}_{A} {J}{}_{B}.\]
	So, under the assumption $v(s) = s \iota\jmath  u(s)$, the condition \eqref{Eqn:BigR} becomes
	\begin{align*}
		0&= X^B Z{}^{C}{}_{c} Z{}^{D}{}_{d} X^E\left( u(s) J{}^{A} R{}_{AB[CD} I{}_{E]} + \iota\jmath su(s)J{}^{A} R{}_{AB[CD} J{}_{E]} \right) \\
		 &= \iota u(s)X^B Z{}^{C}{}_{c} Z{}^{D}{}_{d} X{}^{E} J{}^{A} R{}_{AB[CD} (s\jmath J+\iota I)_{E]} \\
		 &= -\iota u(s)\sigma^{-1}J^AX^B Z{}^{C}{}_{c} Z{}^{D}{}_{d} X{}^{E} X{}^{F} R{}_{AB[CD} \pi{}_{E]F}.
	\end{align*}
	Taking $R = ((J\vee U)\owedge \pi)$ for any $U\in \Gamma(I^{\perp}\cap J^{\perp})$ gives $R\in I^{\perp}\setminus J^{\perp}$ and 
	\[
		J{}^{A} X^BX^F R{}_{AB[CD} \pi{}_{E]F}= X^B\pi{}_{B[D} U{}_{C} \pi{}_{E]F}X^F= 0,
	\]

	In both cases, we showed that $R\in I^{\perp}$ satisfies
	\begin{equation}\label{Eqn:Ex2:JR}
		J^AX^B Z{}^{C}{}_{c} Z{}^{D}{}_{d} X{}^{E} X{}^{F} R{}_{AB[CD} \pi{}_{E]F}
	\end{equation}

	We have
	\[
		J^AR{}_{AB[CD} \pi{}_{EF]}=0\quad \text{and} \quad J^AR{}_{A[BCD} \pi{}_{E]F}=0,
	\] 
	which forces $J^AR{}_{AB[CD} \pi{}_{E]F}$ to have symmetries given by 
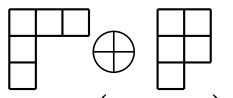
\begin{picture}(15, 15)(0,4)
	\put(0, 15){\line(1,0){15}}
	\put(0, 10){\line(1,0){15}}
	\put(0, 5){\line(1,0){5}}
	\put(0, 0){\line(1,0){5}}
	\put(0, 15){\line(0,-1){15}}
	\put(5, 15){\line(0,-1){15}}
	\put(10, 15){\line(0,-1){5}}
	\put(15, 15){\line(0,-1){5}}
\end{picture}$\oplus$ \begin{picture}(10, 15)(0,4)
	\put(0, 15){\line(1,0){10}}
	\put(0, 10){\line(1,0){10}}
	\put(0, 5){\line(1,0){10}}
	\put(0, 0){\line(1,0){5}}
	\put(0, 15){\line(0,-1){15}}
	\put(5, 15){\line(0,-1){15}}
	\put(10, 15){\line(0,-1){10}}
\end{picture}\;. Projection onto the first Young diagram is equivalent to symmetrising the pair $(BF)$ and projection onto the second Young diagram is equivalent to skew-symmetrising the pair $[BF]$, since
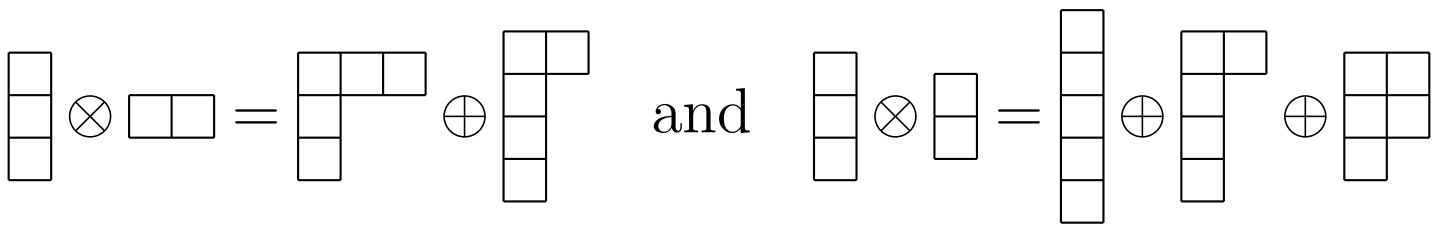
\[
\begin{picture}(8, 16)(0,9)
	\put(0, 24){\line(1,0){8}}
	\put(0, 16){\line(1,0){8}}
	\put(0, 8){\line(1,0){8}}
	\put(0, 0){\line(1,0){8}}
	\put(0, 24){\line(0,-1){24}}
	\put(8, 24){\line(0,-1){24}}
\end{picture}\otimes \begin{picture}(16, 8)(0,1)
	\put(0, 8){\line(1,0){16}}
	\put(0, 0){\line(1,0){16}}
	\put(0, 8){\line(0,-1){8}}
	\put(8, 8){\line(0,-1){8}}
	\put(16, 8){\line(0,-1){8}}
\end{picture} = \begin{picture}(24, 16)(0,9)
	\put(0, 24){\line(1,0){24}}
	\put(0, 16){\line(1,0){24}}
	\put(0, 8){\line(1,0){8}}
	\put(0, 0){\line(1,0){8}}
	\put(0, 24){\line(0,-1){24}}
	\put(8, 24){\line(0,-1){24}}
	\put(16, 24){\line(0,-1){8}}
	\put(24, 24){\line(0,-1){8}}
\end{picture}\oplus \begin{picture}(16, 16)(0,13)
	\put(0, 32){\line(1,0){16}}
	\put(0, 24){\line(1,0){16}}
	\put(0, 16){\line(1,0){8}}
	\put(0, 8){\line(1,0){8}}
	\put(0, 0){\line(1,0){8}}
	\put(0, 32){\line(0,-1){32}}
	\put(8, 32){\line(0,-1){32}}
	\put(16, 32){\line(0,-1){8}}
\end{picture}\quad \text{and} \quad \begin{picture}(8, 24)(0,9)
	\put(0, 24){\line(1,0){8}}
	\put(0, 16){\line(1,0){8}}
	\put(0, 8){\line(1,0){8}}
	\put(0, 0){\line(1,0){8}}
	\put(0, 24){\line(0,-1){24}}
	\put(8, 24){\line(0,-1){24}}
\end{picture}\otimes \begin{picture}(8, 16)(0,5)
	\put(0, 16){\line(1,0){8}}
	\put(0, 8){\line(1,0){8}}
	\put(0, 0){\line(1,0){8}}
	\put(0, 16){\line(0,-1){16}}
	\put(8, 16){\line(0,-1){16}}
\end{picture} = \begin{picture}(8, 16)(0,17)
	\put(0, 40){\line(1,0){8}}
	\put(0, 32){\line(1,0){8}}
	\put(0, 24){\line(1,0){8}}
	\put(0, 16){\line(1,0){8}}
	\put(0, 8){\line(1,0){8}}
	\put(0, 0){\line(1,0){8}}
	\put(0, 40){\line(0,-1){40}}
	\put(8, 40){\line(0,-1){40}}
\end{picture}\oplus \begin{picture}(16, 16)(0,13)
	\put(0, 32){\line(1,0){16}}
	\put(0, 24){\line(1,0){16}}
	\put(0, 16){\line(1,0){8}}
	\put(0, 8){\line(1,0){8}}
	\put(0, 0){\line(1,0){8}}
	\put(0, 32){\line(0,-1){32}}
	\put(8, 32){\line(0,-1){32}}
	\put(16, 32){\line(0,-1){8}}
\end{picture}\oplus \begin{picture}(16, 16)(0,9)
	\put(0, 24){\line(1,0){16}}
	\put(0, 16){\line(1,0){16}}
	\put(0, 8){\line(1,0){16}}
	\put(0, 0){\line(1,0){8}}
	\put(0, 24){\line(0,-1){24}}
	\put(8, 24){\line(0,-1){24}}
	\put(16, 24){\line(0,-1){16}}
\end{picture}
\]
\\
The projectors in \eqref{Eqn:Ex2:JR} take the top slot of the completely trace-free component of
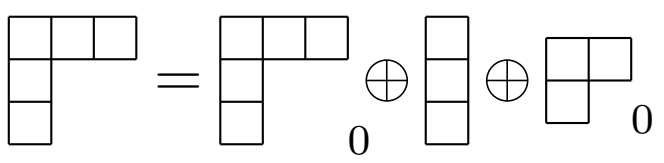
\[\begin{picture}(24, 16)(0,9)
	\put(0, 24){\line(1,0){24}}
	\put(0, 16){\line(1,0){24}}
	\put(0, 8){\line(1,0){8}}
	\put(0, 0){\line(1,0){8}}
	\put(0, 24){\line(0,-1){24}}
	\put(8, 24){\line(0,-1){24}}
	\put(16, 24){\line(0,-1){8}}
	\put(24, 24){\line(0,-1){8}}
\end{picture} = \begin{picture}(24, 16)(0,9)
	\put(0, 24){\line(1,0){24}}
	\put(0, 16){\line(1,0){24}}
	\put(0, 8){\line(1,0){8}}
	\put(0, 0){\line(1,0){8}}
	\put(0, 24){\line(0,-1){24}}
	\put(8, 24){\line(0,-1){24}}
	\put(16, 24){\line(0,-1){8}}
	\put(24, 24){\line(0,-1){8}}
	\put(24,0){${}_0$}
\end{picture}\oplus \begin{picture}(8, 16)(0,9)
	\put(0, 24){\line(1,0){8}}
	\put(0, 16){\line(1,0){8}}
	\put(0, 8){\line(1,0){8}}
	\put(0, 0){\line(1,0){8}}
	\put(0, 24){\line(0,-1){24}}
	\put(8, 24){\line(0,-1){24}}
\end{picture}\oplus \begin{picture}(16, 16)(0,5)
	\put(0, 16){\line(1,0){16}}
	\put(0, 8){\line(1,0){16}}
	\put(0, 0){\line(1,0){8}}
	\put(0, 16){\line(0,-1){16}}
	\put(8, 16){\line(0,-1){16}}
	\put(16, 16){\line(0,-1){8}}
	\put(16,0){${}_0$}
\end{picture}\]
which then must vanish. The first trace part is given by $J^AR{}_{AB[CD} \pi{}_{E]}{}^{B} =0$. The projector onto the top slot of the final component is given by
\[
	J^A X^B Z{}^{C}{}_{c} X^D Y^E X^F R{}_{AB[CD} \pi{}_{E]F} = \sigma^{-1}J^A X^B Z{}^{C}{}_{c} X^D I^E X^F R{}_{AB[CD} \pi{}_{E]F} = 0
\] 
Hence, $J^A R{}_{AB[CD} \pi{}_{E]F} $ is skew symmetric in the pair $BF$. Then, in case \ref{Thm:Ex2:1:IJ}	\begin{align*}
		0&=3\operatorname{tr}_{B,D}(J^A R{}_{AB[CD} \pi{}_{E]F}+J^A R{}_{AF[CD} \pi{}_{E]B})\\ 
		 &= 2J{}^{A} R{}_{AB[C}{}^{B} \pi{}_{E]F} +J^A R{}_{ABEC} \pi{}^{B}{}_{F}  +2J{}^{A} R{}_{AF[C}{}^B\pi{}_{E]B}  +nJ{}^{A} R{}_{AFEC}\\ 
		 &= 2J{}^{A} R{}_{AB[C}{}^{B}\pi{}_{E]F}  + J^A R{}_{AFEC} +2J^A R{}_{AFCE} -2J^A R{}_{AF[C}{}^{B} I{}_{E]} J_B+n J^A R{}_{AFEC} \\ 
		 &= (n-1)J^A R{}_{AFEC} + 2J{}^{A} R{}_{AB[C}{}^{B}\pi{}_{E]F}+2J^AJ^B R{}_{AFB[C} I{}_{E]},
	\end{align*}
	where in the case \ref{Thm:Ex2:2:IJ}, the final term is replaced by $2J^AJ^B R{}_{AFB[C} J{}_{E]}$. So we just need to show that $J^AJ^C R{}_{ABCD} $ is proportional to $\pi{}_{BD} $. Then, in either case, contracting $J{}^{C} $ into the above gives
	\[
		0 = (n-1)J^AJ^C R{}_{AFEC} + J^AJ^C R{}_{ABC}{}^{B}\pi{}_{EF} - J^AJ^B R{}_{AFBE}
	\] 
	So $nJ^AJ^C R{}_{ABCD}  = J^AJ^C R{}_{AEC}{}^{E} \pi{}_{BD}$ and in the case \ref{Thm:Ex2:1:IJ},
	\[
		J^A R{}_{ABCD}  = \pi{}_{B[D} U{}_{C]}, \text{ where } U_C = \frac{2}{n-1}\left( J^A R{}_{ABC}{}^{B} +\frac{1}{n}J^AJ^D R{}_{AB}{}^{B}{}_{D} I_C \right),
	\] 
	and in the case \ref{Thm:Ex2:2:IJ}, the final term above becomes $J^AJ^D R{}_{AB}{}^{B}{}_{D} J_C$.
	
	Now we apply Lemma \ref{Lem:COM} to find the constant of motion. Since $J^A R{}_{ABCD}  = 2\pi{}_{B[D} U{}_{C]} $, we set $\phi = \sigma^{-2}X^BX^D \pi{}_{BD} $ and $\psi = \sigma^{-1}X^CU_C$. 
	
	For the case \ref{Thm:Ex2:1:IJ}, we have $\phi = -2 s$ and $V'(s) = -\frac{\lambda s}{2|s|^{5 / 2} }$, so $r+W$ is a constant of motion if and only if
	\[
		dW = \frac{\lambda s}{|s|^{5 /2 }}(s d\psi-\frac{1}{2}\psi ds) = \lambda\left( \frac{d\psi}{\sqrt{|s|}} - \frac{\psi sds}{2|s|^{5/2}}  \right) = \lambda d\left( \frac{\psi}{\sqrt{|s|} } \right).
	\] 
	So $r_q(p,p)+\lambda \frac{\psi(q)}{\sqrt{|s|} }$ is a constant of motion and $\mathcal{H}$ is second-order maximally superintegrable.

	For the case \ref{Thm:Ex2:2:IJ}, we let $D = 1+\iota\jmath s^2$. Then $\phi = -\iota-\jmath s^2 = -\iota D$ and $V'(s) = \frac{\lambda (D)}{|D|^{5 / 2} }$, so $r+W$ is a constant of motion if and only if

	\[
		dW=V'(s)\bigl(\phi\,d\psi-\tfrac12\psi\,d\phi\bigr)=-\lambda\iota\,d\left(\frac{\psi}{\sqrt{|D|}}\right),
	\]
	So $r_q(p,p)-\lambda \iota\frac{\psi(q)}{\sqrt{|D|} }$ is a constant of motion and $\mathcal{H}$ is second-order maximally superintegrable.
\end{proof}

So we have constructed Hamiltonians that are second-order maximally superintegrable, we will calculate some examples and the associated constants of motion in Section \ref{Sec:FM}. The objective of the next section is to show that the above two Theorems classify all possible second-order maximally superintegrable potentials, under the given assumptions on $V$.

\section{Uniqueness}\label{Sec:Uni}

In this section we will show that, when $\operatorname{span}(I,J)$ is non-degenerate, the systems presented in Section \ref{Sec:Ex} are the only second-order maximally superintegrable systems on conformally flat Einstein manifolds which have a potential parameterised by an Einstein rescaling factor. 

Fix a regular value $c$ such that $\jmath-2c\kappa+c^2\iota\ne0$ for $s$ and consider the level set $\Sigma_c = \{q\in M:s(q)=c\} $. Observe that $s = \tau\sigma^{-1} = c$ is equivalent to $\tau-c\sigma=0$. So $\Sigma{}_{c} $ is the zero set of an almost Einstein scale. The study of tractor calculus for conformal hypersurfaces and submanifolds in general has been widely studied, \cite{curry_conformal_2025, gover_submanifold_2015}. Our situation consisting of the zero set of an almost Einstein scale is a particularly nice case, such hypersurfaces were considered in \cite{gover_almost_2010} and they are a special case of a holonomy reduction \cite{cap_holonomy_2014}. Associated to a conformal hypersurface is a normal tractor, $N_A = D_A \nu$, which in our case is the unit length multiple of $D_A(\tau-c\sigma) = J_A-cI_A$. Note that
\[
	(J-cI)_A(J-cI)^A = \jmath-2c\kappa+c^2\iota.
\] 
The intrinsic tractor bundle of $\Sigma_c$ is canonically identified with $N^{\perp}|_\Sigma = (J-cI)^{\perp}|_\Sigma$. Since $N$ is parallel, the tractor connection of $M$ preserves $N^{\perp}$. This connection agrees with the intrinsic tractor connection of $\Sigma{}_{c} $, which can be defined even for $dim(\Sigma_c)<3$ \cite{curry_conformal_2025}. The tractor metric of $\Sigma_c$ is the tractor metric of $M$ restricted to $N^{\perp}$. In particular, Lemma \ref{Lem:BGG} applies to tractors on $\Sigma{}_{c} $, which becomes the following lemma.

\begin{lemma}\label{Lem:SMBGG}
	Let $(M,\boldsymbol{c})$ be a connected conformal manifold, let $N_A = D_A \nu$ be a parallel tractor of non-zero length such that the zero set $\Sigma = \mathcal{Z}(\nu)$ is non-empty, then $\Sigma$ is an embedded hypersurface \emph{\cite[Theorem 3.5]{cap_holonomy_2014}}. Let $U\in \Gamma(\mathcal{V})$ for some irreducible trace-free tractor bundle $\mathcal{V}$ of $M$. Suppose that $U|_\Sigma$ is orthogonal to $N$, $\Pi(U|_\Sigma)=0$ and $U$ is parallel along $\Sigma$, that is, $\nabla_vU=0$ whenever $v\in T\Sigma$. Then $U|_\Sigma=0$.
\end{lemma}
\begin{proof}
	Since $U$ is orthogonal to $N$, $U|_\Sigma$ is a parallel section the subbundle of $\mathcal{V}|_\Sigma$, which is an irreducible tractor bundle on $\Sigma$. The projection $\Pi$ consists of a combination of $X$ and $Z$ contracted into $U$. As in \cite{curry_conformal_2025}, the projectors $X_\Sigma$ and $Z_\Sigma$ for tractors in $N^{\perp}$ are simply $X{}_{\Sigma}{}^{A}  = X^A$ and $Z{}_{\Sigma}{}^{A}{}_{a} = Z{}^{A}{}_{a} - \epsilon Z{}^{A}{}_{b} n^b n_a$, where $n_a\in \Gamma(\mathcal{E}_a[1])$ has length $n^an_a=\epsilon\in \{-1,+1\} $.

	From this we conclude that $\Pi_\Sigma(U)$, the top slot of $U$ taken with the projectors $X{}_{\Sigma}$ and $Z_\Sigma$ is $\Pi(U)$, projected onto $T\Sigma$. However, since $\Pi(U)=0$, we have $\Pi_\Sigma(U)=0$. In fact, one can use that $U$ is orthogonal to $N$ to show that $\Pi(U)$ already lies in a tensor power of $T\Sigma$, but we don't need that here. Hence, by Lemma \ref{Lem:BGG} applied to the conformal manifold $\Sigma$, $U|_\Sigma=0$.
\end{proof}
Similarly to our comment after Lemma \ref{Lem:BGG}, any specific case of the above lemma can be worked out manually. Let $N_A\in \Gamma(\mathcal{E}_A)$ and let $U_A\in \Gamma(\mathcal{E}_A)$, both satisfy the assumptions of the above lemma. Fix a scale $g = \sigma^{-2}\boldsymbol{g}\in \boldsymbol{c}$ to calculate in. Since $X^AU_A=0$ on $\Sigma$, we have 
\[
	0 = \nabla_v(X^AU_A) = v^aZ{}^{A}{}_{a} U_A
\] 
for all $v\in T\Sigma$. Then write $\nu = t\sigma$, so $\Sigma = \mathcal{Z}(t)$. Then on $\Sigma$
\[
Z{}^{A}{}_{a} N_A = \nabla_a\nu = \sigma\nabla_a t\propto n_a.
\] 
Since $X^AU_A = 0=X^AN_A$ on $\Sigma$, we also have
\[
0 = h{}^{AB} U_AN_B = Z{}^{A}{}_{a} Z{}^{Ba} U_AN_B = n^a Z{}^{A}{}_{a} U_A.
\] 
So $Z{}^{A}{}_{a}  U_A|_\Sigma=0$. Finally then
\[
0 = \nabla_v( Z{}^{A}{}_{a}  U_A) = -v_aY^A U_A
\] 
for all $v\in T\Sigma$, so $Y^AU_A=0$ on $\Sigma$, so $U|_\Sigma=0$. 

In this section we will need the above Lemma applied to a section of the tractor bundle with trace free \begin{picture}(10, 10)
	\put(0, 10){\line(1,0){10}}
	\put(0, 5){\line(1,0){10}}
	\put(0, 0){\line(1,0){5}}
	\put(0, 10){\line(0,-1){10}}
	\put(5, 10){\line(0,-1){10}}
	\put(10, 10){\line(0,-1){5}}
\end{picture} symmetries. The direct proof for this tractor bundle is similar to the above with more steps.
\begin{lemma}\label{Lem:U+}
	Let $(M,\boldsymbol{c})$ be conformally flat. Let $\sigma,\tau\in \Gamma(\mathcal{E}[1])$ be almost Einstein scales and let $I_A = D_A \sigma$ and $J_A = D_A \tau$ be the corresponding parallel tractors. Let $s = \tau\sigma^{-1}\in C^\infty (M\setminus\mathcal{Z}(\sigma))$, let $\tilde{\sigma} = f(s)\sigma$ for some non-constant $f:\mathbb{R}\to \mathbb{R}$ such that $\tilde{\sigma}$ is non-zero on an open dense set and let $K_A = J_A-sI_A$ on $M\setminus\mathcal{Z}(\sigma)$. Suppose that the space spanned by $I, J$ is non-degenerate, or equivalently $\kappa^2-\iota\jmath \neq 0$, where $\iota = I^AI_A$, $\jmath  = J^AJ_A$ and $\kappa = I^AJ_A$.

	Suppose that $r_{bd}\in \Gamma(\mathcal{E}_{(bd)} [4])$ is a Killing tensor for $\sigma$ and suppose that $\tilde{\sigma}$ is a Killing scale for $k{}_{bd}  = r{}_{(bd)_0} $. Then there is a function $\beta\colon \mathbb{R}\to \mathbb{R} $ such that whenever $K_AK^A= \iota s^2 -2\kappa s+\jmath \neq 0$,
	\begin{equation}\label{Eqn:U+}
		(\iota  u(s)+\kappa v(s))J^AK^C R{}_{AB[CD}\hat{K}{}_{E]} + w(s) J^AK^C R{}_{AB[CD} \hat{I}{}_{E]} =\beta(s)\hat{h}{}_{B[D} \hat{K}{}_{E]},
	\end{equation}
	where $w(s) = (s\iota-\kappa)u(s)+(s\kappa-\jmath)v(s)$, $\hat{K}_B = 2K^A I{}_{[A} J{}_{B]} $, $\hat{I}_B = 2I^A I{}_{[A} J{}_{B]} $, $\hat{h}{}_{BD}  = h{}_{BD} -\frac{1}{K_AK^A}K_BK_D$ and ${R}{}_{ABCD}$ is the full prolongation of $r{}_{bd} $ in the scale $\sigma$.
\end{lemma}

\begin{proof}
	We begin by calculating
	\[
		\hat{K}_B = 2K^A I{}_{[A} J{}_{B]}  = (\kappa-s\iota)J_B-(\jmath -s\kappa)I_B\quad \text{and} \quad \hat{I}_B = 2I^A I{}_{[A} J{}_{B]} =\iota J_B-\kappa I_B.
	\] 
	So we can express $I$ and $J$ as a linear combination of $\hat{K}$ and $\hat{I}$.
	\[
	\iota\hat{K}_B-(\kappa-s\iota)\hat{I}_B = (\kappa^2-\iota\jmath )I_B \quad \text{and} \quad \kappa\hat{K}_B-(\jmath -s\kappa)\hat{I}_B = (\kappa^2-\iota\jmath )J_B.
	\]
	So by \eqref{Eqn:BigR},
	\begin{align*}
		0&=(\kappa^2-\iota\jmath){X}{}^{B}Z{}^{C}{}_{c} Z{}^{D}{}_{d}  {X}{}^{E}\left(  u(s){J}{}^{A}{R}{}_{AB[CD}{I}{}_{E]}+v(s){J}{}^{A}{R}{}_{AB[CD}{J}{}_{E]}\right)\\ 
		 &=X{}^{B} Z{}^{C}{}_{c} Z{}^{D}{}_{d} X{}^{E} \big((\iota  u(s)+\kappa v(s))J^A R{}_{AB[CD}\hat{K}{}_{E]} + w(s) J^A R{}_{AB[CD} \hat{I}{}_{E]}\big).
	\end{align*}
	We have $K_A = \nabla^{\sigma}_a\tau Z{}^{a}{}_{A} \mod X$, so by the above,
	\[
	0 = X{}^{B} K{}^{C}  Z{}^{D}{}_{d} X{}^{E} \big((\iota  u(s)+\kappa v(s))J^A R{}_{AB[CD}\hat{K}{}_{E]} + w(s) J^A R{}_{AB[CD} \hat{I}{}_{E]}\big)
	\]
	Fix a value $s=\alpha$ such that $K_AK^A = \sigma^2\boldsymbol{g}(ds,ds) = \iota \alpha^2-2\alpha\kappa+\jmath \neq 0$, so the level set $\Sigma_\alpha = \{q\in M\setminus\mathcal{Z}(\sigma):s(q)=\alpha\} $ is a hypersurface in $M\setminus\mathcal{Z}(\sigma)$. Then consider the tractor 
	\[
	 V{}_{BDE}  = K{}^{C}  \big((\iota  u(s)+\kappa v(s))J^A R{}_{AB[CD}\hat{K}{}_{E]} + w(s) J^A R{}_{AB[CD} \hat{I}{}_{E]}\big)
	\] 
	Then $V{}_{BDE} \in K^{ \perp}$, since $J{}^{A} K{}^{B} R{}_{ABCD} =0$, which also gives
	\[
		4V{}_{[BDE]}  = K{}^{C}  \big((\iota  u(s)+\kappa v(s))J^A R{}_{A[BCD}\hat{K}{}_{E]} + w(s) J^A R{}_{A[BCD} \hat{I}{}_{E]}\big)=0.
	\] 
	However, $V{}_{BCD} $ is not trace-free. So we take the trace-free part of $V$ with respect to the metric on $K^{\perp}$, $\hat{h}{}_{BD}  = h{}_{BD} -\frac{1}{K_AK^A}K_BK_D$. That is, we let
	\[
		V{}^{\circ}{}_{BDE}  = V{}_{BDE} - \frac{2}{n}V{}_{A}{}^{A}{}_{[E} \hat{h}{}_{D]B},
	\]
	which is orthogonal to $K$ and trace-free for the standard tractor metric, since $V\in K^{\perp}$. So $V^{\circ}$ is in the irreducible tractor bundle $\begin{picture}(16, 6.5)(0, 1.5)
	\put(-0.2, 10){\line(1,0){10.4}}
	\put(-0.2, 5){\line(1,0){10.4}}
	\put(-0.2, 0){\line(1,0){5.4}}
	\put(0, 10){\line(0,-1){10}}
	\put(5, 10){\line(0,-1){10}}
	\put(10, 10){\line(0,-1){5}}
	\put(10, 0){${}_{\,0}$}
\end{picture}$. Since $V^{\circ}$ is a linear combination of parallel tractors with coefficients that are functions of $s$, it follows that $V^{\circ}$ is parallel along $\Sigma{}_{\alpha}$. Lastly, we have $X{}^{B} Z{}^{D}{}_{d} X{}^{E} V^{\circ }{}_{BDE} =0$. So, by Lemma \ref{Lem:SMBGG}, $V^{\circ}=0$ on each $\Sigma_\alpha$. By Lemma \ref{Lem:Regs} the hypersurfaces for regular values $s=\alpha$ fill a dense subset of $M\setminus\mathcal{Z}(\sigma)$. Hence, $V^{\circ }=0$ on $M\setminus\mathcal{Z}(\sigma)$.

	Finally, since $\hat{K}\in K^{\perp} $, we have
	\[
		0 = \hat{K}{}^{B} V^{\circ }{}_{BDE} =\frac{2}{n}V{}_{A}{}^{A}{}_{[E} \hat{h}{}_{D]B} \hat{K}^B=\frac{2}{n}V{}_{A}{}^{A}{}_{[E} \hat{K}{}_{D]} 
	\] 
	So $V{}_{A}{}^{A}{}_{E} \propto \hat{K}_E$. Then since $\hat{K}$ and $V{}_{A}{}^{A}{}_{E} $ are parallel along each $\Sigma{}_{\alpha} $, it follows that they are related by a function of $s$.
\end{proof}

Suppose that the subbundle spanned by $I$ and $J$ is non-degenerate, then there is a projection $\pi\colon \mathcal{E}_A\to I^{\perp}\cap J^{\perp} $. Recall the tractor $G{}_{ABCD}  = (h\owedge (II-\frac{\iota}{2}h))$ from Lemma \ref{Lem:Indfromg}. If $U_A\in I^{\perp}\cap J^{\perp}$, then
\[
J^AJ^C U^B G{}_{ABCD}  = \kappa^2 U_D- \iota\jmath U{}_{D} .
\] 
Hence, we can take $\pi{}_{BD}  = \frac{1}{\kappa^2-\iota\jmath }J^AJ^C G{}_{ABCD} $ as the parallel projection onto $I^{\perp}\cap J^{\perp}$. Additionally, we can calculate
\[
	J^A G{}_{ABCD} = 2\kappa h{}_{B[D} I{}_{C]} -2\iota h{}_{B[D} J{}_{C]} +2I_B J{}_{[C} I{}_{D]} = -2h{}_{B[D} \hat{I}{}_{C]} +2I_B J{}_{[C} I{}_{D]},
\] 
so $J^AG_{AB[CD}\hat I_{E]}=0$, since $\hat I\in\operatorname{span}(I,J)$. Contracting with $J^C$ gives
\[
	0= 3J{}^{A} J{}^{C} G{}_{AB[CD} \hat{I}{}_{E]} = (\kappa^2-\iota\jmath)\big(2\pi{}_{B[D} \hat{I}{}_{E]} - J^A G{}_{ABDE}\big),
\]
that is, $J^AG_{ABDE} = 2\pi_{B[D}\hat I_{E]}$.

\begin{theorem}\label{Thm:U'}
	Let $(M,\boldsymbol{c})$ be a conformally flat, connected manifold of dimension $n\ge 3$. Let $\sigma,\tau\in \Gamma(\mathcal{E}[1])$ be almost Einstein scales and let ${I}{}_{A}$, ${J}{}_{A}$ be the corresponding parallel tractors. Let $s = \tau\sigma^{-1}$ and let $\tilde{\sigma} = f(s)\sigma$ for some smooth $f\colon \mathbb{R}\to \mathbb{R} $ such that $\tilde{\sigma}$ is non-zero on an open dense subset of $M$.

	Suppose that the space spanned by $I, J$ is non-degenerate, or equivalently $\kappa^2-\iota\jmath \neq 0$. If the Hamiltonian
	\begin{equation}\label{Eqn:HB:1'}
		\mathcal{H}(q,p) = g^{-1}_q(p,p)+\frac{1}{f(s)^2} = g^{-1}_q(p,p)+V(s),
	\end{equation}
	where $g = \sigma^{-2}\boldsymbol{g}$ and $V(s) = f(s)^{-2}$, is second-order maximally superintegrable, then $f$ solves either
	\begin{tabel}
		\item\label{Thm:U':1} $\omega(s) = \iota u(s)+\kappa v(s)=0$ or
		\item\label{Thm:U':2} $w(s) = (s\iota-\kappa)u(s)+(s\kappa-\jmath)v(s)=0$,
	\end{tabel}
	where $u(s)$ and $v(s)$ are as in Lemma \ref{Lem:Big}.
\end{theorem}
\begin{proof}
	Since the conclusion is local, we may replace $M$ by a simply connected open subset, and so assume that $(M,\boldsymbol{c})$ has trivial conformal holonomy.
	Suppose that $\mathcal{H}$ is second-order maximally superintegrable. Nondegeneracy of the subspace spanned by $I$ and $J$ implies that $(\iota,\kappa)\neq (0,0)$, so by Theorem \ref{Thm:2Ein} and Lemma \ref{Lem:U+} there is a parallel $R{}_{ABCD} $ with curvature symmetries in $I^{\perp}$ such that $R{}_{ABCD} -\lambda G{}_{ABCD} \not\in J^{\perp}$ for all $\lambda \in  \mathbb{R}$, and such that \eqref{Eqn:U+} holds
	\begin{equation}\label{Eqn:U':1}
	(\iota  u(s)+\kappa v(s))J^AK^C R{}_{AB[CD}\hat{K}{}_{E]} + w(s) J^AK^C R{}_{AB[CD} \hat{I}{}_{E]} =\beta(s)\hat{h}{}_{B[D} \hat{K}{}_{E]},
	\end{equation} 
	where $K_A = J_A-sI_A$, $\hat{K}_B = 2K^A I{}_{[A} J{}_{B]} $, $\hat{I}_B = 2I^A I{}_{[A} J{}_{B]} $, and $\hat{h}{}_{BD}  = h{}_{BD} - \frac{1}{K_AK^A}K_BK_D$ for all points $q\in M$ such that $K_AK^A = \iota s^2-2\kappa s+\jmath \neq 0$. The form $2I{}_{[A} J{}_{B]} $ acts like a volume element on the subspace spanned by $I$ and $J$, with determinant $\iota\jmath -\kappa^2$. So the square length of $\hat{K}$ is $(\iota\jmath -\kappa^2)(\iota s^2-2\kappa s+\jmath)\neq 0$. We rewrite \eqref{Eqn:U':1} using $K^A \hat{K}_A=0$ and $I^A \hat{I}_A=0$. This gives
	\begin{equation}\label{Eqn:U':1a}
		2(\iota u(s)+\kappa v(s))J^AK^C R{}_{ABC[D} \hat{K}{}_{E]} +3w(s)J^AJ^C R{}_{AB[CD} \hat{I}{}_{E]}  = 3\beta(s)\hat{h}{}_{B[D} \hat{K}{}_{E]}.
	\end{equation} 

	Observe that $J^AK^C R{}_{AB[CD} \hat{I}{}_{E]} $ is orthogonal to $K$ and $I$ in all indices, so it is orthogonal to the span of $I$ and $J$, including $\hat{K}$. Let $U{}_{BD} = 3\beta(s)\hat{h}{}_{BD} -2(\iota u(s)+\kappa v(s))J^AK^C R{}_{ABCD} $ so that
	\[
		3w(s)J^AJ^C R{}_{AB[CD} \hat{I}{}_{E]} = 2U{}_{B[D} \hat{K}{}_{E]} \in \hat{K}^{\perp}.
	\] 
	Then $\hat{K}{}^{B} U{}_{B[D} \hat{K}{}_{E]} =0$, so $\hat{K}{}^{B} U{}_{BD} \propto \hat{K}_D$. Then contracting $\hat{K}^E$ into the above gives
	\[
		0 = \hat{K}_E \hat{K}^E U{}_{BD} -\hat{K}^EU{}_{BE}  \hat{K}_D\implies U{}_{BD} \propto \hat{K}_B \hat{K}_D,
	\]
	since $\hat{K}_A \hat{K}^A\neq 0$ and $\hat{K}{}^{B} U{}_{BD} \propto \hat{K}_D$. Hence,
	\[
		2U{}_{B[D} \hat{K}{}_{E]} = 3w(s) J^AJ^C R{}_{AB[CD} \hat{I}{}_{E]}=0
	\] 
	So either $w(s)=0$ on all of $M\setminus\mathcal{Z}(\sigma)$ or there is a point of $M\setminus\mathcal{Z}(\sigma)$ such that
	\[
		J^AJ^C R{}_{AB[CD} \hat{I}{}_{E]} = \iota J{}^{A} J{}^{C} R{}_{AB[CD}  J{}_{E]} -\kappa J{}^{A} J{}^{C} R{}_{AB[CD}  I{}_{E]}=0,
	\] 
	which then holds on all of $M$ since it is parallel and $M$ is connected. Then, by \eqref{Eqn:U':1a},
	\[
	2(\iota  u(s)+\kappa v(s))J^AJ^C R{}_{ABC[D}\hat{K}{}_{E]} = 3\beta(s)\hat{h}{}_{B[D} \hat{K}{}_{E]}
	\]
	Now suppose that $(\iota  u(s)+\kappa v(s))\neq 0$ at some point $q\in M$. Then at $q$, there are $c_{1}$ and $c_{2}$ such that
	\[
	J^AJ^C R{}_{ABCD}  = c_{1}\hat{h}{}_{BD} + c_{2}\hat{K}_B \hat{K}_D
	\] 
	at $q$. The right hand side maps $U\in I^{\perp}\cap J^{\perp}$ to $c_{1}U$, while the left hand side is orthogonal to $I$ and $J$. This completely determines $J^AJ^C R{}_{ABCD} \in \operatorname{Hom}(\mathcal{E}^B,\mathcal{E}_D)|_q$. So, at $q$,
	\[
	J^AJ^C R{}_{ABCD} = c_{1}\pi{}_{BD},
	\] 
	where $\pi{}_{BD} =\frac{1}{\kappa^2-\iota\jmath }J^AJ^C G{}_{ABCD} $ is the parallel projection $\pi\colon \mathcal{E}_B\to I^{\perp}\cap J^{\perp} $, then this holds on all of $M$ since both tractors are parallel. Combining our two conclusions gives
	\[
		0=3J^AJ^C R{}_{AB[CD} \hat{I}{}_{E]} =2c_{1}\pi{}_{B[D} \hat{I}{}_{E]} + (\iota\jmath -\kappa{}^{2})J^A R{}_{ABDE} 
	\] 
	\[
		\implies J^A R{}_{ABDE} = -\frac{2c_{1}}{\iota\jmath -\kappa^2}\pi{}_{B[D} \hat{I}{}_{E]}  = -\frac{c_{1}}{{\iota\jmath -\kappa^2}}J^A G{}_{ABDE} .
	\] 
	However, this contradicts the assumptions on $R$ from Theorem \ref{Thm:2Ein}, that is, the equivalent conditions of Lemma \ref{Lem:Indfromg}, since $R+\frac{c_{1}}{{\iota\jmath -\kappa^2}}G$ is orthogonal to $J$. 
\end{proof}
\begin{corollary}\label{Cor:U}
	Let $(M,\boldsymbol{c})$ be a conformally flat, simply connected manifold of dimension $n\ge 3$. Let $\sigma,\tau\in \Gamma(\mathcal{E}[1])$ be almost Einstein scales and let ${I}{}_{A}$, ${J}{}_{A}$ be the corresponding parallel tractors. Let $s = \tau\sigma^{-1}$ and let $\tilde{\sigma} = f(s)\sigma$ for some smooth $f\colon \mathbb{R}\to \mathbb{R} $ such that $\tilde{\sigma}$ is non-zero on an open dense subset of $M$. Suppose that the subspace spanned by $I$ and $J$ is non-degenerate. Then without loss of generality, we consider two cases (see Theorem \ref{Thm:Ex1}). 

	Suppose that $I$ and $J$ are null and that ${I}{}_{A}{J}{}^{A}=1$. Then the Hamiltonian
	\begin{equation}\label{Eqn:HB:1}
		\mathcal{H}(q,p) = g^{-1}_q(p,p)+\frac{1}{f(s)^2} = g^{-1}_q(p,p)+V(s),
	\end{equation}
	where $g = \sigma^{-2}\boldsymbol{g}$ and $V(s) = f(s)^{-2}$, is second-order maximally superintegrable if and only if $f$ solves either
	\[
		v(s) = \frac{f(s)^4}{2}V''(s) = 0, \quad \text{in which case}\quad V(s) = as+b
	\]
	or
	\[
		u(s)-sv(s) = -\frac{f(s)^4}{2}(2sV''(s)+3 V'(s)) = 0, \quad \text{in which case}\quad V(s) = \frac{a}{\sqrt{|s|} }+b
	\] 

	Let $\iota=I_AI^A$ and $\jmath  = J_AJ^A$. Suppose that $I$ and $J$ are orthogonal and $\iota,\jmath \in \{-1,+1\} $. Then the Hamiltonian
	\begin{equation}\label{Eqn:HB:2}
		\mathcal{H}(q,p) = g^{-1}_q(p,p)+\frac{1}{f(s)^2} = g^{-1}_q(p,p)+V(s),
	\end{equation}
	is second-order maximally superintegrable if and only if $f$ solves either
	\[
		u(s) = -\frac{f(s)^4}{2}(sV''(s)+3V'(s)) = 0, \quad \text{in which case}\quad V(s) = \frac{a}{s^2}+b
	\] 
	or
	\[
		v(s)-\iota\jmath su(s)=\frac{f(s)^4}{2}((\iota\jmath s^2+1)V''(s)+3\iota\jmath sV'(s))=0,\quad \text{in which case}\quad V(s) = a \frac{s}{\sqrt{|\iota\jmath s^2+1|} }+b.
	\] 
\end{corollary}The choice of parameter for our potential, $V(s)$ is inspired by two classical results. The first being the second-order analogue of Bertrand’s theorem, which becomes a special case, and the second is Evans' classification of second-order maximally superintegrable systems in three dimensions \cite{evans_superintegrability_1990}; see Section \ref{Sec:FM}. 

Consider all possible scale tractors which are combinations of ${I}{}_{A}$, ${J}{}_{A}$ and ${X}{}_{A}$, with functional coefficients. That is, a scale density $\tilde{\sigma}$ such that 
\[
f_{1}{I}{}_{A}+f_{2}{J}{}_{A} = {D}{}_{A}\tilde{\sigma} \mod {X}{}_{A}
\]
Contracting ${X}{}^{A}$ implies
\[
f_{1}\sigma+f_{2} \tau = \tilde{\sigma},
\] 
and instead contracting ${Z}{}^{A}{}_{a}$, the projector in the scale $\sigma$, gives
\[
f_{2}{\nabla}{}^{\sigma}{}_{a}\tau = {\nabla}{}^{\sigma}{}_{a}\tilde{\sigma} = \sigma (df_{1})_a+\tau(df_{2})_a+f_{2}{\nabla}{}^{\sigma}{}_{a}\tau
\] 
Therefore, $\sigma(df_{1})_a+\tau(df_{2})_a=0$ and, away from the zero set of $\sigma$, $(df_{1})_a = -s (df_{2})_a$. Taking the exterior derivative gives $0={(ds)}{}_{[b}{(df_{2})}{}_{a]}$, so both $df_{1}$ and $df_{2}$ are multiples of $ds$. Away from the zero set of $ds$, which is nowhere dense, this implies that $f_{1}$ and $f_{2}$ can be written as functions of $s$. This implies
\[
\tilde{\sigma} = (f_{1}(s)+sf_{2}(s))\sigma = f(s)\sigma,
\] 
so the class of potentials we have considered in this section are general in the sense of taking linear combinations of $I, J$ and $X$ with functional coefficients.

\section{Families of Hamiltonians}\label{Sec:Ham}
Superintegrable systems are usually presented as linear families of potentials \cite{evans_superintegrability_1990,kress_algebraic_2024,miller_classical_2013,verrier_new_2008}. One classical example is the isotropic harmonic oscillator on $\mathbb{R}^n$,
\begin{equation}\label{Eqn:IO}
	\mathcal{H} = \frac{1}{2}|p|^2+\omega^2|q|^2+\frac{\lambda_{1}}{q_{1}^2}+\dots + \frac{\lambda_n}{q_n^2} + \lambda_{0}.
\end{equation}
A $(N+1)$ \textit{linear family of Hamiltonians} is a collection of functions on $T^*M$ parameterised by $\lambda = (\lambda_{0},\lambda_{1},\dots ,\lambda_N)\in \mathbb{R}^{N+1}$ of the form
\begin{equation}\label{Eqn:HL}
	\mathcal{H}^\lambda(q,p) = g_q^{-1}(p,p)+\lambda_1 V_{1}(q)+\dots +\lambda_N V_N(q)+\lambda_{0},
	\end{equation} 
	where $g$ is a fixed metric on $M$. A $N+1$ \textit{linear family of quadratic constants of motion} for $\mathcal{H}^{\lambda}$ is a collection of functions on $T^*M$ parameterised by $\lambda = (\lambda_{0},\lambda_{1},\dots ,\lambda_N)\in \mathbb{R}^{N+1}$ of the form
	\begin{equation}\label{Eqn:RL}
	r^{\lambda}(q,p) = r_q(p,p)+\lambda_{1} W_{1}(q)+\dots +\lambda_N W_N(q) +c(\lambda)
	\end{equation} 
	where \(c\) is an affine function of the parameters, such that $\left\{\mathcal{H}^{\lambda},r^{\lambda}\right\}=0$ for all $\lambda\in \mathbb{R}^{N+1}$. 
	
	The dependence of the family of Hamiltonians and constants on $\lambda$ is often suppressed in the symbols $\mathcal{H}$ and $r$.
\begin{theorem}\label{Thm:HFam}
	Let $(M,\boldsymbol{c})$ be a conformal manifold, let $\sigma\in \Gamma(\mathcal{E}[1])$ be a scale, let $g = \sigma^{-2}\boldsymbol{g}$, and let $V_{1},\dots ,V_N$ be potentials and consider the family $\mathcal{H}^{\lambda}$ \eqref{Eqn:HL}. Then there is a linear isomorphism between $N+1$ parameter families of quadratic constants of motion for $\mathcal{H}^\lambda$, modulo constant multiples of \(\mathcal H^\lambda\) and additive affine functions of \(\lambda\), and rank two conformal Killing tensors for which $g,V_{1}g,\dots ,V_N g$ are all Killing scales. Linear constants of motion for $\mathcal{H}^{\lambda}$ are in one-to-one correspondence with conformal Killing vectors for which $g,V_{1}g,\dots ,V_N g$ are all Killing scales.
\end{theorem}
\begin{proof}
	Define a linear map $\Phi$ from the space of linear families of quadratic constants of motion as in \eqref{Eqn:RL} to the space of conformal Killing tensors given by $\Phi(r^{\lambda}) = \operatorname{trace-free}(r)$. Fixing some family $r^{\lambda}$, taking $\lambda_i=1$ and leaving the other parameters zero for each $i$ shows that 
	\[
	\left\{g^{-1}+V{}_{i} ,r+W_i\right\}=0
	\] 
	for each $i\in \{1,\dots ,N\} $. So by Theorem \ref{Thm:QCns}, each of $V_{1}g,\dots ,V_Ng$ are Killing scales for $k$. 

	Conversely, suppose that $k$ be a conformal Killing tensor for each of of $g,V_{1} g,\dots ,V_N g$ are Killing scales. Then by Theorem \ref{Thm:QCns}, there exist functions $W_{1},\dots ,W_N$ such that for each $i$,
	\[
		\left\{g^{-1}+V{}_{i} ,r+W_i\right\} = \left\{g^{-1},r\right\}+(\left\{V{}_{i} ,r\right\}+\left\{g^{-1},W{}_{i} \right\})=0,
	\] 
	with terms of each order vanishing independently. So
	\[
		\left\{\mathcal{H}^{\lambda},r^{\lambda}\right\} = \left\{g^{-1},r\right\}+\sum_{i=1}^{N} \lambda_i\left( \left\{V{}_{i} ,r\right\}+\left\{g^{-1},W{}_{i} \right\} \right) =0
	\]
	So $\Phi$ is a surjective linear map from the space of linear families of constants of motion to conformal Killing tensors for which $g,V_{1}g,\dots ,V_Ng$ are Killing scales. A linear family of constants of motion in the kernel of $\Phi$ must have $r = ag^{-1}$ for some $a\in \mathbb{R}$ and $dW_i = a\,dV_i$ for each $i$. Hence, $r^{\lambda} = a\mathcal{H}^{\lambda}+c(\lambda)$, where $c$ is some affine function of $\lambda$.

	Linear constants of motion can have no trace adjustment. If $k$ is a conformal Killing vector which induces a constant of motion for $\mathcal{H}$, then fixing each $\lambda_i=1$ and others zero shows that $g,V_{1}g,\dots ,V_Ng$ are Killing scales by Proposition \ref{Prp:PBKT}. The same Proposition combined with linearity of $\left\{,\right\}$ shows that any conformal Killing vector for which $g,V_{1}g,\dots ,V_Ng$ are Killing scales induces a linear constant of motion for $\mathcal{H}$
\end{proof}
The linearity observed in Theorem \ref{Thm:HFam} leads to an unusual result for Killing scales. The general term for a metric $\sigma^{-2}\boldsymbol{g}$, where $\sigma\in \Gamma(\mathcal{E}[1])$ is a scale is an \textit{almost metric}.
\begin{corollary}\label{Cor:AddM}
	Let $(M,\boldsymbol{c})$ be a conformal manifold, let $k$ be a conformal Killing tensor of rank one or two, and let $g$ be a Killing scale for $k$. Then the set of functions $V$ such that $Vg$ is a Killing scale for $k$ is a linear subspace. In other words, being a Killing scale for $k$ is $\mathbb{R}$-linear in the almost metric.
\end{corollary}
\begin{proof}
	For rank one, by Proposition \ref{Prp:PBKT}, $Vg$ is a Killing scale if and only if $\{V,k\}=0$, which is linear in $V$. For rank two, let $r=k+Tg^{-1}$ be Killing for $g$. By Theorem \ref{Thm:QCns}, $Vg$ is a Killing scale if and only if there is a function $U$ with $\{V,r\}+\{g^{-1},U\}=0$, and this condition is linear in the pair $(V,U)$.
\end{proof}Due to Corollary \ref{Cor:AddM}, families of constants of motion for a family of Hamiltonians
	\[
		\mathcal{H} = g^{-1}_q(p,p)+\lambda_{1}V_{1}+\dots +\lambda_NV_N + \lambda_0.
	\] 
	correspond to conformal Killing tensors for which every nonzero almost metric in the subspace spanned by $g,V_{1}g,\dots ,V_Ng$ is a Killing scale. If we choose a nonzero $\tilde{V}g = \tilde{g}\in \operatorname{span}\left( g,V_{1}g,\dots ,V_Ng \right)$, we can form a new family of Hamiltonians,
	\begin{equation}\label{Eqn:Htilde}
		\tilde{\mathcal{H}} = \tilde{g}_q^{-1}(p,p)+\lambda_{1}\frac{V_{1}}{\tilde{V}}+\dots +\lambda_N \frac{V_N}{\tilde{V}}+\lambda_{0}\frac{1}{\tilde{V}}.
	\end{equation}
	The constant term in this new Hamiltonian is hidden, since $\tilde{V}$ is a linear combination of $V_{1},\dots ,V_N$. Since $\operatorname{span}(\tilde{g},(V_{1} / \tilde{V})\tilde{g},\dots ,(V_N / \tilde{V})\tilde{g},g) =\operatorname{span}(g,V_{1}g,\dots ,V_N g)$, it follows from Theorem \ref{Thm:HFam} that there is a linear isomorphism between families of constants of motion for $\mathcal{H}$ and families of constants of motion for $\tilde{\mathcal{H}}$.
	
	A second-order maximally superintegrable family is called \textit{non-degenerate} if $N=n+1$, that is, if it has $n+2$ parameters \cite[Definition 4.2]{kress_algebraic_2023}. An example is \eqref{Eqn:IO}. 

It then follows that (maximal) quadratic superintegrability of ${\mathcal{H}}$ and that of $\tilde{\mathcal{H}}$ are equivalent. One example is the superintegrable system given in \cite{ballesteros_maximally_2008}.

We now show how our main results, Theorems \ref{Thm:QMS}, \ref{Thm:Ex1} and \ref{Thm:Ex2} can be extended to quasimaximally and second-order maximally superintegrable families of Hamiltonians.

\begin{theorem}\label{Thm:QMSFam}
	Let $(M,\boldsymbol{c})$ have trivial conformal holonomy. Let ${I}{}_{A} = D_A \sigma$ and $J_A = D_A \tau$ for $\sigma,\tau\in \Gamma(\mathcal{E}[1])$ be linearly independent parallel tractors such that $\operatorname{span}(I,J)$ is non-degenerate, let $s=\tau\sigma^{-1}\in C^\infty(M\setminus\mathcal{Z}(\sigma))$ and let $I^1_A,\dots ,I^n_A$ be a parallel pairwise orthogonal basis for $I^{\perp}\cap J^{\perp}$. Let $F\colon \mathbb{R}\to \mathbb{R} $ be smooth and consider the linear family of Hamiltonians
	\begin{equation}\label{Eqn:TilH}
	\mathcal{H}(q,p) = g^{-1}_q(p,p)+\lambda F\left( s \right) + \lambda_{1}\frac{\sigma^2}{\sigma_{1}^2} +\dots +\lambda_n \frac{\sigma^2}{\sigma_n^2} +\lambda_{0}
	\end{equation} 
	where $g = \sigma^{-2}\boldsymbol{g}$ and $\sigma_i = X^AI^i_A$ for $i\in \{1,\dots ,n\} $. Then $\mathcal{H}$ has $2n-2$ functionally independent families of constants of motion; the Hamiltonian and \(2n-3\) further families whose quadratic parts are symmetric squares of vector fields which are Killing for $g$ and $\hat{g} = \tau^{-2}\boldsymbol{g}$.
\end{theorem}

\begin{proof}
	For each $i,j\in \{1,\dots ,n\} $, let $K^{i,j}_{AB} = 2I{}^{i}{}_{[A} I{}^{j}{}_{B]} \in I^{\perp}\cap J^{\perp}$ and let $k{}^{i,j}{}_{b}  = X{}^{A} Z{}^{B}{}_{b} K^{i,j}_{AB} $. By Theorem \ref{Thm:KT} \ref{Prp:ConPro:1:1}, each $k^{i,j}_a$ is Killing for $g$ and $\hat{g}$. Then, by Proposition \ref{Prop:fKS}, each $k^{i,j}_a$ is Killing for $F(s)g$. Let ${r}^{i,j}{}_{ab} = {k}^{i,j}{}_{(a}{k}^{i,j}{}_{b)}$. Then each $r^{i,j}$ is Killing for $g$, $\hat{g}$, and $F(s)g$. Let 
	\begin{equation}
		\begin{aligned}
			{R}{}^{i,j}{}_{ABCD} = (I^{i\,2}\owedge I^{j\,2}) 
			&= {I}{}^{i}{}_{A}{I}{}^{j}{}_{B}{I}{}^{i}{}_{C}{I}{}^{j}{}_{D}+{I}{}^{j}{}_{A}{I}{}^{i}{}_{B}{I}{}^{j}{}_{C}{I}{}^{i}{}_{D}-{I}{}^{i}{}_{A}{I}{}^{j}{}_{B}{I}{}^{j}{}_{C}{I}{}^{i}{}_{D}-{I}{}^{j}{}_{A}{I}{}^{i}{}_{B}{I}{}^{i}{}_{C}{I}{}^{j}{}_{D}\\
			&= K^{i,j}_{AB}K^{i,j}_{CD}
		\end{aligned}
	\end{equation}
	for each $i,j\in \{1,\dots ,n\} $. So, $X{}^{A} Z{}^{B}{}_{b} X{}^{C} Z{}^{D}{}_{d} R{}^{i,j}{}_{ABCD} = k^{i,j}{}_{(a}k^{i,j}{}_{b)} = r{}^{i,j}{}_{bd} $. Also
\[
	{I}{}^{\ell\,A}{R}{}^{i,j}{}_{AB[CD}{I}{}^{\ell}{}_{E]}=0.
\] 
for all $i,j,\ell\in \{1,\dots ,n\} $. So, by Theorem \ref{Thm:IRI}, $g_\ell = \sigma_\ell^{-2}\boldsymbol{g} = \frac{\sigma^2}{\sigma_\ell^2}g$ is a Killing scale for the trace-free part of $r{}^{i,j}{}_{ab} $ for all $i,j,\ell\in \{1,\dots ,n\} $. So, by Theorem \ref{Thm:HFam}, for each $i,j\in \{1,\dots ,n\} $ there is a linear family of constants of motion for the family of Hamiltonians $\mathcal{H}^{\lambda}$ with quadratic part $r^{i,j}$. Then by Proposition \ref{Prp:Find} and Lemma \ref{Lem:gind}, $g^{-1},r^{1,2},\dots ,r^{1,n}, r^{2,n},\dots ,r^{n-1,n}$ are functionally independent, so the corresponding families of constants of motion are functionally independent by Proposition \ref{Prp:Ind}.\end{proof}
\begin{theorem}\label{Thm:MSFam}
	Let $(M,\boldsymbol{c})$ have trivial conformal holonomy. Let $\sigma,\tau\in \Gamma(\mathcal{E}[1])$ be almost Einstein scales and let $I_A = D_A \sigma$ and $J_A = D_A \tau$ be the corresponding parallel tractors. Let $s = \tau\sigma^{-1}\in C^\infty (M\setminus\mathcal{Z}(\sigma))$. Suppose that the subspace spanned by $I$ and $J$ is non-degenerate, then without loss of generality we restrict to cases \ref{Thm:Ex1:1:IJ} and \ref{Thm:Ex1:2:IJ} of Theorems \ref{Thm:Ex1} and \ref{Thm:Ex2}.

	If $V(s)$ is the potential from the corresponding case of Theorem \ref{Thm:Ex1} and $I^1,\dots ,I^n$ is a parallel pairwise orthogonal basis for $I^{\perp}\cap J^{\perp}$, then the $n+2$ parameter linear family of Hamiltonians
	\begin{equation}\label{Eqn:MSFam:1}
		\mathcal{H}^{\lambda}(q,p) = g^{-1}_q(p,p)+\lambda V(s) +\lambda_{1}V_{1}(q)+\dots +\lambda_n V_{n}(q)+\lambda_{0},
	\end{equation} 
	where $\sigma_i = X{}^{A} I{}^{i}{}_{A} $ and $V_i = \frac{\sigma^2}{\sigma_i^2} $, is second-order maximally superintegrable.

	If $V(s)$ is the potential from the corresponding case of Theorem \ref{Thm:Ex2} and $I^1,\dots ,I{}^{n-1} $ is a parallel pairwise orthogonal linearly independent subset of $I^{\perp}\cap J^{\perp}$, then the $n+1$ parameter linear family of Hamiltonians
	\begin{equation}\label{Eqn:MSFam:2}
		\mathcal{H}^{\lambda}(q,p) = g^{-1}_q(p,p)+\lambda V(s) +\lambda_{1}V_{1}(q)+\dots +\lambda{}_{n-1}  V_{n-1}(q)+\lambda_{0},
	\end{equation} 
	where $\sigma_i = X{}^{A} I{}^{i}{}_{A} $ and $V_i = \frac{\sigma^2}{\sigma_i^2} $, is second-order maximally superintegrable.
\end{theorem}
\begin{proof}
	For the cases contained in \eqref{Eqn:MSFam:1} we apply Theorem \ref{Thm:QMSFam} with $F = V$ for the first $2n-2$ constants of motion. Then we choose the final constant of motion such that $\sigma_{1},\dots ,\sigma_n$ are Killing scales.

	For the case in Theorem \ref{Thm:Ex1} \ref{Thm:Ex1:1:IJ}, we take $R = (I^2\owedge (I^i)^2)\in I^{\perp}$ for any $i\in \{1,\dots ,n\} $ for the final Killing tensor, which satisfies
	\[
		R{}_{AB[CD} \pi{}_{E]}{}^{A} = 0  \quad \text{and} \quad I{}^{j}{}^{A} R{}_{AB[CD}  I{}^{j}{}_{E]}  = 0
	\] 
	for all $j\in \{1,\dots ,n\} $.

	Similarly, in the case in Theorem \ref{Thm:Ex1} \ref{Thm:Ex1:2:IJ}, we can take $R = (J^2\owedge (I^i)^2)$ for any $i\in \{1,\dots ,n\} $.

	For the cases contained in \eqref{Eqn:MSFam:2}, we first suppose that the subspace spanned by $I^{1},\dots ,I{}^{n-1} $ is non-degenerate. Then, since $I^{\perp}\cap J^{\perp}$ is non-degenerate, let $I{}^{n} $ span the space of tractors orthogonal to $I^{1},\dots ,I{}^{n-1} $ inside $I^{\perp}\cap J^{\perp}$. Then $I^1,\dots ,I^n$ are pairwise orthogonal and we proceed similarly to \eqref{Eqn:MSFam:1}.

	For the case in Theorem \ref{Thm:Ex2} \ref{Thm:Ex2:1:IJ}, we take $R = ((I \vee I^n)\owedge \pi)\in I^{\perp}$. Then
	\[
		J^A R{}_{ABCD} = 2\pi{}_{B[D} I{}^{n}{}_{C]} \quad \text{and} \quad I{}^{j}{}^{A} R{}_{AB[CD} I{}^{j}{}_{E]} = 0
	\] 
	for all $ j\neq n$, so we can add the first $n-1$ potential terms as in Theorem \ref{Thm:QMSFam}.
	
	Similarly, for the case in Theorem  \ref{Thm:Ex2} \ref{Thm:Ex2:2:IJ}, we take $R = ((J \vee I^n)\owedge \pi)$.
	
	Now suppose that the subspace spanned by $I^{1},\dots ,I{}^{n-1} $ is degenerate. Then exactly one of the $I^1,\dots ,I^{n-1}$ is null, and spans the space of vectors in $I^{\perp}\cap J^{\perp}$ which are orthogonal to $I^{1},\dots ,I{}^{n-1} $. This is because the $I^{1},\dots ,I{}^{n-1} $ are pairwise orthogonal and $I^{\perp}\cap J^{\perp}$ is non-degenerate. Suppose $I{}^{1} $ is null and choose a null $I^n\in I^{\perp}\cap J^{\perp}$ such that $h(I^1,I^n)=1$ and $h(I^i,I^n)=0$ for $i\in \{2,\dots ,n-1\} $. 

	As in the proof of Theorem \ref{Thm:QMSFam}, curvature tensors of the form
	\[
		R^{i,j}_{ABCD} = 4I{}^{i}{}_{[A} I{}^{j}{}_{B]} I{}^{i}{}_{[C} I{}^{j}{}_{D]}
	\]
	for $i,j\in \{1,\dots ,n-1\} $ and $i=1,j=n$ are compatible with $I^{1},\dots ,I^{n-1}$. But $R^{i,n}$ for $i>1$ is not compatible with $I{}^{1}$. So instead we take
	\[
		R{}^{i}{}_{ABCD}  = 2I{}^{1}{}_{[A} I{}^{i}{}_{B]} I{}^{n}{}_{[C} I{}^{i}{}_{D]} +2I{}^{n}{}_{[A} I{}^{i}{}_{B]} I{}^{1}{}_{[C} I{}^{i}{}_{D]}
	\]
	for $i\in \{2,\dots ,n-1\} $ which satisfies
	\[
		\quad I{}^{j}{}^{A} R{}^{i}{}_{AB[CD} I{}^{j}{}_{E]} =0.
	\]
	for all $j\in \{1,\dots ,n-1\} $. Hence, there are families of quadratic constants of motion  with quadratic part 
	\[
		r{}^{i,j}{}_{bd}  = X{}^{A} Z{}^{B}{}_{b} X{}^{C} Z{}^{D}{}_{d} R{}^{i,j}{}_{ABCD} = k^{i,j}{}_ak{}^{i,j}{}_{b}\quad \text{and} \quad \]\[r{}^{i}{}_{ab}  =  X{}^{A} Z{}^{B}{}_{b} X{}^{C} Z{}^{D}{}_{d} R{}^{i}{}_{ABCD} = k{}^{1,i}{}_{(a} k{}^{i,n}{}_{b)} 
	\] 
	The functions induced by $g^{-1},k^{1,2},k^{1,3},\dots ,k^{1,n},k^{2,n},\dots ,k^{n-1,n}$ are functionally independent by Proposition \ref{Prp:Find}, which is equivalent to functional independence of the functions induced by $g^{-1},r^{1,2},r^{1,3},\dots ,r^{1,n},r^2,\dots ,r^{n-1}$, whence we are again in the position of Theorem \ref{Thm:QMSFam} and need only check the final constant of motion. For the case in Theorem \ref{Thm:Ex2} \ref{Thm:Ex2:1:IJ} we take $R = ((I \vee I^1)\owedge \pi)$ and for the case in Theorem \ref{Thm:Ex2} \ref{Thm:Ex2:2:IJ}, we take $R = ((J \vee I^1)\owedge \pi)$.
\end{proof}The results derived here match classical results for the Harmonic oscillator and Kepler system, where to maintain quadratic superintegrability, they can become $n+2$ and $n+1$ parameter families respectively. However, the Kepler system can also include the final potential term when the final constant of motion is allowed to be quartic in the momenta \cite{verrier_new_2008}. It would be interesting to see how this constant arises in tractors.

\section{Calculations}\label{Sec:FM}

\subsection{Tractors to Proper Constants of Motion}\label{Sec:TrCM}

Throughout our theoretical development of superintegrable systems through tractors, we have been using the equivalences between constants of motion, Killing tensors and parallel tractors and proving superintegrability by constructing tractors satisfying certain requirements. Here we give an overview of how to recover explicit formulae in coordinates for constants of motion, starting with the corresponding tractors.

We recap what has already been explained in previous sections. If $I_A = D_A\sigma$, then the scale density is recovered by taking $\sigma = X^AI_A$, and the corresponding metric is $g = \sigma^{-2}\boldsymbol{g}$. If $I$ is parallel and $R{}_{ABCD}\in \Gamma(\begin{picture}(10, 6.5)(0, 2)
	\put(-0.2, 10){\line(1,0){10.4}}
	\put(-0.2, 5){\line(1,0){10.4}}
	\put(-0.2, 0){\line(1,0){10.4}}
	\put(0, 10){\line(0,-1){10}}
	\put(5, 10){\line(0,-1){10}}
	\put(10, 10){\line(0,-1){10}}
\end{picture}(\mathcal{E}_A))$ or $K{}_{AB} \in \Gamma(\mathcal{E}{}_{[AB]} )$ are parallel tractors orthogonal to $I$, then the Killing tensors are recovered by taking the invariant top slots,
\[
r{}_{bd}  = X{}^{A} Z{}^{B}{}_{b} X{}^{C} Z{}^{D}{}_{d} R{}_{ABCD} \quad \text{and} \quad k_b = X{}^{A} Z{}^{B}{}_{b} K{}_{AB}.
\] 
Now suppose that $V$ is a potential and consider the Hamiltonian
\[
\mathcal{H}(q,p) = g^{-1}_q(p,p)+V(q) = \sigma^2\boldsymbol{g}^{-1}_q(p,p)+V(q).
\]
If $Vg$ is a Killing scale for $k{}_{b}$, then the corresponding linear constant of motion given by Proposition \ref{Prp:PBKT} is just $f(q,p) = k_q(p)$. For quadratic constants of motion, we will use Lemma \ref{Lem:COM}. In particular, for Einstein scales, represented by $V(s) = s^{-2}$, the constant of motion is given by the following proposition.

\begin{proposition}\label{Prp:Comp}
	Let $(M,\boldsymbol{c})$ be conformally flat. Let $R{}_{ABCD} $ be a parallel tractor with curvature symmetries. Let $I_A = D_A\sigma $ and $J_A = D_A\tau $ be parallel scale tractors and let $g = \sigma^{-2}\boldsymbol{g}$. Suppose that $I^A R{}_{ABCD} =0$ and that 
	\[
		J^A R{}_{ABCD}  = 2S{}_{B[D} J{}_{C]} 
	\] 
	let $r{}_{bd}  = X{}^{A} Z{}^{B}{}_{b} X{}^{C} Z{}^{D}{}_{d} R{}_{ABCD} $. Let $\phi = \sigma^{-2}X^BX^D S{}_{BD} $. Then the function on $T^*M$ given by $f(q,p) = r_q(p,p)+ \frac{\lambda\phi(q)}{s^2}$ is a constant of motion for the Hamiltonian $\mathcal{H}(q,p) = g^{-1}_q(p,p)+\frac{\lambda}{s^2}$.
\end{proposition}

\begin{proof}
	Let $\psi = \sigma^{-1}X^CJ_C=s$ and let $V(s) = \frac{\lambda}{s^2}$. Then $V'(s) = -\frac{2\lambda}{s^3}$. So, by Proposition \ref{Lem:COM}, $r+W$ is a constant of motion if and only if
	\[
	dW = V'(s)\left( \phi d\psi-\tfrac{1}{2}\psi d\phi \right) = \lambda\left( \frac{d\phi}{s^2}-\frac{2\phi ds}{s^{3}} \right) = \lambda d\left( \frac{\phi}{s^2} \right) ,
	\] 
	which is certainly solved by $W = \lambda\frac{\phi}{s^2}$.

\end{proof}

\subsection{Flat Model}
Here we briefly discuss how the tractor bundle arises on the flat model of conformal geometry, and how this naturally gives a basis of tractors on many conformally flat spaces, in particular $\mathbb{R}{}^{P,Q} $. Our presentation follows Section 3 of \cite{curry_introduction_2018}. The model for conformal geometry in signature $(P,Q)$ is constructed by taking a non-degenerate symmetric bilinear form of signature $(P+1,Q+1)$ on $\mathbb{R}^{n+2}$,
\begin{equation*}
	B = \begin{pmatrix} 0&0&1\\0&I&0\\1&0&0 \end{pmatrix}
\end{equation*}
where $I$ is the diagonal matrix with $P$ entries of $+1$ and $Q$ entries of $-1$. The punctured null cone $\mathcal{N}$ of $B$ is the set of $(q_{0},\dots ,q_{n+1})\in \mathbb{R}^{n+2}\setminus \{0\} $ such that
\[
	0=\mathbf{q}^TB\mathbf{q} = \sum_{i=1}^{P} q_i^2- \sum_{i=P+1}^{n} q_i^2 +2q_{0}q_{n+1},
\]
The ray projectivisation $\mathcal{N} / \mathbb{R}_{>0} $ is the model space, which we will call the conformal sphere of signature $(P,Q)$ and denote by $\mathbb{S}{}^{P,Q} $. In Riemannian signature there are two connected components to this model, and usually one connected component is chosen. The conformal geometry of $\mathbb{S}^{P,Q}$ arises by choosing a section of the projectivisation $\mathcal{N}\to \mathbb{S}{}^{P,Q} $ and taking the metric induced by $B$ on this section as an element of the conformal class. In Riemannian signature, a connected component of the model is diffeomorphic to $\mathbb{S}^n$ and conformally equivalent to $(\mathbb{S}^n,\boldsymbol{c})$ where the round sphere metric is in $\boldsymbol{c}$. 

Tractors on the model correspond to tangent vectors on $\mathcal{N}$ which are parallel up the rays. The tractor metric is given by $B$. Parallel cotractors on the model then correspond to one forms on $\mathcal{N}$ which extend to parallel vector fields on $\mathbb{R}{}^{n+2}$. There is an obvious parallel basis of the cotractor bundle, given by
\[
	I = p_0, J = {p_{n+1}}, I^1 = p_1,\dots ,I^n = p_n,
\] 
The tractors $I$ and $J$ are null and have product $1$, and the remaining tractors have length $\pm 1$ and are pairwise orthogonal. The canonical tractor $X$ is given by the Euler vector field on $\mathbb{R}^{n+2}$, $X = \sum_{i=0}^{n+1} q_i\partial_i$. Weight $w$ densities correspond to functions on $\mathcal{N}$ which are homogeneous of degree $w$ up the rays of $\mathcal{N}_+$, so we can see that $X$ has weight one as expected.

If $\sigma$ is a one-density, then the section of $\mathcal{N}\to \mathbb{S}^{P,Q}$ given by points where $\sigma=1$ is gives the pseudo Riemannian space with metric $\sigma^{-2}\boldsymbol{g}$. The length of a scale tractor $D_A\sigma$ is equal to $-\frac{R}{n(n-1)}$, where $R = R{}_{ab}{}^{ab} $ is the scalar curvature of the corresponding metric. So the tractors $I$ and $J$ correspond to flat pseudo-Euclidean space, $\mathbb{R}^{P,Q}$, while the other tractors correspond to spaces of constant positive or negative curvature, all are hyperbolic in the Riemannian case.

So setting $q_{0} = B(X,I) = \sigma$, which is a weight one density, to represent the base Euclidean scale, we find that the zero set of $\sigma$ is the set of rays $(0,\tilde{q},q_{n+1})$, where $\tilde{q}\in \mathbb{R}^{P,Q}$ is null. On the rest of $\mathcal{N}$, the Euclidean space represented by $I$ is given by the geometry of the submanifold $\sigma = 1$, that is, points of the form
\[
	(1,q_{1},\dots ,q_n,q_{n+1}) = \left(1,q_{1},\dots ,q_n,-\frac{r^2}{2}\right),
\] 
where $r^2 = q_{1}^2+\dots+q_P^2 -q_{P+1}^2-\dots -q_n^2$. This shows that 
\[
\tau = B(X,J) = -\frac{r^2}{2}\sigma, \quad \sigma_{1} = B(X,I_1) = q_{1}\sigma,\quad \dots ,\quad  \sigma_n = B(X,I_n) = q_n\sigma,
\] 
where $(q_{1},\dots ,{q}{}_{n})$ are the coordinates on Euclidean space. Hence, in the scale $\sigma$, the system of tractors is
\[
I_A = \begin{pmatrix} \sigma\\ 0 \\ 0 \end{pmatrix}, I^i_A = \begin{pmatrix} \sigma q_i\\  \sigma p_i  \\ 0 \end{pmatrix}, J_A = \begin{pmatrix} -\sigma\frac{r^{2}}{2}\\ -\sigma e \\ \sigma^{-1} \end{pmatrix},
\]
where $e = \sum_{i=1}^{P}  q_ip_i -\sum_{i=P+1}^{n} q_ip_i $.

We begin by focusing on the pair $I,J$. In Theorems \ref{Thm:Ex1} and \ref{Thm:Ex2}, this pair fits into the case \ref{Thm:Ex1:1:IJ}. Theorem \ref{Thm:Ex1} allows potentials of the form 
\begin{equation}\label{Eqn:FM1}
V(q) = \lambda r^2+\lambda_{0},
\end{equation} 
which match the isotropic Harmonic Oscillator and Theorem \ref{Thm:Ex2} allows potentials of the form
\begin{equation}\label{Eqn:FM2}
V(q) = \frac{\lambda}{r}+\lambda_{0},
\end{equation} 
which match the Kepler system. We use the basis $I^1,\dots ,I^n$ of $I^{\perp}\cap J^{\perp}$. Then the Killing vectors $k^{i,j}_b = 2X{}^{A} Z{}^{B}{}_{b} I{}^{i}{}_{[A} I{}^{j}{}_{B]} = \sigma^2(q_ip_j-q_jp_i)$ are the angular momenta. 
By Theorem \ref{Thm:Ex1} \ref{Thm:Ex1:1}, the final constant of motion for \eqref{Eqn:FM1} comes from some ${R}{}_{ABCD} = {(II\owedge U)}{}_{ABCD}$, where $U$ is a symmetric parallel tractor orthogonal to $I$ and $J$. Then the constant of motion is (note that we are using $\mathcal{H} = |p|^2+\lambda r^2 = |p|^2-2\lambda s$ with $s = -r^2 / 2$)
 \[
 u_q(p,p) + {\lambda}\phi(q),
 \] 
 where $u{}_{bd}  = \sigma^2Z{}^{B}{}_{b} Z{}^{D}{}_{d} U{}_{BD} $ and $\phi = \sigma^{-2}X^BX^D U{}_{BD} $.

So the superintegrable system is
\[
	f_1 = \mathcal{H} = |p|^2+\lambda r^2, f_2 = L^{1,2},\dots f_n = L^{1,n},f_{2n-2} = L^{n-1,n}, f_{2n-1} = u_q(p,p) + {\lambda}\phi(q)
\] 
To extend to the family of Hamiltonians,
\[
\mathcal{H}^{\lambda}(q,p) = |p|^2+\lambda r^2 +\frac{\lambda_{1}}{q_{1}^2}+\dots +\frac{\lambda_n}{q_n^2} +\lambda_{0},
\] 
we apply Proposition \ref{Prp:Comp}. We need to shift to the rank two Killing tensors, $r^{i,j}{}_{ab} = k^{i,j}{}_ak^{ij}{}_b$ and the tractors $R^{i,j}{}_{ABCD} = ((I^i)^2\owedge (I^j)^2)$. Clearly $R^{i,j}$ is orthogonal to $I{}^{\ell} $ for $\ell\neq i,j$. Then
\[
	I^j{}^AR{}^{i,j}{}_{ABCD}  = 2I{}^{i}{}_{B} I{}^{j}{}_{[C} I{}^{i}{}_{D]}.
\] 
By Proposition \ref{Prp:Comp}, the constant of motion for $|p|^2+\lambda_j\frac{1}{q_j^2}$ is $L_{i,j}^2+\lambda_j \frac{q_i^2}{q_j^2}$. So the first $2n-2$ constants of motion for $\mathcal{H}^{\lambda}$ are
\[
	f_{1} = \mathcal{H}^{\lambda},f_{2} = (L^{1,2})^2+\lambda_{1}\frac{q_{2}^2}{q_{1}^2}+\lambda_{2}\frac{q_{1}^2}{q_{2}^2},\dots ,f_{2n-2} = (L^{n-1,n})^2+\lambda_{n-1}\frac{q_{n}^2}{q_{n-1}^2}+\lambda_{n}\frac{q_{n-1}^2}{q_{n}^2}
\] 
The final Killing tensor is the top slot of $R{}_{ABCD}  = I^2\owedge U$. In order to ensure compatibility with the other scales, we choose $U = I^1_AI^1_B$. Then we just compute 
\[
	I^1{}^A R{}_{ABCD} = 2I{}_{B} I{}^{1}{}_{[C} I{}_{D]}.
\] 
So the final constant of motion is, again by Proposition \ref{Prp:Comp},
\[
	f_{2n-1} = p_1^2+\lambda q_1^2 + \frac{\lambda_{1}}{q_{1}^2}.
\] 

For a second example, we consider the pair $I^1,I^2$. The base metric for the Hamiltonian in this case has scale density $X{}^{A} I{}^{1}{}_{A}  = q_{1}\sigma$, where $\sigma$ is the Euclidean scale. So, on each of the open half-spaces, the metric is given by $g_1 = q_{1}^{-2}g$, where $g$ is the Euclidean metric. This is the upper half space model of hyperbolic space. In Theorems \ref{Thm:Ex1} and \ref{Thm:Ex2}, this pair fits into the case \ref{Thm:Ex1:2:IJ}. For the following calculations take \(q_1>0\). Theorem \ref{Thm:Ex1} allows potentials of the form 
\begin{equation}\label{Eqn:FM3}
V(q) = \lambda \frac{q_{1}^2}{q_{2}^2} +\lambda_{0},
\end{equation} 

Theorem \ref{Thm:Ex2} allows potentials of the form
\begin{equation}\label{Eqn:FM4}
	V(q) = \frac{\lambda s}{\sqrt{s^2 +1}}+\lambda_{0}= \frac{\lambda q_{2}}{\sqrt{q_{1}^2+q_{2}^2} }+\lambda_{0},
\end{equation} 
where $s = q_{2} / q_{1}$. We will compute the details for \eqref{Eqn:FM4}. The space of tractors orthogonal to $I{}^{1} $ and $I^{2}$ is spanned by $I,J,I^3,\dots ,I{}^{n} $, but this basis is not pairwise orthogonal. Instead it is natural to take $I+J,I-J,I^{3},\dots I^n$. This means that the constants of motion we construct will be compatible with the spherical metric, represented by $I-J$, and an additional hyperbolic metric represented by $I+J$, but not either Euclidean metric. But, as we show later, it can be adapted to be compatible with $I$. 

The Killing vectors for this Hamiltonian are $k^{i,j}_b = 2X{}^{A} Z{}^{B}{}_{b} I{}^{i}{}_{[A} I{}^{j}{}_{B]} = \sigma^2(q_ip_j-q_jp_i)$ for $3\le i<j\le n $ and 
\[
k^{\pm i}{}_{b}  = 2X{}^{A} Z{}^{B}{}_{b} (I{}^{i}{}_{[A}(I\pm J){}_{B]} ) = \sigma^2(\mp eq_i-p_i\pm \tfrac{r^2}{2}p_i),\] 
where $e = \sum_{i=1}^{n} q_ip_i$. 
For the final Killing tensor, we can take $R{}_{ABCD} = (\pi\owedge S)_{ABCD}$, where $S{}_{AB} = 2I{}^{2}{}_{(A}V{}_{B)}  $ for any parallel $V$ orthogonal to $I^1$ and $I{}^{2}$, where $\pi{}_{AB} = h{}_{AB} -I{}^{1}{}_{A} I{}^{1}{}_{B} -I{}^{2}{}_{A} I{}^{2}{}_{B} $. We then choose $V = I^3$, so 
\[
	 X{}^{A} Z{}^{B}{}_{b} X{}^{C} Z{}^{D}{}_{d}R{}_{ABCD}   =  \sigma^{4}\left(q_{2}q_{3}|p|^2-L^{1,2}L^{1,3} \right)\quad \text{and} \quad
\]
\[
	I{}^{2A} R{}_{ABCD} = \pi{}_{B[D} I{}^{3}{}_{C]}.
\] 

So, by Theorem \ref{Thm:Ex2}, the second-order maximally superintegrable Hamiltonian is
\begin{equation}\label{Eqn:12Kep}
\mathcal{H}(q,p) = q_{1}^2|p|^2+\frac{q_{2}}{\sqrt{q_{1}^2+q_{2}^2}}
\end{equation} 
and the constants of motion are, when $n\ge 4$,
\[
	f_{1} = \mathcal{H}, f_{2} = k^{+3}, f_{3} = k^{-3}, f_{4} = L^{3,4},\dots ,f_n = L^{3,n}, f_{n+1} = k^{+n}, f_{n+2} = k^{-n},
\]\[f_{n+3} = L^{4,n},\dots ,f_{2n-2} = L^{n-1,n}\quad \text{and} \quad f_{2n-1} = q_{2}q_{3}|p|^2-L{}^{1,2} L{}^{1,3}-\frac{q_{3}}{2\sqrt{q_{1}^2+q_{2}^2} }.
\] When $n = 3$, we can take $k^{-3}$, $k^{+3}$ and $e = \sum q_ip_i $ together with the quadratic integral.

To construct a second-order maximally superintegrable linear family, it would be natural to choose the orthonormal basis $I+J,I-J,I^3,\dots ,I^n$ for $I^1{}^{\perp}\cap I^2{}^{\perp}$, which gives
\[
	\mathcal{H}^{\lambda}(q,p) = q_{1}^2|p|^2+\lambda_{1} \frac{q_{2}}{\sqrt{q_{1}^2+q_{2}^2}}+\lambda_{2}\frac{q_{1}^2}{(1- \tfrac{r^2}{2})^2}+\lambda_{3}\frac{q_{1}^2}{q_{3}^2}+\dots +\lambda_n\frac{q_{1}^2}{q_n^2}+\lambda_{0},
\] 
or another choice of $n-1$ scales from $q_{3},\dots ,q_n,1-\frac{r^2}{2},1+\frac{r^2}{2}$ and this is second-order maximally superintegrable by Theorem \ref{Thm:MSFam}.

However, we can also take $I,I^3,\dots ,I^n$ in Theorem \ref{Thm:MSFam}, giving the family of Hamiltonians
\begin{equation}\label{Eqn:H21}
	\mathcal{H}^{\lambda}(q,p) = q_{1}^2|p|^2+\lambda_{1} \frac{q_{2}}{\sqrt{q_{1}^2+q_{2}^2}}+\lambda_{2}q_{1}^2+\lambda_{3}\frac{q_{1}^2}{q_{3}^2}+\dots +\lambda_n\frac{q_{1}^2}{q_n^2}+\lambda_{0},
\end{equation}
which after a St\"ackel transform, dividing by $q_{1}^2$ gives
\[
	\tilde{\mathcal{H}}^{\lambda}(q,p) = |p|^2+\lambda_{1} \frac{q_{2}}{q_{1}^2\sqrt{q_{1}^2+q_{2}^2}}+\lambda_{0}\frac{1}{q_{1}^2}+\lambda_{3}\frac{1}{q_{3}^2}+\dots +\lambda_n\frac{1}{q_n^2}+\lambda_{2},
\]
which provides a generalisation to any dimension of the third second-order maximally superintegrable system in Table I of \cite{evans_superintegrability_1990}. 

We briefly outline the calculation of the constants of motion for \eqref{Eqn:H21}. From the proof of Theorem \ref{Thm:MSFam}, the constants of motion will have quadratic parts induced by various Killing tensors, note that the labels on the tractors there differ from those here. 

First, $r^{I,3},\dots ,r^{I,n},r^{I,J}$, induced by tractors $4I{}_{[A} I{}^{i}{}_{B]} I{}_{[C} I{}^{i}{}_{D]} $ and $4I{}_{[A} J{}_{B]} I{}_{[C}J{}_{D]} $, so the squares of linear momenta, $p_{3}^2,\dots ,p_n^2$ and the square of the Euler vector field $(\sum_{j=1}^{n} p_jq_j)^2$. Secondly, Killing tensors that are the top slots of
\[
	R{}^{i}{}_{ABCD}  = 2I{}_{[A} I{}^{i}{}_{B]} J{}_{[C} I{}^{i}{}_{D]} +2J{}_{[A} I{}^{i}{}_{B]} I{}_{[C} I{}^{i}{}_{D]},
\]
which are $p_i\left(q_ie-\frac{r^2}{2}p_i\right)$. To find the scalar terms, we first calculate
\[
	4I{}^{A} I{}_{[A} J{}_{B]} I{}_{[C} J{}_{D]} = 2I{}_{B} J{}_{[C} I{}_{D]}  \quad \text{and} \quad 4I{}^{i}{}^{A} I{}_{[A} I{}^{i}{}_{B]} I{}_{[C} I{}^{i}{}_{D]} = 2I{}_{B} I{}^{i}{}_{[C} I{}_{D]}
\] 
So, by Proposition \ref{Prp:Comp}, the first $n$ constants of motion for $\mathcal{H}^{\lambda}$ are
\[
f_{1}=\mathcal{H}^{\lambda}, f_{2} = \left(\sum_{j=1}^{n} q_jp_j\right)^2+\lambda_{2}r^2, f_{3} = p_{3}^2+\frac{\lambda_{3}}{q_{3}^2},\dots ,f_n = p_n^2+\frac{\lambda_n}{q_n^2}
\]
Then we calculate
\[
	I{}^{A} R{}^{i}{}_{ABCD} =I{}^{i}{}_{B} I{}_{[C} I{}^{i}{}_{D]} \quad \text{and} \quad I{}^{i}{}^{A} R{}^{i}{}_{ABCD} = I{}_{B} I{}^{i}{}_{[C} J{}_{D]} +J{}_{B} I{}^{i}{}_{[C} I{}_{D]} 
\] 
so the following $n-2$ constants of motion are
\[
	f_{n-2+i} = p_{i}\left( q_{i}e - \frac{r^2}{2}p_{i} \right)+\lambda_{2}\frac{q_{i}^2}{2}-\lambda_{i}\frac{r^2}{2q_{i}^2},
\] 
for $i\in \{3,\dots ,n\} $. The final Killing tensor is the top slot of
\[
	R{}_{ABCD} = 2I{}_{[A} \pi{}_{B][D} I{}^{2}{}_{C]} +2I{}^{2}{}_{[A} \pi{}_{B][D} I{}_{C]},
\]
which is $\sigma^{4}(q_{2}\boldsymbol{g}{}_{bd} + p{}^{1}{}_{(b} L^{1,2}{}_{d)})$. We calculate
\[
	I{}^{A} R{}_{ABCD} = I{}_{B} I{}_{[C} I{}^{2}{}_{D]} \quad \text{and} \quad I{}^{i}{}^{A} R{}_{ABCD}  = I{}_{B} I{}^{i}{}_{[C} I{}^{2}{}_{D]} +I{}^{2}{}_{B} I{}^{i}{}_{[C} I{}_{D]} 
\]
for $i\in \{3,\dots ,n\} $ and
\[
	I{}^{2A} R{}_{ABCD} = \pi{}_{B[D} I{}_{C]},
\] 

so the final constant of motion as in \eqref{Eqn:Ex2:2:CM} and using Proposition \ref{Prp:Comp} is 
\[
	f_{2n-1} = q_{2}|p|^2+p_{1}L{}^{1,2} +\lambda_{2}q_{2}+\lambda_{3}\frac{q_{2}}{q_{3}^2}+\dots +\lambda_n \frac{q_{2}}{q_n^2}-\frac{\lambda_{1}}{2\sqrt{q_{1}^2+q_{2}^2} }.
\] 

\newpage
\bibliographystyle{amsplain}
\bibliography{../DGThesis.bib}
\end{document}